\documentclass[12pt]{article}

\usepackage[english]{babel}

\usepackage{amsfonts,amsmath,amsxtra,amsthm,amssymb,latexsym}
\usepackage{color}

\newcommand{\ack}{\section*{Acknowledgments}}
\newtheorem{theorem}{Theorem}[section]
\newtheorem{corollary}[theorem]{Corollary}
\newtheorem{lemma}[theorem]{Lemma}

\newtheorem{definition}[theorem]{Definition}
\newtheorem{remark}[theorem]{Remark}
\newtheorem{example}[theorem]{Example}
\newtheorem{proposition}[theorem]{Proposition}

\numberwithin{equation}{section}

\DeclareMathOperator{\supp}{supp} \DeclareMathOperator{\im}{Im}
\DeclareMathOperator{\re}{Re} 
 \DeclareMathOperator{\dom}{dom}
\DeclareMathOperator{\ran}{ran} \DeclareMathOperator{\Ext}{Ext}
\DeclareMathOperator{\Span}{span}

\DeclareMathOperator{\Rank}{rank}
\DeclareMathOperator{\diag}{diag}\DeclareMathOperator{\loc}{loc}
\DeclareMathOperator{\comp}{comp}\DeclareMathOperator{\per}{per}
\DeclareMathOperator{\op}{op}

\newcommand\R{{\mathbb{R}}}
\newcommand\C{{\mathbb{C}}}

\newcommand\N{{\mathbb{N}}}
\newcommand\Z{{\mathbb{Z}}}

\newcommand\gH{{\mathfrak{H}}}
\newcommand{\gG}{{\Gamma}}
\newcommand{\gT}{{\Theta}}

\newcommand{\gd}{{d}}
\newcommand{\gA}{{\alpha}}
\newcommand{\gB}{{\beta}}

\newcommand\cH{{\mathcal{H}}}
\newcommand\cI{{\mathcal{I}}}

\newcommand\cC{{\mathcal{C}}}
\newcommand\cN{{\mathfrak{N}}}

\newcommand\rH{{\rm{H}}}
\newcommand\rF{{\rm{F}}}

\newcommand\rD{{\rm{d}}}
\newcommand\I{{\rm{i}}}
\newcommand\gr{{\rm{gr}}}

\def\Ext{{\rm Ext\,}}

\def\mul{{\rm mul\,}}
\def\wt#1{{{\widetilde #1} }}

\title{1--D Schr\"odinger operators with local interactions\\ on a discrete set}

 \author{Aleksey Kostenko and Mark Malamud}

\date{}

\begin{document}
\maketitle
\begin{abstract}
Spectral properties of 1-D Schr\"odinger operators
$\mathrm{H}_{X,\alpha}:=-\frac{\mathrm{d}^2}{\mathrm{d} x^2} + \sum_{x_{n}\in
X}\alpha_n\delta(x-x_n)$ with local point interactions on a discrete
set $X=\{x_n\}_{n=1}^\infty$ are well studied
when $d_*:=\inf_{n,k\in\N}|x_n-x_k|>0$. Our paper is devoted to the case
$d_*=0$. We consider  $\mathrm{H}_{X,\alpha}$ in the
framework of extension theory of symmetric operators
by applying the technique of boundary triplets and the
corresponding Weyl functions.

We  show that the spectral properties of
$\mathrm{H}_{X,\alpha}$  like self-adjointness, discreteness, and lower
semiboundedness correlate with the corresponding spectral properties of certain classes  of Jacobi
matrices. Based on this connection, we obtain necessary and
sufficient conditions for the operators $\mathrm{H}_{X,\alpha}$ to be
self-adjoint, lower-semibounded, and discrete in the case
$d_*=0$.

The operators with $\delta'$-type interactions are investigated
too.
The obtained results demonstrate that in the case $d_*=0$, as
distinguished from the case $d_*>0$, the spectral properties of
the operators with $\delta$ and $\delta'$-type interactions are
substantially different.
\end{abstract}

\quad

\noindent {\bf Keywords: } Schr\"odinger operator, local point interaction, self-adjointness, lower semiboundedness, discreteness

\quad

 \noindent {\bf AMS Subject classification: } 34L05, 34L40, 47E05, 47B25, 47B36, 81Q10

\newpage

\tableofcontents

\newpage

\section{ Introduction}\label{intro}

Differential operators with point interactions arise in various physical applications as exactly solvable models that describe complicated physical phenomena (numerous results as well as a comprehensive list of references may be found in  \cite{Alb_Ges_88, Alb_Kur_00, Exn_04}). An important class of such operators is formed by the differential operators with the coefficients having singular support on a disjoint set of points. The most known examples are the operators $\rH_{X,\gA,q}$ and $\rH_{X,\gB,q}$ associated with the formal differential expressions
\begin{equation}\label{I_01}
\ell_{X,\gA,q}:=-\frac{\rD^2}{\rD x^2}+q(x)+\sum_{x_{n}\in X}\gA_n\delta_n,\qquad
 \ell_{X,\gB,q}:=-\frac{\rD^2}{\rD x^2}+q(x)+\sum_{x_{n}\in X}\gB_n(\cdot,\delta'_n)\delta'_n, 
\end{equation}
where $\delta_n:=\delta(x-x_n)$ and $\delta$ is  a Dirac delta-function.
These operators describe $\delta$- and $\delta'$-interactions, respectively, on a discrete set $X=\{x_n\}_{n\in I}\subset\R$, and the coefficients $\gA_n,\ \gB_n\in\R$ are called the strengths of the interaction at the point $x=x_n$. Investigation of these models was originated by Kronig and Penney \cite{Kro_Pen} and Grossmann et. al. \cite{Gro} (see also \cite{Ges_Hol_87}), respectively. In particular, the "Kronig--Penney model" ($\ell_{X,\gA,q}$ with $X=\Z$, $\gA_n\equiv \gA$, and $q\equiv 0$) provides a simple model for a nonrelativistic electron moving in a fixed crystal lattice.

There are several ways to associate the operators with $\ell_{X,\gA,q}$ and $\ell_{X,\gB,q}$. For example, a $\delta$-interaction at a point $x=x_0$ may be defined using the \emph{form method}, that is the operator $-\frac{\rD^2}{\rD x^2}+\gA_0\delta(x-x_0)$
is defined as an operator associated in $L^2(\R)$ with the quadratic form 
\[
\mathfrak{t}[f]=\int_{\R}|f'(t)|^2dt+\gA_0|f(x_0)|^2,\qquad f\in W_2^1(\R).
\]
Another way to introduce a local interaction at $x_0$ is to consider a symmetric operator $\rH_{\min}:=\rH_{\min}^-\oplus \rH_{\min}^+$, where $\rH_{\min}^-$ and $\rH_{\min}^+$ are the minimal operators generated by $-\frac{\rD^2}{\rD x^2}$ in $L^2(-\infty,x_0)$ and $L^2(x_0,+\infty)$, respectively, and to impose \emph{boundary conditions} connecting $x_0+$ and $x_0-$.

Both these methods have disadvantages if the set $X$ is infinite. The form method works only for the case of lower semibounded operators. If we apply the method of boundary conditions, then the corresponding minimal operator $\rH_{\min}$ has infinite deficiency indices and the description of self-adjoint extensions of $\rH_{\min}$ is rather complicated problem in this case.

An alternative approach was proposed recently in \cite{bsw} (see also \cite{Min_86} for the case of $\delta$-type interactions). Namely, the operators with general local interactions on a discrete set $X$ were defined as self-adjoint extensions such that the Lagrange brackets $[f,g]:=\overline{f(x)}g'(x)-\overline{f'(x)}g(x)$ are continuous on $\R$ for arbitrary elements $f,g$ from the domain.  
It was shown in \cite{bsw, Min_86} that classical Sturm--Liouville theory with all its fundamental objects can be generalized to include local point interactions. In particular, the Weyl's alternative has been established in this case.


Nevertheless, to the best of our knowledge there are only a few results that describe the spectral properties of operators with local interactions in the case $\gd_*=0$, where
\begin{equation}\label{d_*}
\gd_*:=\inf_{i,j\in I}|x_{i}-x_j|=0.
\end{equation}

Let us present a brief historical overview. Note that we are interested in the case when the set $X$ is infinite (the case $|X|<\infty$ is considered in great detail in \cite{Alb_Ges_88}).
First we need some notation. Let $\cI$ be the semi axis, $\cI=[0,+\infty)$, and let
$X=\{x_n\}_{n=1}^\infty\subset \cI$ be a strictly increasing
sequence, $x_{n+1}>x_{n},\ n\in \N$, such that $x_n\to +\infty$.
We denote $\gd_n:=x_{n}-x_{n-1}$, $x_0:=0$, and assume $q\in L^2_{loc}[0,+\infty)$.  
In $L^2(\cI)$, the minimal symmetric  operators $\rH_{X,\gA,q}$  and
$\rH_{X,\gB,q}$ are naturally associated  with \eqref{I_01}. Namely, define the operators $\rH^0_{X,\gA,q}$  and
$\rH^0_{X,\gB,q}$ by the differential expression
\begin{equation}
\tau_q:=-\frac{\rD^2}{\rD x^2}+q(x)\label{I_03}
\end{equation}
on the domains, respectively,
\begin{gather}
\dom(\rH^0_{X,\gA,q})=\left\{f\in W^{2,2}_{\comp}(\cI\setminus X): f'(0)=0,\ \begin{array}{c}  f(x_n+)=f(x_n-)\\
 f'(x_n+)-f'(x_n-)=\gA_n f(x_n)\end{array},\ n\in\N \right\} ,\label{I_04}\\
\dom(\rH^0_{X,\gB,q})=\left\{f\in W^{2,2}_{\comp}(\cI\setminus X): f'(0)=0,\ \begin{array}{c}  f'(x_n+)=f'(x_n-)\\
 f(x_n+)-f(x_n-)=\gB_n f'(x_n)\end{array},\ n\in\N \right\}.\label{I_05}
\end{gather}
Let $\rH_{X,\gA,q}$ and $\rH_{X,\gB,q}$ be the closure of
$\rH^0_{X,\gA,q}$ and $\rH^0_{X,\gB,q}$, respectively. In general,
the operators $\rH_{X,\gA,q}$ and $\rH_{X,\gB,q}$ are symmetric but not automatically self-adjoint, even in the case $q\equiv 0$.

Spectral analysis of 
$\rH_{X,\gA,q}$ and $\rH_{X,\gB,q}$ consists (at least partially) of the following problems:
      \begin{description}
\item $(a)$ Finding self-adjointness criteria for $\rH_{X,\gA,q}$
and $\rH_{X,\gB,q}$ and description of self-adjoint extensions  if
the deficiency indices $\rH_{X,\gA,q}$ and $\rH_{X,\gB,q}$ are
nontrivial.
\item $(b)$ Lower semiboundedness of the operators
$\rH_{X,\gA,q}$ and $\rH_{X,\gB,q}$.
\item $(c)$ Discreteness of
the spectra  of the operators  $\rH_{X,\gA,q}$ and
$\rH_{X,\gB,q}$.
\item $(d)$ Characterization of continuous,
absolutely continuous, and singular parts of the spectra of the
operators $\rH_{X,\gA,q}$ and $\rH_{X,\gB,q}$.
\item $(e)$ Resolvent comparability of the operators $\rH_{X,\gA^{(1)},q}$ and
$\rH_{X,\gA^{(2)},q}$ with  $\gA^{(1)}\neq \gA^{(2)}$.
      \end{description}

 In the present paper, we confine ourselves to the case of bounded potentials $q\in L^\infty(\cI)$.
 Let us note that the case of unbounded $q$ was studied in \cite{Bra85, bsw, Ges_Kir_85, Chr_Sto_94}
 and  the case of $q$ being a $W^{2,-1}_{\loc}(\cI)$ distribution
 was studied in  \cite{HryMik_01, HryMik_02, SavShk_99, SavShk_03} (see also the references therein).
More precisely, it is shown in \cite{bsw} (see also \cite{Min_86}) that $\mathrm{n}_\pm(\rH_{X,\gA,q})\leq 1$
and the deficiency indices may be characterized in terms of the limit point and the limit circle classification
for the endpoint $x=+\infty$. Brasche  \cite[Theorem 1]{Bra85} proved that $\rH_{X,\gA,q}$ is self-adjoint and
lower semibounded if the potential $q$ is lower semibounded and the strengths $\gA_n,\ n\in \N,$ are nonnegative.
Assuming the condition $\gd_*>0$, Gesztesy and Kirsch \cite{Ges_Kir_85}, Christ and Stolz \cite{Chr_Sto_94}
(see also \cite{bsw}) established self-adjointness of $\rH_{X,\gA,q}$ for several classes of unbounded potentials $q$.
In particular, Gesztesy and Kirsch \cite[Theorem 3.1]{Ges_Kir_85} proved that $\rH_{X,\gA,q}=\rH_{X,\gA,q}^*$
if $q\in L^\infty(\cI)$ and $\gd_*>0$ (other proofs are given in \cite{Koc_79} and \cite{Chr_Sto_94}).
Moreover, Christ and Stolz \cite[pp. 495--496]{Chr_Sto_94} showed that the condition $\gd_*>0$ cannot
be dropped there even if $q\equiv 0$. More precisely, they proved that
$\mathrm{n}_\pm(\rH_{X,\gA,0})=1$ if $d_n=\frac{1}{n}$ and  $\alpha_n=-2n-1$, $n\in \N$.
Note also that self-adjointness of $\rH_{X,\gA,0}$ with arbitrary $X=\{x_n\}_{n=1}^\infty\subset\cI$
was erroneously stated without proof in \cite{Mih_93}.

Finally, we emphasize that 
in contrast to $\delta$-type interactions the operator $\rH_{X,\gB,0}$ is self-adjoint for
arbitrary $\{\gB_n\}_{n=1}^\infty\subset\R$ 
(see \cite[Theorem 4.7]{bsw}). Let us also mention the recent papers \cite{AlbNiz, Niz_03} dealing with spectral properties of Hamiltonians with $\delta'$-interactions on compact subsets of $\R$ with Lebesgue measure zero.

In the present paper, we investigate problems $(a)-(c)$ and $(e)$
in the case  $d_*=0$ 
and $q\in L^{\infty}(\R)$ (we postpone the study of the case of unbounded
$q$ as well as the problem $(d)$ to our forthcoming paper).
We consider the operators with point interactions in the framework
of extension theory of symmetric operators. This approach allows
one to treat the operators $\rH_{X,\gA,q}$  and $\rH_{X,\gB,q}$ as
self-adjoint (or symmetric) extensions of the minimal operator
\begin{equation}\label{I_06}
\rH_{\min} :=\oplus_{n\in \N}\rH_{n},\qquad
\rH_n=-\frac{\rD^2}{\rD x^2}+q(x),\qquad
\dom(\rH_{n})=W^{2,2}_0[x_{n-1},x_n],
\end{equation}
being  a direct sum of symmetric operators $\rH_n$ with deficiency
indices $n_{\pm}(\rH_n)=2$.

We investigate these operators by applying the technique of
boundary triplets and the corresponding Weyl functions (see
Section \ref{Sec_II_Prelim} for precise definitions). This new
approach to extension theory of symmetric operators has been
appeared and elaborated  during the last three decades
(see \cite{Gor84, DM91, DM95, BGP07} and references therein).
The main ingredient 
is an abstract version of the
Green formula for the adjoint $A^*$ of a symmetric operator $A$ (see formula \eqref{II.1.2_green_f}). A boundary triplet for $A^*$ always exists whenever $n_+(A)= n_-(A)$, though it is not unique. Its
role in extension theory is similar to  that of a coordinate system in analytic
geometry. It enables one to describe self-adjoint extensions
in terms of (abstract) boundary conditions in place of the second
J. von Neumann formula, though this description is simple and
adequate only for a suitable choice of a boundary triplet. Note that construction of a suitable
boundary triplet is a rather difficult problem if $n_\pm(A)=\infty$.

This approach was first applied to the spectral analysis of
$\rH_{X,\gA,q}$ by Kochubei in \cite{Koc_89}.
More precisely,
he proved that in the case $\gd_*>0$ (and $q\in L^\infty(\cI)$) a boundary triplet
$\Pi$ for $\rH_{\min}^*$ can be chosen as a direct sum of triplets $\Pi_n$ defined by \eqref{IV.1.1_05}, that is
$\Pi:=\{\cH,\Gamma_0,\Gamma_1\}:=\oplus_{n=1}^\infty \Pi_n$, where
\begin{equation}\label{I_07}
\cH :=\oplus_{n\in\N}\cH_n,\qquad \Gamma_0
:=\oplus_{n\in\N}\Gamma_0^{(n)},\qquad \Gamma_1
:=\oplus_{n\in\N}\Gamma_1^{(n)}.
\end{equation}
Based on this construction, he gave an alternative proof of the self-adjointness of $\rH_{X,\gA,0}$
(see \cite[Theorem 3.1]{Ges_Kir_85}) and
investigated the problem $(e)$ as well.

The main difficulty in extending this approach  to the case $\gd_*=0$
(or unbounded $q$) is the  construction of a suitable boundary
triplet for the operator $\rH_{\min}^*$ (see \cite{Koc_79,
Koc_89}).
It looks natural that the triplet
$\Pi=\{\cH,\Gamma_0,\Gamma_1\}$ defined by \eqref{I_07} and \eqref{IV.1.1_05}  forms
a boundary triplet for $\rH_{\min}^*$ in this case too.  Indeed,
Green's identity holds for $f,g\in \dom(\rH_{\min}^*)$ with compact supports in $\cI$. However,
$\dom(\Gamma_0)\cap\dom(\Gamma_1)$ 
is only a
proper part of $\dom(\rH_{\min}^*)$ and the boundary mapping
$\Gamma:=\{\Gamma_0,\Gamma_1\}$ cannot be extended onto $\dom(\rH_{\min}^*)$ if $\gd_*=0$.
In this case,  erroneous  construction of
a boundary triplet for $\rH_{\min}^*$ was
announced in  \cite{Mih_93}
(see Remark \ref{rem_Mih}). Note also that the first example the operator \eqref{I_06} with $q\notin L^\infty$ and such that $\Pi$ is not a boundary triplet for $\rH_{\min}^*$ was given in \cite{Koc_79}.

Recently Neidhardt and one of the authors  proved that the
triplet of the form \eqref{I_07} becomes a boundary triplet
after appropriate regularization of the mappings $\Gamma^{(n)}_0$
and $\Gamma^{(n)}_1$, $n\in\N$ (see \cite[Theorem 5.3]{MalNei_08}).
Starting with this result, we investigate the problem in full generality. More precisely, we show that in general
$\Pi=\{\cH,\Gamma_0,\Gamma_1\}$ of the form \eqref{I_07} is only a boundary relation in the sense of \cite{DHMS06} and we find
a criterion for $\Pi$ to form a boundary
triplet for $\rH_{\min}^*$.  Moreover, we present a  general
regularization procedure that enables  us to construct a suitable
boundary triplet $\Pi$ for $\rH_{\min}^*$ in the form
$\Pi=\oplus_{n=1}^\infty \Pi_n$. Namely, in this boundary triplet
the sets of Hamiltonians $\rH_{X,\gA,0}$ and $\rH_{X,\gB,0}$ are
parameterized \emph{by means of certain classes} of Jacobi  (\emph{tri-diagonal})
matrices (the construction from \cite{MalNei_08}
leads to \emph{multi-diagonal} matrices). In turn, the latter leads to a
\emph{correlation} between spectral properties of the Hamiltonians
\eqref{I_01} and the corresponding Jacobi matrices. Note that
another technique for analyzing spectral properties of
$\rH_{X,\gA,0}$ and $\rH_{X,\gB,0}$ by means of second order
difference operators was proposed by Phariseau \cite{Pha_60} (see
also \cite[Chapter III.2.1]{Alb_Ges_88}).

More precisely, in the case of $\delta$-interactions, we show that
the spectral properties of the operator $\rH_{X,\gA,0}$ are closely
connected with the corresponding spectral properties of the  Jacobi matrix
\begin{gather}
B_{X,\gA}=\left(\begin{array}{cccc}
r_1^{-2}\bigl(\alpha_1+\frac{1}{\gd_1}+\frac{1}{\gd_2}\bigr) & (r_1r_2\gd_2)^{-1} & 0&   \dots\\
(r_1r_2\gd_2)^{-1} &r_2^{-2}\bigl(\alpha_2+\frac{1}{\gd_2}+\frac{1}{\gd_3}\bigr) & (r_2r_3\gd_3)^{-1} &  \dots\\
0 & (r_2r_3\gd_3)^{-1} & r_3^{-2}\bigl(\alpha_3+\frac{1}{\gd_3}+\frac{1}{\gd_4}\bigr)&  \dots\\
\dots & \dots & \dots& \dots
\end{array}\right),\label{I_08}
   \end{gather}
   where $r_n=\sqrt{\gd_n+\gd_{n+1}}, \ n\in\N$.
We first show that $\mathrm{n}_\pm(\rH_{X,\gA,0})=\mathrm{n}_\pm(B_{X,\gA})$ (Theorem
\ref{th_delta_sa}) and hence that $\mathrm{n}_\pm(\rH_{X,\gA,0})\leq 1$ (cf. \cite{Min_86, bsw}). Combining  this with the Carleman criterion, we arrive at the following result (see Proposition
\ref{cor_delta_carleman}):

\emph{the operator $\rH_{X,\gA,q}$ with $\delta$-interactions is
self-adjoint for any $\gA=\{\gA_n\}_{n\in\N}\subset \R$ provided
that}
\[
\sum_{n\in\N}\gd_n^2=\infty \quad \text{and}\quad q\in L^\infty(\cI).
\]
This result is sharp. Namely, (see Proposition
\ref{prop_IV.2.2_03}):

\emph{if  $\sum_{n\in\N}\gd_n^2<\infty$ and $X=\{x_n\}_{n\in\N}$  satisfies also
some concave assumptions, then there exists
$\gA=\{\gA_n\}_{n\in\N}$ such that the operator $\rH_{X,\gA,0}$ is
symmetric with $\mathrm{n}_\pm(\rH_{X,\gA,0})=1$}.

Moreover, we show that the equality
$\mathrm{n}_\pm(\rH_{X,\gA,0})=1$ yields that the strengths $\gA_n$ cannot
tend to $\infty$ very fast (Proposition \ref{prop_IV.2.2_04}).
This situation is illustrated by Example \ref{example_IV.2.2_01}.
More precisely,  let $\rH_{X,\gA,0}$ be the minimal
closed symmetric operator associated with the differential
expression $\ell_{X,\gA,0}$, where $\cI=\R_+$ and $X=\{x_n\}_{n\in\N}$ is defined by $\gd_n=x_n-x_{n-1}:=\frac{1}{n},\ n\in\N$. Then
\begin{description}
\item $(i)$ \ $\mathrm{n}_\pm(\rH_{X,\gA,0})=0$\ \ \emph{if
either}\ \ $\alpha_n\leq -(4n+2)+O(n^{-1})$\ \  \emph{or}\ \
$\alpha_n\geq -Cn^{-1}$ with some $C>0$, \item $(ii)$ \
$\mathrm{n}_\pm(\rH_{X,\gA,0})=1$\ \  \emph{if}\ \  $\alpha_n= -
a(4n+2)+O(n^{-1})$ \emph{with} $a\in(0,1)$.
\end{description}
The latter enables us to construct a positive potential $q>0$ (see Section \ref{Sec_V}) such that the operator
$\rH_{X,\gA,q}$ with $\gA_n=-4n-2$ and $\gd_n=x_n-x_{n-1}=1/n$ 
is symmetric with $n_\pm(\rH_{X,\gA,q})=1$. This shows that self-adjointness of $\rH_{X,\gA,0}$ is not stable under positive perturbations in the case $\gd_*=0$ (in the case $\gd_*>0$, it was shown in \cite[Theorem 3.1]{Ges_Kir_85} that self-adjointness of $\rH_{X,\gA,0}$ is stable under perturbations by a wide class of potentials $q$).

Further, in the case $\gd_*=0$ we solve the problems $(b)$ and $(c)$ in terms of the Jacobi operators \eqref{I_08}. Namely, we show that \emph{the operator $\rH_{X,\gA,0}$ is lower semibounded if and only if the operator
$B_{X,\gA}$ is also lower semibounded}.
As for discreteness of the spectrum of $\rH_{X,\gA,0}$, we first note that any self-adjoint extension of $\rH_{X,\gA,0}$ has discrete spectrum whenever $\mathrm{n}_{\pm}(\rH_{X,\gA,0}) = 1$. In the case $\rH_{X,\gA,0}=\rH_{X,\gA,0}^*$, \emph{the operator
$\rH_{X,\gA,0}$  has discrete spectrum if and only if $\gd_n\to 0$
and $B_{X,\gA}$ is discrete} (Theorem \ref{th_disc_d}).

Using recent advances in the spectral theory of unbounded Jacobi
operators (see \cite{JanNab01, JanNab03, CojJan07}), we obtain
necessary and sufficient conditions for discreteness and lower
semiboundedness of the operator $\rH_{X,\gA,0}$ in the case $\gd_*=0$.
We show that condition
\begin{equation}\label{I_11}
\frac{\gA_n}{\gd_n+\gd_{n+1}}\geq
C,\quad n\in \N,\quad\text{ for some}\quad C\in \R,
\end{equation}
is sufficient for semiboundedness. If $\gd_*>0$, then \eqref{I_11}
reads $\inf_{n\in\N}\gA_n>-\infty$ and it is also  necessary (see
\cite{Bra85} and also Corollary \ref{cor_IV.2.5_01}). If
$\gd_*=0$, then the situation becomes more complicated. In
Proposition \ref{prop_IV.2.5_01}, we show that the operator
$\rH_{X,\gA,0}$ might be non-semibounded even if $\gA_n\to 0$.

Further (see Proposition
\ref{prop_chihara_1}), \emph{the operator $\rH_{X,\gA,0}=
\rH_{X,\gA,0}^*$ is discrete provided that }   
\begin{equation}\label{I_10}
\lim_{n\to\infty}\gd_n=0,\qquad \lim_{n\to\infty}\frac{|\gA_n|}{\gd_n}=\infty,\quad \text{\emph{and}}\quad \lim_{n\to\infty}\frac{1}{\gd_n\gA_n}>-\frac14.
\end{equation}
The third condition in \eqref{I_10} is sharp (cf. Remark \ref{rem_IV.2.5_01}). Besides, \eqref{I_10} implies that
$\rH_{X,\gA,0}$ may be discrete if $\gA=\{\gA_n\}_{n\in\N}$ is bounded. 
Also \eqref{I_10} enables us to construct operators, 
which are discrete but not lower semibounded. For instance, \emph{the operator}
$\rH_C=-\frac{\rD^2}{\rD
x^2}-\sum_{n\in\N}C\sqrt{n}\ \delta(x-\sqrt{n})$
\emph{with $C>8$ has discrete spectrum though it is not lower semibounded}.

Let us stress that the spectral properties of the operators $\rH_{X,\gA,0}$ and $\rH_{X,\gB,0}$ are completely different in the case $\gd_*=0$. This becomes clear because of the structure of the boundary operators $B_{X,\gA}$ and $B_{X,\gB}$ that parameterize the Hamiltonians $\rH_{X,\gA,0}$ and $\rH_{X,\gB,0}$, respectively. Namely, we show that the spectral properties of the operator with $\delta'$-interactions are closely connected with the Jacobi matrix
\begin{equation}\label{I_12}
B_{X,\gB}:=R_X^{-1/2}(I+U^*)B_\gB^{-1}(I+U)R_X^{-1/2},\qquad B_\gB=\diag(-\gB_n-\gd_n),\quad R_X=\diag(\gd_n),
  \end{equation}
and $U$ is  unilateral shift on $l_2(\N)$. On the other hand, the
operator \eqref{I_12} is closely connected with the Krein string
spectral theory (see Subsection \ref{ss_II_krein}). Namely, in the
case when $\gB_n+\gd_n>0,\ n\in\N$, the difference expression
associated with \eqref{I_12} describes the motion of the
nonhomogeneous string with the mass distribution
\[
\mathcal{M}_\gB(x)=\sum_{x_{n-1}<x}\gd_n,\quad x\geq 0;\qquad x_n-x_{n-1}=\gB_n+\gd_n,\quad x_0=0.
\]
Based on this connection, we obtain the following criteria for the
operator $\rH_{X,\gB,0}$ to be self-adjoint, lower semibounded, and
discrete\footnote{Here we can consider the case when $\cI$ is a
bounded interval}
(Theorem \ref{th_delta'_sa} and Propositions \ref{prop_d'_disc1}, \ref{prop_d'_disc2} and \ref{cor_IV.3.4_01}) 

$(a)$ \emph{$\rH_{X,\gB,0}$ is self-adjoint if and only if either $\cI=\R_+$ or}
\[
\sum_{n\in\N}\bigl[\gd_{n+1}\sum_{i=1}^n(\gB_i+\gd_i)^2\bigr]=\infty.
\]

$(b)$ \emph{For the operator $\rH_{X,\gB,0}$ to be
lower semibounded it is necessary that}
\[
\frac{1}{\gB_n}\geq -C_1 \gd_n- \frac{1}{\gd_{n}},\quad\mathrm{and}\quad
\frac{1}{\gB_n}\geq -C_1 \gd_{n+1} - \frac{1}{\gd_{n+1}},\qquad n\in\N,
\]
\emph{and it is sufficient that}
\[
\frac{1}{\gB_n}\geq -C_2\min\{\gd_{n},\gd_{n+1}\},\qquad n\in\N
\]
with some positive constants $C_1,\ C_2>0$ independent of $n\in\N$.

$(c1)$ \emph{Let $\cI=\R_+$. The spectrum of $\rH_{X,\gB,0}$ is not discrete if one of the following conditions hold}
\begin{description}
\item $\bullet$ \quad$\lim_{n\to\infty}x_n\sum_{j=n}^\infty\gd_j^3>0$,
\item  $\bullet$ \quad $\gB_n\geq -C\gd_n^3$,\quad $n\in\N$,\quad $C>0$,
\item $\bullet$ \quad $\gB_n^{-}\le -C(\gd_n^{-1}+\gd_{n+1}^{-1})$, \quad $n\in\N$, \qquad ($\gB_n^{-}:=\gB_n$ if $\gB_n<0$ and $\gB_n^{-}:=-\infty$ if $\gB_n>0$).
\end{description}

 $(c2)$ \emph{If $\gd_n+\gB_n\geq 0$ for all $n\in\N$, then the spectrum of $\rH_{X,\gB,0}$ is discrete if and only if}
\[
\lim_{n\to\infty}x_n\sum_{j=n}^\infty\gd_j^3=0\quad
\text{and}\quad
\lim_{n\to\infty}x_n\sum_{j=n}^\infty(\gB_j+\gd_j)=0.
\]

Note that $(a)$ and $(c2)$ follow, respectively, from Hamburger's theorem and Kac--Krein discreteness criterion for the operator \eqref{I_12}. The results are demonstrated  by Example \ref{example_IV.3.4_01}.

In conclusion let us briefly describe the content of the paper. 

Section \ref{Sec_II_Prelim} is preparatory. It contains necessary
definitions and statements on theory of boundary triplets of symmetric operators and the Krein string spectral theory.

In Section \ref{Sec_III_direct_sums}, for arbitrary family of symmetric operators $\{S_n\}_{n\in\N}$, we investigate a direct sum $\Pi=\oplus_{n=1}^\infty \Pi_n$ of boundary triplets $\Pi_n$ for
$S_n^*$, $n\in \N$. We obtain two criteria
for 
$\Pi$ to form a boundary triplet for the operator $A^* =
\oplus_1^\infty S^*_n$ and regularization procedures for
$\Pi_n$ are given. 

Sections \ref{Subsec_IV.1_bt_interactions}--\ref{sec_delta'} are devoted to the spectral
analysis of operators with $\delta-$ and $\delta'-$interactions on
a discrete set $X$. We confine ourselves to the case $q\in
L^\infty$. 
In Section \ref{Subsec_IV.1_bt_interactions}, we construct
boundary triplets for the operator $\rH_{\min}^*$. Spectral
analysis of the Hamiltonians $\rH_{X,\gA,0}$ and $\rH_{X,\gB,0}$
are provided in Sections \ref{Subsec_IV.2_delta} and
\ref{sec_delta'}, respectively. More precisely, we study
self-adjointness of the minimal operators $\rH_{X,\gA,0}$ and
$\rH_{X,\gB,0}$, discreteness of their spectra, and their lower
semiboundedness.

In Section \ref{Sec_V}, we show that self-adjointness of the
operator $\rH_{X,\gA,q}$ with $\delta-$interactions is not stable
under perturbation by  positive unbounded potentials if $\gd_*=0$.

\textbf{Notation.} $\mathfrak{H}$, $\cH$ stand for the separable Hilbert spaces. $[\mathfrak{H}, \cH]$ denotes the set of bounded operators from $\mathfrak{H}$ to $\cH$; $[\mathfrak{H}]:=[\mathfrak{H},\mathfrak{H}]$ and ${\mathfrak S}_p(\mathfrak{H}),\ p\in(0,\infty)$, is the Neumann-Schatten ideal in $[\mathfrak{H}]$. $\mathcal{C}(\mathfrak{H})$ and $\widetilde{\mathcal{C}}(\mathfrak{H})$ are the sets of closed operators and linear relations in $\mathfrak{H}$, respectively.
Let $T$ be a linear operator in a Hilbert space $\mathfrak{H}$. In what follows, $\dom (T)$, $\ker (T)$, $\ran (T)$ are the domain, the kernel, the range of $T$, respectively; $\sigma(T)$, $\rho(T)$, and $\widehat{\rho}(T)$ denote the spectrum, the resolvent set, and the set of regular type points of $T$, respectively; $ R_T \left(\lambda \right):=\left( T-\lambda I\right)^{-1} $, $\lambda \in \rho(T)$, is the resolvent of $T$.

Let $X$ be a discrete subset of $\cI\subseteq\R$. By $W^{2,2}(\cI\setminus X)$, $W^{2,2}_0(\cI\setminus X)$, and
$W^{2,2}_{\loc}(\cI\setminus X)$ we denote the Sobolev spaces 
\begin{gather*}
W^{2,2}(\cI\setminus X):=\{f\in L^2(\cI): f,
f'\in AC_{loc}(\cI\setminus X), f''\in
  L^2(\cI)\},\\
W^{2,2}_0(\cI\setminus X):=\{f\in
W^{2,2}(\cI): f(x_k)=f'(x_k)=0,\,
\mbox{for all }  x_k\in X\},\\
W^{2,2}_{\comp}(\cI\setminus X):=\{f\in W^{2,2}(\cI\setminus X): \supp f\ \text{is compact in}\ \cI\}.
\end{gather*}
Let $I$ be a subset of $\Z$, $I\subseteq\Z$. We denote by $l_2(I,\cH)$ the Hilbert space of $\cH$-valued sequences such that $\|f\|^2=\sum_{n\in I}\|f_n\|_{\cH}^2<\infty$; $l_{2,0}(I,\cH)$ is a set of sequences with only finitely many values being nonzero; we also abbreviate $l_2:=l_2(\N,\C)$, $l_{2,0}:=l_{2,0}(\N,\C)$.

\section{Preliminaries}\label{Sec_II_Prelim}

\subsection{Boundary triplets and Weyl functions}\label{Subsec_II.1_btrips}

In this section we briefly review the notion of abstract boundary
triplets and associated Weyl functions in the extension theory of
symmetric operators (we refer to  \cite{DM91, DM95, Gor84} for a detailed study of
boundary triplets).

\subsubsection{Linear relations, boundary triplets, and self-adjoint extensions}
\label{sss_II.1.1_lr}

\textbf{1.} \ The set $\widetilde\cC(\cH)$ of closed linear relations in $\cH$ is the set of closed linear subspaces of $\cH\oplus\cH$.
Recall that $\dom(\Theta) =\bigl\{
f:\{f,f'\}\in\Theta\bigr\} $, $\ran(\Theta) =\bigl\{
f^\prime:\{f,f'\}\in\Theta\bigr\} $, and $\mul(\Theta) =\bigl\{
f^\prime:\{0,f'\}\in\Theta\bigr\} $ are the domain, the range, and the multivalued part of $\Theta$. A closed linear operator $A$ in $\cH$ is identified
with its graph $\gr(A)$, so that the set  $\cC(\cH)$  of closed linear
operators in $\cH$ is viewed as a subset of $\widetilde\cC(\cH)$.
In particular, a linear relation $\Theta$ is an operator if and
only if 
$\mul(\Theta)$ is trivial. For the
definition of the inverse, the resolvent set and the spectrum
of linear relations we refer to \cite{DS87}.
We recall that the adjoint relation
$\Theta^*\in\widetilde\cC(\cH)$ of $\Theta\in \widetilde\cC(\cH)$ is defined by
\begin{equation*}
\Theta^*= \left\{
\{h,h^\prime\}: (f^\prime,h)_{\cH}=(f,h^\prime)_{\cH}\,\,\text{for all}\,
\{f,f^\prime\}
\in\Theta\right\}.
\end{equation*}
A linear relation $\Theta$ is said to be {\it symmetric} if
$\Theta\subset\Theta^*$ and self-adjoint if $\Theta=\Theta^*$.

For a symmetric linear relation $\Theta\subseteq\Theta^*$ in $\cH$
the multivalued part $\mul(\Theta)$ is the orthogonal complement
of $\dom(\Theta)$ in $\cH$. Setting $\cH_{\rm
op}:=\overline{\dom(\Theta)}$ and $\cH_\infty=\mul(\Theta)$, one
arrives at the orthogonal decomposition  $\Theta= \Theta^{\rm
op}\oplus \Theta^\infty$, where  $\Theta^{\rm op}$ is a  symmetric
operator in $\cH_{\rm op}$ and is called \emph{the operator part} of $\Theta$, and
$\Theta^\infty=\bigl\{\bigl(\{0,f'\}\bigr):f'\in\mul(\Theta)\bigr\}$ is  a ``pure'' linear relation
 in $\cH_\infty$. 
\textbf{2.} \
Let $A$ be a densely defined closed symmetric operator in the
separable Hilbert space $\gH$ with equal deficiency indices
$\mathrm{n}_\pm(A)=\dim \cN_{\pm \I} \leq \infty,\ \ \cN_z:=\ker(A^*-z)$.

\begin{definition}[\cite{Gor84}]\label{def_ordinary_bt}
A triplet $\Pi=\{\cH,\gG_0,\gG_1\}$ is called a {\rm boundary
triplet} for the adjoint operator $A^*$ if $\cH$ is a Hilbert
space and $\Gamma_0,\Gamma_1:\  \dom(A^*)\rightarrow \cH$ are
bounded linear mappings such that the abstract Green identity
\begin{equation}\label{II.1.2_green_f}
(A^*f,g)_\gH - (f,A^*g)_\gH = (\gG_1f,\gG_0g)_\cH - (\gG_0f,\gG_1g)_\cH, \quad
f,g\in\dom(A^*),
\end{equation}
holds
and the mapping $\gG:=\{\Gamma_0,\Gamma_1\}:  \dom(A^*)
\rightarrow \cH \oplus \cH$ is surjective.
\end{definition}
First note that a boundary triplet 
for $A^*$ exists since the deficiency indices of $A$ are assumed to be
equal. Moreover, $\mathrm{n}_\pm(A) = \dim(\cH)$ and $A=A^*\upharpoonright\left(\ker(\Gamma_0) \cap \ker(\Gamma_1)\right)$ hold. Note also that a boundary triplet for $A^*$ is not unique.

A closed extension $\widetilde{A}$ of $A$ is called \emph{proper}
if $A\subseteq\widetilde{A}\subseteq A^*$. Two proper extensions
$\widetilde{A}_1$ and $\widetilde{A}_2$ of $A$ are called
\emph{disjoint} if
$\dom(\widetilde{A}_1)\cap\dom(\widetilde{A}_2)=\dom(A)$ and
\emph{transversal} if in addition
$\dom(\widetilde{A}_1)\dotplus\dom(\widetilde{A}_2)=\dom(A^*)$.
The set of  proper extensions of $A$ is denoted by $\Ext A.$
Fixing a boundary triplet $\Pi$ one can parameterize the set $\Ext
A$ in the following way.
\begin{proposition}[\cite{DM95}]\label{prop_II.1.2_01}
Let $A$ be as above and let $\Pi=\{\cH,\gG_0,\gG_1\}$ be a boundary
triplet for $A^*$. Then the mapping
     \begin{equation}\label{II.1.2_01A}
(\Ext A\ni)\ \widetilde A \to  \Gamma \dom(\widetilde A)
=\{\{\Gamma_0 f,\Gamma_1f \} : \  f\in \dom(\widetilde A) \} =:
\Theta \in \widetilde\cC(\cH)
     \end{equation}
establishes  a bijective correspondence between the sets $\Ext_A$
and  $\widetilde\cC(\cH)$. We put $A_\Theta :=\widetilde A$ where
$\Theta$ is defined by \eqref{II.1.2_01A}, i.e. $A_\Theta:=
A^*\upharpoonright \Gamma^{-1}\Theta=A^*\upharpoonright
\bigl\{f\in\dom(A^*): \ \{\Gamma_0f,\Gamma_1f\} \in\Theta\bigr\}$.
Then:
%
%
\item $(i)$ $A_\Theta$ is symmetric
(self--adjoint) if and only if $\Theta$ is symmetric, 
and $\mathrm{n}_\pm(A_\Theta)=\mathrm{n}_\pm(\Theta)$ holds.
\item $(ii)$ The extensions $A_\Theta$ and $A_0$ are disjoint
(transversal) if and only if $\Theta\in \cC(\cH)$ $\bigl(\Theta
\in [\cH]\bigr)$. In this case $A_\Theta$ admits a representation
$A_\Theta = A^*\!\upharpoonright\ker(\gG_1- \Theta\gG_0)$.
         \end{proposition}

It follows immediately from Proposition \ref{prop_II.1.2_01} that
the extensions
%
%
$A_0:=A^*\!\upharpoonright\ker(\gG_0)$ and
$A_1:=A^*\!\upharpoonright\ker(\gG_1)$
%
%
are self-adjoint. Clearly, $A_j=A_{\Theta_j} \ (j=0,1)$,
where the subspaces $\gT_0:= \{0\} \times \cH$ and
$\gT_1 := \cH \times \{0\}$ are self-adjoint relations in $\cH$.
Note that $\gT_0$ is a "pure" linear relation.

\subsubsection{Weyl functions, $\gamma$-fields, and Krein type formula for resolvents}\label{sss_II.1.3_weylsec}

\textbf{1.} \
In \cite{DM91, DM95}  the concept of the classical Weyl--Titchmarsh $m$-function
from the theory of Sturm-Liouville operators  was generalized to
the case of symmetric operators  with equal deficiency indices.
The role of abstract Weyl functions in the extension theory  is
similar to that of the classical Weyl--Titchmarsh $m$-function in the spectral
theory of singular Sturm-Liouville operators.

\begin{definition}[{\cite{DM91}}]\label{def_Weylfunc}
Let $A$ be a densely defined closed symmetric operator in $\gH$
with equal deficiency indices and let $\Pi=\{\cH,\gG_0,\gG_1\}$ be
a boundary triplet for $A^*$.
The operator valued functions $\gamma :\rho(A_0)\rightarrow  [\cH,\gH]$ and
$M:\rho(A_0)\rightarrow  [\cH]$ defined by
\begin{equation}\label{II.1.3_01}
\gamma(z):=\bigl(\Gamma_0\!\upharpoonright\cN_z\bigr)^{-1}
\qquad\text{and}\qquad M(z):=\Gamma_1\gamma(z), \qquad
z\in\rho(A_0),
\end{equation}
are called the {\em $\gamma$-field} and the {\em Weyl function},
respectively, corresponding to the boundary triplet $\Pi.$
\end{definition}
The $\gamma$-field $\gamma(\cdot)$ and the Weyl function
$M(\cdot)$ in \eqref{II.1.3_01}
are well defined. 
Moreover, both $\gamma(\cdot)$ and $M(\cdot)$ are holomorphic on
$\rho(A_0)$ and the following relations hold (see \cite{DM91})
\begin{gather}\label{II.1.3_02'}
\gamma(z)=\bigl(I+(z-\zeta)(A_0-z)^{-1}\bigr)\gamma(\zeta),
\\
\label{II.1.3_02}
M(z)-M(\zeta)^*=(z-\overline\zeta)\gamma(\zeta)^*\gamma(z), 
\\
   \label{II.1.3_02B}
\gamma^*(\overline z)
 = \Gamma_1(A_0 - z)^{-1}, \qquad z,\
\zeta\in\rho(A_0).
      \end{gather}
Identity \eqref{II.1.3_02}  yields that $M(\cdot)$ is an
$R_{\cH}$-function (or {\it Nevanlinna function}), that is,
$M(\cdot)$ is an ($[\cH]$-valued) holomorphic function on
$\C\setminus \R$ and
     \begin{equation}\label{II.1.3_03}
 \im z\cdot\im M(z)\geq 0,\qquad  M(z^*)=M(\overline
z),\qquad  z\in \C\setminus \R.
\end{equation}
 Besides, it follows  from \eqref{II.1.3_02}
that  $M(\cdot)$ satisfies $0\in \rho(\im M(z))$ for
$z\in\C\setminus\R$.
Since $A$ is densely defined, $M(\cdot)$ admits
an integral representation (see, for instance, \cite{DM95})
\begin{equation}\label{WF_intrepr}
M(z)=C_0+\int_{\R}\left(\frac{1}{t-z}-\frac{t}{1+t^2}\right)d\Sigma_M(t),\qquad z\in\rho(A_0),
\end{equation}
where $\Sigma_M(\cdot)$ is an
operator-valued Borel measure on $\R$ satisfying
$\int_\R \frac{1}{1 + t^2}d\Sigma_M(t) \in [\cH]$ and $C_0 = C_0^*\in [\cH]$. The integral in (\ref{WF_intrepr}) is understood in the strong sense.

In contrast to spectral measures of self-adjoint operators the
measure $\Sigma_M(\cdot)$ is not necessarily orthogonal. However,
the measure $\Sigma_M$ is uniquely determined by the Nevanlinna
function $M(\cdot)$.  The operator-valued measure $\Sigma_M$ is
called \emph{the spectral measure} of $M(\cdot)$. If $A$ is a
simple symmetric operator, then  the Weyl function $M(\cdot)$
determines the pair $\{A,A_0\}$ up to unitary equivalence (see
\cite{DM95, KL71}). Due to this fact,  spectral
properties of $A_0$ can be expressed in terms of $M(\cdot)$.

\textbf{2.} \  The following result provides a description of resolvents and
spectra of proper extensions of the operator $A$ in terms of the
Weyl function $M(\cdot)$ and the corresponding boundary
parameters.
   \begin{proposition}[\cite{DM91}]\label{prop_II.1.4_spectrum}
For any  $\Theta\in \wt\cC(\cH)$  the following Krein type formula
holds
\begin{equation}\label{II.1.4_01}
(A_\Theta - z)^{-1} - (A_0 - z)^{-1} = \gamma(z) (\Theta -
M(z))^{-1} \gamma^*({\overline z}), \quad z\in \rho(A_0)\cap
\rho(A_\Theta).
\end{equation}
Moreover, if $z\in \rho(A_0)$, then
\[
z\in\sigma_i(A_\Theta) \quad \Leftrightarrow\quad 0\in \sigma_i(\Theta-M(z)),\qquad i\in\{\rm p,\ c,\ r\}.
\]
\end{proposition}
Formula \eqref{II.1.4_01} is a generalization of the well known Krein
formula for canonical resolvents (cf. \cite{Akh_Glz}). We note also that all objects in
\eqref{II.1.4_01} are expressed in terms of the boundary triplet $\Pi$.

The following result is deduced from \eqref{II.1.4_01}
%
%
\begin{proposition}[\cite{DM91}]\label{prop_II.1.4_02}
Let $\Pi=\{\cH,\gG_0,\gG_1\}$  be a boundary triplet for $A^*,$
$\Theta_1,\Theta_2\in \wt\cC(\cH)$.
Then:

\item  $(i)$ for any $z
\in\rho(A_{\Theta_1})\cap\rho(A_{\Theta_2}),\
\zeta\in\rho(\Theta_1)\cap\rho(\Theta_2)$ the following
equivalence holds
\begin{equation}\label{II.1.4_02}
(A_{\Theta_1}-z)^{-1} - (A_{\Theta_2}-z)^{-1}\in{\mathfrak
S}_p(\gH)
\quad \Longleftrightarrow \quad (\Theta_1 - \zeta)^{-1}-
(\Theta_2 - \zeta )^{-1}\in{\mathfrak S}_p(\cH).
\end{equation}

\item $(ii)$\  If, in addition, $\Theta_1, \Theta_2\in\cC(\cH)$ and
$\dom(\Theta_1) = \dom(\Theta_2)$, then
      \begin{equation}\label{II.1.4_03}
\overline{\Theta_1 - \Theta_2} \in{\mathfrak S}_p(\cH)
\Longrightarrow (A_{\Theta_1}-z)^{-1} -
(A_{\Theta_2}-z)^{-1}\in{\mathfrak S}_p(\gH).
      \end{equation}
\item $(iii)$\  Moreover, if  $\Theta_1, \Theta_2\in[\cH]$, then
implication \eqref{II.1.4_03} becomes equivalence.
%
%
%
\end{proposition}
\subsubsection{Extensions of a nonnegative operator}\label{sss_II.1.5_nno}

Assume that  a symmetric operator $A\in \cC(\mathfrak{H})$ is
nonnegative. Then the set $\Ext_A(0,\infty)$ of its  nonnegative
self-adjoint extensions is non-empty (see \cite{Akh_Glz, Kato66}).
Moreover, there is a maximal 
nonnegative extension $A_{\rm F}$ (also called \emph{Friedrichs'}
or \emph{hard} extension) and there is a  minimal nonnegative
extension $A_{\rm K}$ (\emph{Krein's} or \emph{soft} extension)
satisfying
$$
(A_F+x)^{-1} \le (\wt A + x)^{-1} \le (A_K + x)^{-1}, \qquad x\in
(0,\infty), \quad \wt A\in  \Ext_A(0,\infty),
$$
(for detail we refer the reader to \cite{Akh_Glz,Gor84}).
   \begin{proposition}[\cite{DM91}]\label{prop_II.1.5_01}
Let \ $\Pi=\{\cH,\Gamma_0,\Gamma_1\}$ be a boundary triplet for
$A^*$ such that $A_0=A_0^*\geq 0$.
Let $M(\cdot)$ be the corresponding Weyl function.  Then $A_0=A_{\rm F}\ \
(A_0=A_{\rm K})$ if and only if
      \begin{equation}\label{II.1.5_01}
\lim_{x\downarrow-\infty}(M(x)f,f)=-\infty,\qquad \bigl(\lim_{x\uparrow0}(M(x)f,f)=+\infty\bigr),\qquad f\in\cH\setminus\{0\}.
    \end{equation}
       \end{proposition}
It is said that $M(\cdot)$ \emph{uniformly tends to $-\infty$} for
$x\to-\infty$ if for any $a>0$ there exists $x_a<0$ such that
$M(x_a)<-a\cdot I_{\cH}$. In this case we will write
$M(x)\rightrightarrows-\infty, \ x\to-\infty$.
%
        \begin{proposition}[\cite{DM91}]\label{prop_II.1.5_02}
Let $A$ be a non-negative symmetric operator in $\mathfrak{H}$. Assume that
$\Pi=\{\cH,\Gamma_0,\Gamma_1\}$ is a boundary triplet for $A^*$ such
that $A_0=A_{\rm F}$, and let also $M(\cdot)$ be the corresponding Weyl
function. Then the  following assertions \item $(i)$ \ a linear
relation $\Theta\in \widetilde\cC_{\rm self}(\cH)$ is semibounded
below, \item $(ii)$ \ a self-adjoint extension $A_\Theta$ is
semibounded below, \item are equivalent if and only if
$M(x)\rightrightarrows-\infty$ for $x\to-\infty$.
    \end{proposition}

\subsubsection{Generalized boundary triplets and boundary relations}\label{subsec_II.2_genbt}

In many applications the notion of a boundary triplet is too
strong. Therefore it makes sense to relax its definition. To do this
we follow \cite[Section 6]{DM95}.
   \begin{definition}[{\cite{DM95}}]\label{def_II.2.1_generalized_bt}
Let $A$ be a closed densely defined symmetric operator in $\gH$
with equal deficiency indices. Let $A_* \supseteq A$ be a not
necessarily closed extension of $A$ such that $(A_*)^* = A$. A
triplet $\Pi=\{\cH,\gG_0,\gG_1\}$  is called a \emph{generalized
boundary triplet for $A^*$} if $\cH$ is a Hilbert space and
$\gG_j: 
\dom(A_*)\to\cH$, $j=0,1$, are linear mappings
such that

\item $(G1)$ $\gG_0$ is surjective,

 \item $(G2)$
$A_{*0}:=A_*\upharpoonright\ker(\gG_0)$ is a self-adjoint
operator,

\item $(G3)$ Green's formula holds
   \begin{equation}\label{II.2.1_03}
(A_*f, g)_{\gH} - (f, A_*g)_{\gH} = (\gG_1f, \gG_0 g)_{\cH} -
(\gG_0 f, \gG_1 g)_{\cH},\qquad f,g\in\dom(A_*)=\dom(\Gamma).
       \end{equation}
       \end{definition}
Note that one always has $A \subseteq A_* \subseteq A^* =
\overline{A_*}$. The following properties of a generalized
boundary triplet have been  established in \cite{DM95}.
    \begin{lemma}[{\cite{DM95}}]\label{lem_II.2.1_01}
Let  $\Pi=\{\cH,\gG_0,\gG_1\}$ be a generalized
boundary triplet  for $A^*$. 
 Then:
 \item $(i)$
$\cN^*_{z}:=\dom(A_*)\cap\cN_{z}$ is dense in $\cN_{z}$ and
$\dom(A_*) =\dom(A_0) + \cN^*_{z}$.
\item $(ii)$ $\overline{\gG_1\dom(A_0)}=\cH$.
\item $(iii)$ $\ker(\gG) = \dom(A)$\ \ and\ \ ${\overline{\ran(\gG)}} = \cH \oplus \cH$,
where $\gG: = \{\gG_0,\gG_1\}$.
\end{lemma}

For any generalized boundary triplet
$\Pi=\{\cH,\gG_0,\gG_1\}$   we set
$A_{*j}:=A^*\lceil\ker(\gG_j)$, $j= 0,1$.
Note  that the extensions $A_{*0}$ and $A_{*1}$
are always  disjoint but not  necessarily transversal.

Starting with Definition \ref{def_II.2.1_generalized_bt}, one can  introduce concepts of
the (generalized) $\gamma$-field $\gamma(\cdot)$ and the Weyl
function $M(\cdot)$ corresponding to a
generalized boundary triplet $\Pi$ in just the same way as
it was done for (ordinary)
boundary triplets (for detail see \cite{DM95}). Let us mention only
the following proposition (cf. \cite[Proposition 6.2]{DM95}).
     \begin{proposition}[\cite{DM95}]\label{prop2.10}
Let $\Pi=\{\cH,\gG_0,\gG_1\}$ be a generalized boundary triplet
for   $A^*$, $A_{*}=A^*\lceil \dom(\Gamma)$, and let $M(\cdot)$ be the
corresponding Weyl function. Then:

\item $(i)$\  $M(\cdot)$ is an $[\cH]$-valued Nevanlinna function
satisfying $\ker(\im M(z))= \{0\}, \ z\in \C_+$.

\item  $(ii)$\  $\Pi$ is an ordinary boundary triplet  if and only
if $0\in \rho(\im M(i))$.
\end{proposition}

We also need the following definition.
\begin{definition}[\cite{DHMS06}]\label{def_II.2.1_02_boundaryrelation}
Let $A$ be as in Definition \ref{def_II.2.1_generalized_bt}
and let $\cH$ be an auxiliary Hilbert space. A linear relation
(multi-valued mapping) $\Gamma:  \gH\to\cH^2$ is called a boundary
relation for $A^*$ if:

$(i)\  \dom (\Gamma)$ is dense in $\dom(A^*)$, and identity
    \begin{equation}\label{II.2.1_02}
(A_*f, g)_{\gH} - (f, A_*g)_{\gH} = (l', h)_{\cH} - (l, h')_{\cH},
    \end{equation}
where  $A_* = A^*\lceil  \dom(\Gamma)$,  holds for every $\{
f,{\hat l}\},\{g, {\hat h}\}\in\Gamma$,

$(ii)\  \Gamma$ is maximal in the sense that  if  
$\{\hat g, \hat h\} \in\gH^2\oplus\cH^2$ satisfies the identity
$(A_*f, g) - (f, g') = (l', h) - (l, h')$
for every $\{f, \hat l\}\in\Gamma,$  
then $\{g, \hat h\}\in\Gamma$.

Here $f, {g}\in\dom\Gamma(\subset\gH)$, $g'\in \gH,\  \hat
g:=\{g,g'\}$ and ${\hat h}= \{h, h'\}, {\hat l} = \{l, l'\}
\in\ran\Gamma(\subset\cH^2)$.
          \end{definition}

Note that in general $\Gamma$ is multi-valued.   If it is
single-valued, it splits $\Gamma =\{\Gamma_0, \Gamma_1\}$ and
Green's identity \eqref{II.2.1_02} takes usual form \eqref{II.2.1_03}.

\subsection{Nonhomogeneous Krein--Stieltjes string}\label{ss_II_krein}

In this subsection, we collect some facts on Jacobi operators of a special form.  Namely, consider two sequences with positive elements $m=\{m_n\}_{n=1}^\infty$ and $l=\{l_n\}_{n=1}^\infty$, $m_n,\ l_n>0,\ n\in\N$. Next, consider the matrix
\begin{equation}\label{A_01}
J_{m,l}=\left(\begin{array}{cccc}
\frac{1}{m_1}\frac{1}{l_1} & \frac{1}{l_1\sqrt{m_1m_2}} & 0&   \dots\\
\frac{1}{l_1\sqrt{m_1m_2}} & \frac{1}{m_2}\bigl(\frac{1}{l_1}+\frac{1}{l_2}\bigr) & \frac{1}{l_2\sqrt{m_2m_3}}&  \dots\\
0 & \frac{1}{l_2\sqrt{m_2m_3}} & \frac{1}{m_3}\bigl(\frac{1}{l_2}+\frac{1}{l_3}\bigr)&  \dots\\
\dots & \dots & \dots&  \dots
\end{array}\right).
     \end{equation}
With unilateral shift $U$ in $l_2(\N)$, $U\mathrm{e}_n=\mathrm{e}_{n+1}$, $n\in\N$, where $\{\mathrm{e}_n\}_{n\in\N}$ is the standard orthonormal basis in $l_2$, the matrix  $J_{m,l}$  can be  written
as
    \begin{equation}\label{A_02}
J_{m,l}=M^{-1/2}(I+U)L^{-1}(I+U^*)M^{-1/2}, \qquad M=\diag(m_n),\quad L=\diag(l_n).
      \end{equation}
It is known that the difference expression associated with
$J_{m,l}$ has a useful mechanical interpretation, related to the
Krein string theory (for detail we refer the reader to
\cite[Appendix, pp.232--236]{Akh} and \cite{KK71}). Namely, define
the function
   \begin{equation}\label{A_03}
\mathcal{M}(x)=\sum_{x_{n-1}<x}m_n,\quad x\in
[0,\mathcal{L});\qquad \mathcal{L}=\sum_{n=1}^\infty l_n,\quad
x_n-x_{n-1}=l_n,\quad x_0=0.
     \end{equation}
Then the equation of motion of a nonhomogeneous string with the
mass distribution $\mathcal{M}$ is the same as the difference
equation associated with the Jacobi matrix $J_{m,l}$ (strings with
discrete mass distributions are called Stieltjes strings).

Further, associated with the matrix $J_{m,l}$  one introduces the minimal Jacobi operator in $l_2(\N)$ (see \cite{Akh, Ber68}). We denote it also by $J_{m,l}$. By \emph{Hamburger's theorem} \cite[Theorem 0.5]{Akh},
the operator $J_{m,l}$ is self-adjoint if and only if
      \begin{equation}\label{A_04}
\sum_{n=1}^\infty m_{n+1}x_n^2=\infty.
     \end{equation}
A discreteness criterion for the nonhomogeneous string was
obtained by Kac and Krein in \cite{KK58} (see also \cite[\S 11]{KK71}). Applying  their result to the operator \eqref{A_01}, we arrive at the following criterion.
     \begin{theorem}[\cite{KK58}]\label{th_Append_01}
Assume \eqref{A_04} and set $\mathcal{M}(\mathcal{L}):=\lim_{x\uparrow \mathcal{L}}\mathcal{M}(x)=\sum_{n=1}^\infty m_n$. Then $J_{m,l}= J_{m,l}^*$ has
discrete spectrum if and only if
 \begin{description}\label{A_05}
\item in the case $\mathcal{L}=\infty$, \quad $\lim_{n\to \infty}x_n\sum_{j=n}^\infty m_j=0$\quad (the latter yields $\mathcal{M}(\mathcal{L})<\infty$);
\item in the case $\mathcal{M}(\mathcal{L})=\infty$ and $\mathcal{L}<\infty$, \quad $\lim_{n\to \infty}(\mathcal{L}-x_n)\sum_{j=1}^n m_j=0$. 
\end{description}
      \end{theorem}
    \begin{remark}\label{rem_A01}
If condition \eqref{A_04} does not hold, then
$\mathrm{n}_\pm(J_{m,l})=1$ and hence any self-adjoint extension
of $J_{m,l}$ has discrete spectrum.

Note also that for
$J_{m,l}$ to be discrete it is necessary that either
$\{m_n\}_{n=1}^\infty\in l_1$ or $\{l_n\}_{n=1}^\infty\in
l_1$.
\end{remark}
%
%
\section{Direct sums of symmetric operators and boundary triplets}\label{Sec_III_direct_sums}
%
%
%

\subsection{Direct sum of boundary triplets as a boundary relation}\label{subsec_III.1_dirsum}

Let $S_n$ be a densely defined symmetric operator in a Hilbert
space $\mathfrak{H}_n$ with equal deficiency indices,
$\mathrm{n}_+(S_n) = \mathrm{n}_-(S_n) \le \infty,\
n\in \N.$ Consider the operator $A:= \oplus^{\infty}_{n=1}S_n$
acting in a Hilbert direct sum $\mathfrak{H} :=
\oplus_{n=1}^{\infty}\mathfrak{H}_n$ of spaces $\mathfrak{H}_n.$
By definition, $ \mathfrak{H} = \{f= \oplus^{\infty}_{n=1} f_n: \
f_n\in \mathfrak{H}_n, \ \sum^{\infty}_{n=1}\|
 f_n\|^2<\infty\}.$  We also denote by $\mathfrak{H}^0$ the linear manifold
 consisting of vectors   $f= \oplus^{\infty}_{n=1} f_n\in
\mathfrak{H}$ with finitely many nonzero entries. Clearly,
   \begin{equation}\label{III.1_01}
A^* = \oplus^{\infty}_{n=1}S^*_n,\qquad \dom(A^*) =  \{f =
\oplus^{\infty}_{n=1} f_n\in \mathfrak{H}:\
 f_n\in\dom(S^*_n),\ \  \sum^{\infty}_{n=1}\|S^*_n
 f_n\|^2<\infty\}.
   \end{equation}
We provide  the domains $\dom(S^*_n)=: \gH_{n+}$ and
$\dom(A^*)=: \gH_+$ with the graph norms  $\|f_n\|^2_{\gH_{n+}} :=
\|f_n\|^2 + \|S^*_nf_n\|^2$ and  $\|f\|^2_{\gH_+}  := \|f\|^2 + \|A^*f\|^2= \sum_n
\|f_n\|^2_{\gH_{n+}}$, respectively.

Further, let   $\Pi_n=\{\cH_n, \Gamma^{(n)}_0, \Gamma^{(n)}_1\}$
be a boundary triplet  for $S^*_n$, $n\in \N$.
By $\|\Gamma_j^{(n)}\|$ we denote  the norm of the linear mapping $\Gamma^{(n)}_j\in[\gH_{n+},\cH_n]$, $
j=0,1$,  $n\in \N$.

Let $\cH :=\oplus_{n=1}^{\infty}\cH_n$ be a Hilbert direct sum of
$\cH_n$. Define  mappings $\Gamma_0$ and $\Gamma_1$
by setting
   \begin{equation}\label{III.1_02}
\Gamma_j  := \oplus_{n=1}^{\infty} \Gamma^{(n)}_j,\qquad
\dom(\Gamma_j) = \bigl\{f = \oplus^{\infty}_{n=1} f_n
\in\dom(A^*): \ \sum^{\infty}_{n=1}\|\Gamma^{(n)}_j
f_n\|^2_{\cH_n} <\infty\bigr\}.
   \end{equation}
Clearly $\mathfrak{H}_+\cap\mathfrak{H}^0\subset
\dom(\Gamma_j)\subset\dom(A^*),$ and $\dom(\Gamma)
:=\dom(\Gamma_1)\cap\dom(\Gamma_0)$ is dense in $\gH_+$
since $\mathfrak{H}_+\cap\mathfrak{H}^0$ is dense in
$\mathfrak{H}_+.$ Define the operators $S_{nj} := S_{n}^*\lceil
\ker \Gamma^{(n)}_j$ and ${\widetilde A}_j :=
\oplus^{\infty}_{n=1}S_{nj}$, $j=0,1$. Then ${\widetilde
A}_0$ and ${\widetilde A}_1$ are self-adjoint extensions of $A$.
Note that ${\widetilde A}_0$ and ${\widetilde A}_1$ are disjoint
but not necessarily transversal.

Finally, we  set
\begin{equation}\label{III.1_03}
A_* = A^*\lceil\dom(\Gamma)\quad \text{and}\quad A_{*j} :=
A_*\lceil\ker(\Gamma_j),\quad  j=0,1.
   \end{equation}
Clearly, $A_{*j}$ is symmetric (not necessarily self-adjoint or
even closed!) extension of $A$, $A_{*j} \subset {\widetilde
A}_j$, $j=0,1$,  and
\[
\dom (A_{*j})= \{f =
\oplus^{\infty}_{n=1} f_n\in \mathfrak{H}:\
 f_n\in \ker \Gamma^{(n)}_j,\ \  \sum_{n}\bigl(\|S^*_n
 f_n\|^2 + \|\Gamma_{j'}^{(n)}f_n\|^2\bigr) <\infty\},\quad (0':=1,\ 1':=0).
\]

       \begin{definition}\label{def_III.1_01}
Let  $\Gamma_j$ be defined by \eqref{III.1_02} and $\cH
=\oplus_{n=1}^{\infty}\cH_n$. A collection $\Pi=\{\cH, \Gamma_0,
\Gamma_1\}$ will be called  a \emph{direct sum of boundary triplets} and
will be assigned as $\Pi:= \oplus^{\infty}_{n=1}\Pi_n$.
    \end{definition}

By Definition \ref{def_ordinary_bt}, for a direct sum $\Pi=
\oplus^{\infty}_{n=1}\Pi_n$ to form a
boundary triplet for  $A^*=\oplus_{n=1}^\infty S^*_n$ it is necessary (but not sufficient!)  that

  (a)  $A_{*0}$ and $A_{*1}$ are  self-adjoint,

(b)  $A_{*0}$ and $A_{*1}$ are transversal,

(c)  $\dom(\Gamma)= \dom(A^*)$,

(d)  $\Gamma_0$ and $\Gamma_1$ are closed and bounded as
mappings from $\gH_+$ to $\cH$.\\
\noindent It might happen that all of these conditions are
violated for the direct sum $\Pi.$ Nevertheless, we will show that
$\Pi$ is a boundary relation for the operator $A^*$ in the sense
of Definition \ref{def_II.2.1_02_boundaryrelation}.
   \begin{theorem}\label{th_III.1_01}
Let  $\Pi_n=\{\cH_n, \Gamma^{(n)}_0, \Gamma^{(n)}_1\}$ be a
boundary triplet  for $S^*_n$, $M_n(\cdot)$ the corresponding Weyl
function, $n\in \N$. Let also $A^* = \oplus^{\infty}_{n=1}S^*_n$
and $\Pi = \oplus_{n=1}^{\infty}\Pi_n$.
Then:
   \item $(i)$\  $\Pi = \{\cH,\Gamma_0,\Gamma_1\}$  forms  a
boundary relation for $A^*$ with single-valued  $\Gamma =
\{\Gamma_0,\Gamma_1\}$.
   \item $(ii)$\ The corresponding Weyl
function is
        \begin{equation}\label{III.1_05}
M(z)=\oplus^{\infty}_{n=1}M_n(z).
        \end{equation}
\item $(iii)$\   $\ran\Gamma
=\ran\bigl(\{\Gamma_0,\Gamma_1\}\bigr)$ is dense in
$\cH\oplus\cH$.
     \item $(iv)$ The mapping $\Gamma: \gH_+ \to \cH \oplus \cH$
     is closed and the mappings $\Gamma_j:\gH_+ \to \cH $ are closable.  
\item $(v)$ If $\overline\Gamma_j$ is a closure of $\Gamma_j$, then
the following equivalences hold
   \begin{equation}\label{III.1_04}
\dom(\overline\Gamma_j) = \gH_+ \  \Longleftrightarrow \
\overline\Gamma_j\in[\gH_+, \cH] \  \Longleftrightarrow \
\sup_{n\in\N}\|\Gamma^{(n)}_j\| := C_j <\infty,\footnote{$\|\Gamma_j^{(n)}\|$ stands for the the norm of $\Gamma_j^{(n)}$ as a bounded linear mapping from $\gH_{n+}$ to $\cH_n$} \quad j=0,1.
   \end{equation}
In particular, $\dom(\Gamma) = \dom(\Gamma_0)\cap \dom(\Gamma_1)
=\gH_+$ if and only if $\max\{C_0, C_1\} <\infty$.

\item $(vi)$\  The operator $A_{*j}$ $($see \eqref{III.1_03}$)$ 
is essentially self-adjoint and $\overline {A_{*j}} = {\widetilde
A}_j = \oplus^{\infty}_{n=1}S_{nj}$, $j=0,1$.

\item $(vii)$\  $A_{*j}$ 
is self-adjoint, $A_{*j} = {\widetilde A}_j =
\oplus^{\infty}_{n=1}S_{nj},$
 whenever $C_{j'} =\sup_{n\in \N}\|\Gamma_{j'}^{(n)}\|<\infty,\ j=0,1$. 
If in addition  ${\widetilde A}_0$ and ${\widetilde A}_1$ are transversal, then 
$
A_{*j} = (A_{*j})^* \  \Longleftrightarrow   \  C_{j'} =
\sup_{n\in \N}\|\Gamma^{(n)}_{j'}\|  <\infty.
$
               \end{theorem}
   \begin{proof}
$(i)$\  Let us prove Green's identity  \eqref{II.2.1_03}.
By \eqref{III.1_01}--\eqref{III.1_03} and Definition \ref{def_III.1_01},  for
 $f=\oplus^{\infty}_{n=1} f_n,$\ $g =\oplus^{\infty}_{n=1}
g_n\in\dom(A_*) = \dom (\Gamma)$ we get
   \begin{gather}
(A_*f, g)_\gH- (f, A_* g)_\gH = \sum_{n=1}^\infty[(S^*_n f_n, g_n)_{\gH_n}- (f_n,S^*_n g_n)_{\gH_n}] \notag \\
= \sum_{n=1}^\infty\left[(\Gamma^{(n)}_1
f_n,\Gamma^{(n)}_0g_n)_{\cH_n}
 - (\Gamma^{(n)}_0 f_n,\Gamma^{(n)}_1 g_n)_{\cH_n}\right] = (\Gamma_1
f,\Gamma_0 g)_{\cH}-(\Gamma_0 f,\Gamma_1 g)_{\cH}.\label{III.1_06}
   \end{gather}
Note, that the series in the above equality converge due to
\eqref{III.1_01} and \eqref{III.1_02}.

To prove the maximality assumption assume that Green's identity
    \begin{equation}\label{III.1_07}
(A_*f, g)_\gH - (f, g')_\gH = (\Gamma_1 f, h)_\cH - (\Gamma_0 f, h')_\cH
  \end{equation}
holds for every $f\in\dom(A_*)$ and some $g, \ g' \in\gH$, and
$\{h,h'\}\in\cH\oplus\cH$. Let us show that $g\in\dom(A_*)$ and
$\Gamma g = \{\Gamma_0 g,\Gamma_1 g\} = \{h,h'\}$. If
$f\in\dom(A)$, equality \eqref{III.1_07} yields  $g\in\dom(A^*)$
and $g' = A^*g$. Hence $g=\oplus^{\infty}_{n=1} g_n,\
g_n\in\dom(S^*_n)$, and $A^* g=\oplus^{\infty}_{n=1} S^*_n g_n$.
Setting $f=f_n\in \dom(S_n^*)$ in \eqref{III.1_07} and noting that
$h=\oplus^{\infty}_{n=1} h_n, h'=\oplus^{\infty}_{n=1}
h'_n\in\cH$, we get
    \begin{equation}\label{III.1_08}
(S^*_n f_n, g_n)_{\gH_n} - (f_n, S^*_n g_n)_{\gH_n}  = (\Gamma^{(n)}_1 f_n,
h_n)_{\cH_n}\
 -\ (\Gamma^{(n)}_0f_n, h'_n)_{\cH_n}, \qquad n\in \N.
    \end{equation}
Since  $\Pi_n$ is a boundary triplet for $S_n^*,$  $\Gamma^{(n)}_0
g_n=h_n$ and $\Gamma^{(n)}_1 g_n = h'_n, \ n\in \N.$  Moreover,
the inclusion  $\{h,h'\}\in\cH\oplus\cH$ yields
    \begin{equation}\label{III.1_09}
\sum^{\infty}_{n=1}\bigl(\|\Gamma^{(n)}_0 g_n\|^2_{\cH_n} +
\|\Gamma^{(n)}_1 g_n\|^2_{\cH_n} \bigr) =
\sum^{\infty}_{n=1}\bigl(\|h_n \|^2_{\cH_n} + \|h'_n\|^2_{\cH_n}
\bigr) < \infty.
    \end{equation}
Inequality \eqref{III.1_09} means that $g\in\dom(A_*)= \dom(\Gamma)$ and
$\Gamma g = \{\Gamma_0 g,\Gamma_1 g\} = \{h,h'\}$. This proves the
maximality condition.

$(ii)$\ Straightforward.

$(iii)$\  Denote by $\cH^0$ the linear manifolds of vectors $h=
\oplus^{\infty}_{n=1} h_n\in \cH$ having finitely many nonzero
entries. Clearly $\cH^0$ is dense in $\cH$. It remains to note
that $\cH^0=\ran\bigl( \Gamma\lceil
(\mathfrak{H}_+\cap\mathfrak{H}^0)\bigr)\subset\ran(\Gamma),$
since $\ran(\Gamma^{(n)}) = \cH_n\oplus\cH_n$, $n\in {\N}$.

$(iv)$\  Let $f_k=\oplus^{\infty}_{n=1}f_{kn},\  \varphi =
\oplus^{\infty}_{n=1}\varphi_{n}  \in\gH_+$, and $\|f_k-\varphi\|_{\gH_+} \to
0$ and
  \begin{equation}\label{III.1_10}
\lim_{k\to\infty}\Gamma f_k=\lim_{k\to\infty}\{\Gamma_0
f_k,\Gamma_1 f_k\} = \{h,h'\} = \{\oplus^{\infty}_{n=1}h_n,
\oplus^{\infty}_{n=1}h'_n\} \in\cH\oplus\cH.
  \end{equation}
Let us prove that $\varphi\in\dom(A_*)$ and $\Gamma \varphi=\{h,h'\}$. Since
$\Gamma_j f_k=\oplus^{\infty}_{n=1}\Gamma_j^{(n)}f_{kn}$, by
\eqref{III.1_10} we get
    \begin{equation}\label{III.1_11}
\lim_{k\to\infty}\Gamma_0^{(n)}f_{kn}=h_n, \qquad
\lim_{k\to\infty}\Gamma_1^{(n)}f_{kn} = h'_n, \qquad n\in \N.
    \end{equation}
Since $\lim_{k\to\infty}\|f_{kn}-\varphi_n\|_{\gH_{n+}}=0$ 
and the mappings
$\Gamma^{(n)}=\{\Gamma^{(n)}_0,\Gamma^{(n)}_1\}:\gH_{n+}\to\cH_n\oplus\cH_n$
are closed (in fact, continuous), \eqref{III.1_11} yields
  \begin{equation}\label{III.1_12}
\varphi_n\in\gH_{n+}=\dom(S^*_n)\quad \text{and} \quad
\Gamma^{(n)}\varphi_n = \{h_n, h'_n\}.
  \end{equation}
In turn, since $\varphi\in\gH_+ = \dom(A^*)$ and
    \begin{equation}\label{III.1_13}
\sum^{\infty}_{n=1}\bigl(\|\Gamma^{(n)}_0 \varphi_n\|^2_{\cH_n} +
\|\Gamma^{(n)}_1 \varphi_n\|^2_{\cH_n} \bigr) =
\sum^{\infty}_{n=1}\bigl(\|h_n \|^2_{\cH_n}  + \|h'_n\|^2_{\cH_n}
\bigr) < \infty,
    \end{equation}
we obtain $\varphi\in\dom(A_*)$ and $\Gamma \varphi=\{\Gamma_0
\varphi,\Gamma_1 \varphi\}=\{h,h'\}$. Hence $\Gamma$ is
closed.

$(v)$\ By $(iv)$, the mapping $\Gamma$ is closed. Hence $(v)$ is implied
by the closed graph theorem.

$(vi)$\  Clearly, $\mathfrak{H}_+\cap\mathfrak{H}^0\subset
\dom(\widetilde{A}_j)$. Hence $\dom(A_{*j})$ is dense in
$\dom(\widetilde{A}_j)$ (in the graph topology).

$(vii)$\ Let  $C_1<\infty.$  Let us prove the self-adjointness of
$A_{*0}.$   
Since $A_{*0} \subset \widetilde{A}_0$, it suffices to show that
$\dom(\widetilde{A}_0) \subset \dom(A_*).$
 Let $f = \oplus^{\infty}_{n=1} f_n\in\dom(\widetilde{A}_0)$.
 Clearly $f\in\dom(\Gamma_0)$ since $f_n\in \ker \Gamma^{(n)}_0$. Let us show that $f\in\dom(\Gamma_1)$. According to the
second J. von Neumann formula,
  \begin{equation}\label{III.1_14}
f_n = f_{S_n} + (I + U_n)f_n(\I),  \qquad f_{S_n}\in \dom (S_n),
\quad  f_n(\I)\in \cN^{(n)}_\I := \cN_\I(S_n),
  \end{equation}
where $U_n$ is an isometry from $\cN^{(n)}_\I$ onto
$\cN^{(n)}_{-\I}$. Since $f\in\dom(A^*)$, it follows form
\eqref{III.1_14} that
   \begin{eqnarray*}  \label{III.1_15}
4\sum_{n=1}^\infty\|f_n(\I)\|_{\mathfrak{H}_n}^2 = \sum_{n=1}^\infty\|(I
+ U_n)f_n(\I)\|_{\gH_{n+}}^2
\le\sum_{n=1}^\infty\bigl(\|f_{S_n}\|_{\gH_{n+}}^2 + \|(I +
U_n)f_n(\I)\|_{\gH_{n+}}^2\bigr) =
\sum_{n=1}^\infty\|f_n\|_{\mathfrak{H}_+}^2 <\infty.
   \end{eqnarray*}
Hence $f(\I):=\oplus^{\infty}_{n=1} f_n(\I)\in\dom(A^*)$. Combining this
fact with the assumption   $C_1<\infty$, we get from \eqref{III.1_14}
   \begin{equation}\label{III.1_16}
\sum_{n=1}^\infty\|\Gamma_1^{(n)}f_n\|_{\cH_n}^2  =
\sum_{n=1}^\infty\|\Gamma_1^{(n)}(I + U_n)f_n(\I)\|_{\cH_n}^2 \le
4C_1^2 \sum_{n=1}^\infty\|f_n(\I)\|_{\mathfrak{H}_{n+}}^2 \le 8C_1^2
\sum_{n=1}^\infty\|f_n(\I)\|_{\mathfrak{H}_n}^2,
     \end{equation}
that is $f\in\dom(\Gamma_1)$. Thus, $f\in \dom (A_*) =
\dom(\Gamma) = \dom(\Gamma_0)\cap  \dom(\Gamma_1).$

Further, let us prove the converse statement assuming that ${\widetilde A}_0$ and
${\widetilde A}_1$ are transversal. Note that $A_{*0}={\widetilde A}_0$ if
$A_{*0}=A^*_{*0}$. Hence \eqref{III.1_03} yields $\dom({\widetilde A}_0)=\dom(A_{*0})\subset\dom(A_*)\subset\dom(\Gamma_1)$. On the
other hand, $\dom(A^*)=\dom({\widetilde A}_0)+\dom({\widetilde
A}_1)$ since ${\widetilde A}_0$ and ${\widetilde A}_1$ are
transversal. Thus $\Gamma_1$ admits an extensions on
$\gH_+=\dom(A^*)$, since $\dom({\widetilde
A}_1)\subset\dom(\Gamma_1)$. By $(v)$, $C_1<\infty$.
           \end{proof}

Next 
 we find a criterion for  a direct sum  $\Pi = \oplus_{n=1}^\infty\Pi_n$  to form a generalized boundary triplet.
      \begin{proposition}\label{prop_III.1_01}
Let $\Pi_n=\{\cH_n, \Gamma^{(n)}_0, \Gamma^{(n)}_1\}$ be a
boundary triplet  for $S^*_n$ and $M_n(\cdot)$ the corresponding
Weyl function,  $n\in \N$.   
Then the following conditions are
equivalent:
      \item $(i)$\   A direct sum  $\Pi =
\oplus_{n=1}^\infty\Pi_n =\{\cH,\Gamma_0,\Gamma_1\}$  is a
generalized boundary triplet for $A^*$,
    \item $(ii)$\
$\ran(\Gamma_0)=\cH  = \oplus_{n=1}^\infty {\cH}_n$,
 \item  $(iii)$\  $\sup_n\|M_n(\I)\|=:C_3<\infty$.
             \end{proposition}
\begin{proof}
$(i)\Rightarrow(ii)$ This implication  is immediate from
Definition \ref{def_II.2.1_generalized_bt}.

$(ii)\Rightarrow(i)$  
By Theorem \ref{th_III.1_01}$(i)$,  $\Pi$ is a boundary relation. Therefore,
by \cite[Lemma 4.10 $(iii)$]{DHMS06}, $A_{*0}$ is closed since
 $\ran(\Gamma_0) (=\cH)$ is closed.
On the other hand, by Theorem \ref{th_III.1_01}$(vi)$, $A_{*0}$ is essentially
self-adjoint. Thus $A_{*0} = (A_{*0})^*$ and the assumption $(iii)$ of
Definition \ref{def_II.2.1_generalized_bt} is verified.

$(ii)\Rightarrow (iii)$. Let $\ran(\Gamma_0)=\cH$. According to the implication $(ii)\Rightarrow (i)$, $\Pi$ is a generalized boundary triplet for
$A^*.$  Therefore, by \cite[Propostion 6.2]{DM95}, the
corresponding Weyl function $M$ takes values in $[\cH]$. By Theorem \ref{th_III.1_01} $(ii)$,
$M(z)=\oplus^{\infty}_{n=1}M_n(z)$  hence $M(\I)\in
[\cH]$ precisely when  $C_3=\sup_n\|M_n(\I)\| <\infty$.

$(iii)\Rightarrow (ii)$.
Let $\gamma_n$ be
the $\gamma$-field of the boundary triplet $\Pi_n$. Then \eqref{II.1.3_02} implies
     \begin{equation}\label{III.1_17}
\im
M_n(\I)=\bigl(M_n(\I)-M^*_n(\I)\bigr)/2\I=\gamma_n(\I)^*\gamma_n(\I),
\qquad n\in \N.
     \end{equation}
Since $\sup_n\|M_n(\I)\|= C_3<\infty$,  equality \eqref{III.1_17} yields
    \begin{equation}\label{III.1_18}
\sup_n\|\gamma_n(\I)\|^2  =   \sup_n\|\im M_n(\I)\| = C_3 <\infty.
    \end{equation}
Let $h=\oplus^{\infty}_{n=1} h_n\in\cH$. Then $f_n(\I):=
\gamma_n(\I)h_n\in\cN_\I(S_n^*)$ and, by \eqref{III.1_18},
     \begin{equation}\label{III.1_19}
\sum^{\infty}_{n=1}\|f_n(\I)\|^2  =
\sum^{\infty}_{n=1}\|\gamma_n(\I)h_n\|^2 \le
{C_3} \sum^{\infty}_{n=1}\|h_n \|^2 <\infty.
     \end{equation}
Hence $f(\I):=\oplus^{\infty}_{n=1}f_n(\I)\in\cN_\I(A^*)=
\oplus_{n=1}^\infty \cN_\I(S_n^*)$ and
 $\Gamma_0 f(\I) = \oplus_{n=1}^\infty\Gamma^{(n)}_0 f_n(\I) =
\oplus_{n=1}^\infty h_n = h$.
Thus $f(\I)\in\dom(\Gamma_0)$ and $\ran(\Gamma_0) = \cH.$ The proof is completed. 
        \end{proof}
      \begin{corollary}\label{cor_III.1_01}
Let $\Pi_n=\{\cH_n, \Gamma^{(n)}_0, \Gamma^{(n)}_1\}$ be a
boundary triplet  for $S^*_n$, $n\in \N$, and let $\Gamma_1$ be defined
by \eqref{III.1_02}. Then the following conditions are equivalent:
\item $(i)$\  $\sup_n\|M_n(\I)^{-1}\| = C_4 < \infty$,
\item $(ii)$\  $\ran(\Gamma_1) = \cH =  \oplus^{\infty}_{n=1} {\cH}_n$.
             \end{corollary}
\begin{proof}
Alongside the boundary triplet $\Pi_n$ we consider a triplet
$\widetilde{\Pi}_n =\{\cH_n, -\Gamma^{(n)}_1, \Gamma^{(n)}_0\},\
 n\in \N$.  The corresponding Weyl function is
$\widetilde{M}_n(\cdot) = - M_n(\cdot)^{-1},\ n\in \N$. To
complete the proof it remains to apply Proposition
\ref{prop_III.1_01}.
     \end{proof}

     \begin{remark}\label{rem_III.1_01}
By Theorem \ref{th_III.1_01} $(ii)$,
$\ker\bigl(\im M(z)\bigr)=\{0\}$, $z\in{\C}_+$, and hence
$M(\cdot)\in R^s(\cH)$. According to \eqref{III.1_05}, the
inequality $\sup_n\|M_n(\I)\|<\infty$ is equivalent to the
inclusion $M(\I)\in[\cH]$, that is $M(\cdot)\in R^s[\cH]$. Hence,
the implication $(iii)\Rightarrow(i)$ in Proposition
\ref{prop_III.1_01}  is immediate from \cite[Theorem 6.1]{DM95}.
However we prefer  a direct proof because of its simplicity.

Here  $R^s(\cH)$ and $R^s[\cH]$ are the Nevanlinna subclasses (definitions may be found in
\cite[Section 2.6]{DHMS06}).
   \end{remark}

Next we present  sufficient conditions for a direct sum $\Pi =
\oplus_{n=1}^{\infty}\Pi_n$ to be  a generalized boundary
triplet for $A^*$. These conditions are formulated only in terms of the mappings $\Gamma_j^{n}$.

\begin{proposition}\label{prop_III.1_02}
Assume the conditions of Theorem \ref{th_III.1_01} hold. Then:
\item $(i)$\ A direct sum  $\Pi = \oplus_{n=1}^{\infty}\Pi_n
=\{\cH,\Gamma_0,\Gamma_1\}$ of boundary triplets $\Pi_n$  is a
generalized boundary triplet for $A^*$ provided that
$C_1 =\sup_{n}\|\Gamma^{(n)}_1\|  <\infty$.
\item $(ii)$ If in addition ${\widetilde A}_0=\oplus_{n=1}^\infty
S_{n0}$ and ${\widetilde A}_1=\oplus_{n=1}^\infty S_{n1}$  are
transversal, then condition $C_1<\infty$ is necessary and
sufficient for $\Pi$ to be a generalized boundary triplet
for $A^*$.
  \end{proposition}
 \begin{proof}
$(i)$\  Condition $(G3)$ of Definition
\ref{def_II.2.1_generalized_bt} is immediate from Theorem
\ref{th_III.1_01} $(i)$. Moreover, by Theorem \ref{th_III.1_01}
$(vii)$, condition $C_1<\infty$ yields $A_{*0}= (A_{*0})^*$, hence
condition $(G2)$ of Definition \ref{def_II.2.1_generalized_bt}.
Let us check condition $(G1)$. Since $\gamma_n^*(\overline z) =
\Gamma_1^{(n)}(S_{n0} - z)^{-1}$ (see  \eqref{II.1.3_02B}), we get
that for any $n\in \N$
     \begin{eqnarray}\label{III.1_23}
\| \gamma_n(\overline z)^*f\|^2 =
\|\Gamma_1^{(n)}(S_{n0}-z)^{-1}f\|^2 \le
C_1^2\|(S_{n0}-z)^{-1}f\|_{\frak H_+}^2 \nonumber \\
= C_1^2( \|S_{n0}(S_{n0} - z)^{-1}f\|_{\frak H_+}^2 +
\|(S_{n0}-z)^{-1}f\|_{\frak H_+}^2) \le 2C^2_1(1 + (|z|^2 +1)/|\im
z|^2), 
               \end{eqnarray}
and hence $\| \gamma_n(\pm \I)\| = \| \gamma_n^*(\pm \I)\|\le
C_1\sqrt 6,\ n\in \N.$ Since $M_n(z) = \Gamma_1^{(n)} \gamma_n(
z)$ (see \eqref{II.1.3_01}), we have
$$
\|M_n(\I)h\| \le \|\Gamma_1^{(n)}\|\cdot \|\gamma_n(\I)h\|_{\frak
H_+} \le C_1 \sqrt 2 \|\gamma_n(\I)h\| \le  C_1^2 \sqrt {12},
\quad n\in \N.
$$
Hence, by Proposition \ref{prop_III.1_01},
$\ran(\Gamma_0) = \cH$.

$(ii)$\ 
Follows from Theorem \ref{th_III.1_01} $(vii)$.
       \end{proof}
\begin{corollary}\label{cor_III.1_02}
Assume  the conditions of Proposition  \ref{prop_III.1_01}. 
Then:
\item $(i)$ A direct sum  $\widetilde{\Pi} =
\oplus_{n=1}^{\infty}\widetilde{\Pi}_n$ of boundary triplets
$\widetilde{\Pi}_n =
\{\cH,\widetilde{\Gamma}_0^{(n)},\widetilde{\Gamma}_1^{(n)}\} = \{\cH_n,
-\Gamma^{(n)}_1, \Gamma^{(n)}_0\}$  is a generalized boundary
triplet for $A^*$ whenever
$C_0  = \sup_{n}\|\Gamma^{(n)}_0\|  <\infty$.
\item $(ii)$ If in addition ${\widetilde
A}_0=\oplus_{n=1}^{\infty} S_{n0}$ and ${\widetilde
A}_1=\oplus_{n=1}^{\infty} S_{n1}$  are transversal, then
condition $C_0<\infty$ is necessary and sufficient for
$\widetilde{\Pi}$ to be a generalized boundary triplet for
$A^*$.
  \end{corollary}

\subsection{When direct sum of boundary triplets is a boundary triplet?}\label{subsec_III.2}

\textbf{1. General case.}\\ 
As it was already mentioned, the direct sum
$\Pi=\oplus_{n=1}^\infty\Pi_n$ is not a boundary triplet
without additional restrictions (cf. Theorem \ref{th_III.1_01}).
We start with the following result.
      \begin{proposition}\label{prop3.8}
Assume the conditions of Theorem \ref{th_III.1_01}. Then the direct sum
$\Pi = \oplus_{n=1}^{\infty}\Pi_n$
is an ordinary  boundary triplet for $A^*$ if and only if
   \begin{equation}\label{3.21}
\max\{C_0,\  C_1\} <\infty, \qquad C_j  =
\sup_{n\in\N}\|\Gamma^{(n)}_j\|.
   \end{equation}
        \end{proposition}
\begin{proof}
Necessity is immediate from \eqref{III.1_02} and Definition \ref{def_ordinary_bt}.

Sufficiency.   Consider  $\gH^2 :=\gH\oplus\gH$ and $\cH^2:=
\cH\oplus\cH$ as Krein spaces with the fundamental symmetries
$J_{\gH}=i
\begin{pmatrix}
0&-I_{\gH}\\
I_{\gH}&0
\end{pmatrix}$ and $J_{\cH}=i
\begin{pmatrix}
0&-I_{\cH}\\
I_{\cH}&0
\end{pmatrix}$,
respectively.  Now  identity \eqref{III.1_06} can be rewritten as
     \begin{equation}
(J_{\gH}{\hat f},{\hat g})_{\gH^2} = (J_{\cH}\Gamma \hat f,\Gamma
\hat g)_{\cH^2},
     \end{equation}
where ${\hat f}:=\{f,A^* f\},\ {\hat g}:=\{g,A^* g\}$ and $\Gamma
\hat f := \Gamma  f$. This means that  $\Gamma: \gH^2\to \cH^2$ is
an isometry from the Krein space $\{\gH^2, J_\gH\}$ to the Krein
space $\{\cH^2, J_\cH\}$. By Theorem \ref{th_III.1_01} $(v)$,
$\dom(\Gamma) = \gr(A^*)$, the graph of $A^*$. Since $\dom(\Gamma)$
is closed in $\gH^2$, $\ran(\Gamma)$ is closed too (see
\cite[Proposition 2.3]{DHMS06}). On the other hand, by Theorem
\ref{th_III.1_01} $(iii)$, $\ran(\Gamma)$ is dense in $\cH^2$ and hence
$\ran(\Gamma)=\cH^2$.
         \end{proof}

\begin{remark}
Proposition  \ref{prop3.8} shows that condition
\eqref{3.21} is sufficient (but not necessary!)  for
transversality of the extensions  $A_{*0}$ and $A_{*1}$ defined by
\eqref{III.1_03}.  This fact complements Theorem
\ref{th_III.1_01}(vii). Moreover, it shows that in the case of a
special boundary relation $\Pi = \oplus_{n=1}^{\infty}\Pi_n
$, condition  $(d)$ after Definition
\ref{def_III.1_01} is sufficient for $\Pi = \oplus_{n=1}^{\infty}\Pi_n$ to be  an ordinary
boundary triplet. Besides, $(d)$ and $(c)$ are equivalent and
yield the previous conditions $(a), (b)$.
      \end{remark}

Now we are ready to state  the main results of this section.
    \begin{theorem}\label{th_criterion(bt)}
Let  $\Pi_n=\{\cH_n, \Gamma_0^{(n)}, \Gamma_1^{(n)} \}$ be a
boundary triplet for $S_{n}^*$ and  $M_n(\cdot)$  the
corresponding Weyl function, $n\in \N $. A direct sum
$\Pi=\oplus_{n=1}^{\infty}\Pi_n$ forms  an ordinary boundary
triplet for the operator $A^* =\oplus_{n=1}^{\infty}S_n^*$ if and only if
        \begin{equation}\label{WF_criterion}
C_3= \sup_n\|M_n(\I)\|_{\mathfrak{H}_n}  < \infty
\quad\text{and}\quad   C_4= \sup_n\|(\im
M_n(\I))^{-1}\|_{\mathfrak{H}_n} < \infty.
   \end{equation}
         \end{theorem}
\begin{proof}


By  Proposition  \ref{prop_III.1_01}, the first
inequality in \eqref{WF_criterion} is equivalent to the fact
that $\Pi = \oplus_{n=1}^{\infty}\Pi_n$ is a generalized boundary
triplet for the operator $A^*$. By Theorem \ref{th_III.1_01},
$(ii)$ the corresponding (generalized) Weyl function is
$M(\cdot)=\oplus^{\infty}_{n=1}M_n(\cdot).$ Therefore, the second
inequality in \eqref{WF_criterion} is equivalent to $C_4= \|(\im
M(\I))^{-1}\|_{\mathfrak{H}} < \infty,$ that is to the condition
$0\in \rho(\im M(i)).$ To complete the proof it remains to apply
Proposition \ref{prop2.10}.
        \end{proof}

 Theorem \ref{th_criterion(bt)} makes it possible to
construct an ordinary boundary triplet starting
with an arbitrary boundary relation $\Pi=\oplus_{n=1}^\infty \Pi_n$.
    \begin{theorem}[\cite{MalNei_08}]\label{th_III.2.1_02}
Let $S_n$ be a symmetric operator in $\mathfrak{H}_n$ with
deficiency indices $\mathrm{n}_{\pm}(S_k) = \mathrm{n}_n\le \infty$ and $S_{n0}=
S_{n0}^*\in \Ext S_n, \  n\in \N.$ Then for any $n\in \N$ there
exists a boundary triplet $\Pi_n=\{\mathcal{H}_n, \Gamma_0^{(n)},
\Gamma_1^{(n)} \}$  for $S_{n}^*$ such that $\ker\Gamma_0^{(n)}=
\dom(S_{n0})$ and $\Pi=\oplus_{n=1}^{\infty}\Pi_n$ forms an
ordinary  boundary triplet for $A^* =\oplus_{n=1}^{\infty}S_n^*$
satisfying $\ker\Gamma_0 = \dom(\widetilde{A}_{0}):=
\oplus_{n=1}^{\infty}S_{n0}.$
     \end{theorem}
\begin{proof}
By \cite[Chapter III.1.4]{Gor84}, there exists a boundary triplet
$\widetilde{\Pi}_n=\{\mathcal{H}_n, \widetilde{\Gamma}_0^{(n)},
\widetilde{\Gamma}_1^{(n)} \}$ for $S_n^*$ such that
$\dom(S_{n0})=S_n^*\upharpoonright\ker \widetilde{\Gamma}_0^{(n)},
n\in \N$. Let $\widetilde{M}_n(\cdot)$ be  the corresponding Weyl
function. Denote $Q_n :=\re \widetilde{M}_n(\I)$ and choose a
factorization of ${\im} \widetilde M_n(\I)$,  $R_n^*R_n := {\im}
\widetilde M_n(\I),$ such that $R_k\in [\cH_k]$ and $0\in
\rho(R_k).$ Then we  define the mappings
$\Gamma_j^{(n)}:\dom(S_n^*)\to \cH_n$ as follows
   \begin{equation}\label{3.24}
\Gamma_0^{(n)} := R_n\widetilde{\Gamma}_0^{(n)},\qquad
\Gamma_1^{(n)} :=
(R_n^*)^{-1}(\widetilde{\Gamma}_1^{(n)}-Q_n\widetilde{\Gamma}_0^{(n)}),
\qquad n\in \N.
   \end{equation}
 It is easy  to check that $\Gamma_j^{(n)}$ are well defined
 and $\Pi_n=\{\mathcal{H}_n, \Gamma_0^{(n)}, \Gamma_1^{(n)} \}$
forms  a boundary triplet for $S_{n}^*$. Moreover, the  Weyl
function $M_n(\cdot)$ corresponding to $\Pi_n$ satisfies
$M_n(\I)=\I I_{\cH_n},\  n\in \N.$
%
%
 Hence, by Theorem \ref{th_criterion(bt)}, a triplet
$\Pi=\oplus_{n=1}^\infty \Pi_n$ forms a boundary triplet for
$A^*$. The required property $\ker\Gamma_0 = \ker\wt\Gamma_0 =
\dom(A_{0}):= \oplus_{n=1}^{\infty}S_{n0}$ is immediate from
\eqref{3.24}.
           \end{proof}
\begin{remark}
Note that the regularization \eqref{3.24}  of the direct sum
$\wt\Pi=\oplus_{n=1}^\infty \wt\Pi_n = \{\mathcal{H},
\wt{\Gamma}_0, \wt{\Gamma}_1 \}$  has been proposed in \cite[Theorem
5.3]{MalNei_08}.  We emphasize however that condition
\eqref{WF_criterion} is more flexible than the condition
$M_n(\I)=\I I_{\cH_n}$, $n\in \N$, given in \cite[Theorem
5.3]{MalNei_08}. The latter is  very important in applications (cf.
Remark  \ref{rem_III.2.2_02} below).
\end{remark}

\textbf{2. The case of operators with common regular real point.}\\ 
Assume the operator $A=  \oplus^{\infty}_{n=1}S_n$ has a
regular real point, i. e., there exists $a={\overline
a}\in{\hat\rho}(A)$. This is equivalent to the existence of
$\varepsilon> 0$ such that
   \begin{equation}\label{III.2.2_01}
(a-\varepsilon, a + \varepsilon) \subset \cap^{\infty}_{n=1}
{\hat\rho}(S_n).
       \end{equation}
In particular, \eqref{III.2.2_01} holds  whenever the operators $S_n$
are nonnegative, $S_n\geq 0$.
Assuming condition \eqref{III.2.2_01} to be satisfied, we can simplify
conditions \eqref{WF_criterion} of Theorem
\ref{th_criterion(bt)} as follows.

       \begin{theorem}\label{th_III.2.2_01}
Let $\{S_n\}_{n= 1}^\infty$ be a sequence of symmetric operators
satisfying \eqref{III.2.2_01}. Let also $\Pi_n=\{\cH_n,
\Gamma^{(n)}_0, \Gamma^{(n)}_1\}$ be a boundary triplet for
$S^*_n$ such that
$(a-\varepsilon,a+\varepsilon)\subset\rho(S_{n0})$ and
$M_n(\cdot)$ the corresponding Weyl function. Then $\Pi=
\oplus^{\infty}_{n=1} \Pi_n$ is a boundary triplet for $A^*=
\oplus^{\infty}_{n=1}S_n^*$ if and only if
    \begin{equation}\label{III.2.2_02}
C_5 := \sup_{n\in\N}\|M_n(a)\|<\infty \qquad \text{and} \qquad  C_6 :=
\sup_{n\in\N}\|\bigl(M'_n(a)\bigr)^{-1}\|<\infty,
     \end{equation}
where  $M'_n(a):=({dM_n}(z)/{dz})|_{z=a}$.
        \end{theorem}
       \begin{proof}
Necessity is obvious. Indeed, if $\Pi=\oplus_{n=1}^\infty\Pi_n$ is a boundary triplet, then the corresponding Weyl function $M(\cdot)$ is defined by \eqref{III.1_05}. Moreover, $M(\cdot)$ is an $R_{[\cH]}$-function analytic at $z=a$ and hence $M(a)\in [\cH]$. Furthermore, it satisfies $0\in \rho (M'(a))$ and thus \eqref{III.2.2_02} is fulfilled.

Sufficiency. We deduce the proof from Theorem  \ref{th_criterion(bt)}. Namely, we will show that conditions \eqref{WF_criterion} of Theorem \ref{th_criterion(bt)} are implied by the corresponding conditions
in \eqref{III.2.2_02}.

  First note that
$M(\cdot):=\oplus^{\infty}_{n=1}M_n(\cdot)$ is a ${\mathcal
C}(\cH)$-valued Nevanlinna function since for any $z\in\C_+$ the
operator $M(z)$ is closed. Further, $M_n(\cdot)$ is regular on
$(a-\varepsilon, a+\varepsilon)$ since $(a-\varepsilon,
a+\varepsilon)\subset\rho(S_{n 0})$. Due to condition \eqref{III.2.2_01},  $M(\cdot)$ is also holomorphic on
$(a-\varepsilon, a+\varepsilon)$ in the sense of Kato \cite{Kato66}, that is $\bigl(M(z)-\I\bigr)^{-1}$ is bounded and
holomorphic at $z_0=a$, as well as at $z\in {\C}_+\cup{\C}_-\cup (a-\varepsilon, a+\varepsilon)$
(see \cite[Theorem 7.1.3]{Kato66}).
Moreover, due to the first condition in \eqref{III.2.2_02}, $M(\cdot)$ is bounded at $z=a$,
$M(a)\in[\cH]$. By \cite[Section 7.1.2]{Kato66}, $M(z)\in[\cH]$
for $|z-a|$ small enough (see also \cite[Theorem 4.2.23(b)]{Kato66}.
In turn, the latter yields $M(z)\in[\cH]$ for any $z\in{\C}_+$ (see \cite{DHMS06}). In particular,
$M(\I)\in[\cH]$ and the first inequality in \eqref{WF_criterion} is verified.

Further, by \eqref{II.1.3_02},
    \begin{equation}\label{III.2.2_04}
M'_n(a) = ({dM_n}(z)/{dz})|_{z=a} =
\gamma^*_n(a)\gamma_n(a),\qquad n\in \N.
     \end{equation}
According to  \eqref{II.1.3_02'}, $\gamma_n(\I) = [I-(a-\I)(S_{n 0} -
\I)^{-1}]\gamma_n(a)$. Hence
   \begin{equation}\label{III.2.2_06}
\gamma_n^*(\I)\gamma_n(\I)=\gamma^{*}_n(a)[I-(a+\I)(S_{n0}+\I)^{-1}][I-(a-\I)(S_{n 0} - \I)^{-1}]\gamma_n(a).
      \end{equation}
Noting that $\bigl(I-(a-\I)(S_{n0}-\I)^{-1}\bigr)^{-1}=I+(a-\I)(S_{n0}-a)^{-1}$, we get
\[
\inf_{f\in \cH_n}(\gamma_n^*(\I)\gamma_n(\I)f,f)
\geq\|I+(a-\I)(S_{n0}-a)^{-1}\|^{-2}_{\mathfrak{H}_n}\inf_{f\in
\cH_n}(\gamma_n^*(a)\gamma_n(a)f,f).
\]
Since $(a -\varepsilon,  a + \varepsilon)\subset\rho(S_{n 0})$, we
have $\|I+(a-\I)(S_{n0}-a)^{-1}\|\le 1+\frac{\sqrt{1+a^2}}{\varepsilon}=:C$. Combining these inequalities with
\eqref{III.2.2_06} and \eqref{III.1_17}, we obtain
 \[
 \|(\im M_n(\I))^{-1}\|_{\cH_n}\leq 
 C^2\|(M_n'(a))^{-1}\|_{\cH_n},
 \]
and the second inequality in \eqref{WF_criterion} is verified.
      \end{proof}
For operators $A = \oplus_{n=1}^\infty S_n$ satisfying
\eqref{III.2.2_01} we complete Theorem \ref{th_III.2.2_01}  by
presenting a regularization procedure for $\Pi =
\oplus_{n=1}^\infty\Pi_n$ leading to a boundary triplet.
  \begin{corollary}\label{cor_III.2.2_01}
Let $\{S_n\}^{\infty}_{n=1}$ be a sequence of symmetric operators
satisfying \eqref{III.2.2_01}. Let also ${\widetilde \Pi}_n =
\{\cH_n, {\widetilde\Gamma}^{(n)}_0, {\widetilde\Gamma}^{(n)}_1\}$
be a boundary triplet for $S^*_n$ such that $(a-\varepsilon,
a+\varepsilon)\subset\rho(S_{n0})$, $S_{n0}=
S_n^*\lceil\ker({\widetilde\Gamma}^{(n)}_0)$, and ${\widetilde
M}_n(\cdot)$ the corresponding Weyl function,  $n\in{\N}$. Then:
\item $(i)$\  The operator ${\widetilde M}'_n(a)$ is positively
definite, $n\in{\N}$.
 \item $(ii)$\  For any factorization ${\widetilde M}'_n(a)= R_n^*R_n$, where $R_n\in
[\cH_n]$ and $0\in\rho(R_n)$, a triplet
     \begin{equation}\label{III.2.2_08}
\Pi_n=\{\cH_n,\Gamma^{(n)}_0,\Gamma^{(n)}_1\}\quad
\text{with}\quad \Gamma^{(n)}_0 := R_n{\widetilde\Gamma}^{(n)}_0,\qquad
 \Gamma^{(n)}_1:= (R^{-1}_n)^*\bigl({\widetilde\Gamma}^{(n)}_1 -
 \wt M_n(a){\widetilde\Gamma}^{(n)}_0\bigr),
     \end{equation}
is a boundary triplet for $S^*_n$.
\item $(iii)$ A direct sum $\Pi=\oplus_{n=1}^\infty\Pi_n$ forms a
boundary triplet for $A^*$.
    \end{corollary}
        \begin{proof}
$(i)$\ Let  ${\widetilde\gamma}_n$ be  the $\gamma$-field
corresponding to the triplet ${\widetilde \Pi}_n = \{\cH_n,
{\widetilde\Gamma}^{(n)}_0, {\widetilde\Gamma}^{(n)}_1\}$. The
functions  ${\widetilde M}_n(\cdot)$ and
${\widetilde\gamma}_n(\cdot)$ are regular within $(a-\varepsilon,
a+\varepsilon)$ for every $n\in{\N}$ since $(a-\varepsilon,
a+\varepsilon)\subset\rho(S_{n0})$. By \eqref{III.2.2_04},
$\widetilde M'_n(a)>0$ and $0\in\rho\bigl(\widetilde
M'_n(a)\bigr)$ since $\gamma_n(a)$ isomorphically maps $\cH_n$
onto $\cN_a$.

$(ii)$\ By $(i)$, $\widetilde M'_n(a)$ admits a factorization
$\widetilde M'_n(a) = R_n^*R_n$, where $R_n\in [\cH]$ and
$0\in\rho(R_n)$. Therefore, the mappings $\Gamma^{(n)}_0$ and
$\Gamma^{(n)}_1$ are defined correctly and $\Pi_n$ is a boundary
triplet for $S^*_n$.

$(iii)$ Let $M_n(\cdot)$ be the Weyl function corresponding to the
triplet $\Pi_n$. 
It follows from \eqref{III.2.2_08} and the definition of the Weyl
function that
     \begin{equation}\label{III.2.2_09}
M_n(z) = (R^{-1}_n)^*[{\widetilde M}_n(z)-{\widetilde
M}_n(a)]R^{-1}_n, \qquad  n\in{\N}.
   \end{equation}
Hence $M_n(a)=0$ and $M'_n(a) = (R^{-1}_n)^*{\widetilde M}'_n(a)
R^{-1}_n =  I_{\cH_n}$, $n\in{\N}$. Thus,  both conditions in
\eqref{WF_criterion}  are satisfied and, by Theorem
\ref{th_III.2.2_01}, $\Pi=\oplus_{n=1}^\infty\Pi_n$ forms a
boundary triplet for $A^*$.
        \end{proof}
\begin{corollary}\label{cor_III.2.2_02}
Let $\{S_n\}^{\infty}_{n=1}$ be a sequence of symmetric operators
satisfying \eqref{III.2.2_01}.
Let also ${\widetilde \Pi}_n = \{\cH_n,
{\widetilde\Gamma}^{(n)}_0, {\widetilde\Gamma}^{(n)}_1\}$ be a
boundary triplet for $S^*_n$ such that $(a-\varepsilon,
a+\varepsilon)\subset\rho(S_{n0})$, $S_{n0}=
S_n^*\lceil\ker({\widetilde\Gamma}^{(n)}_0)$, and ${\widetilde
M}_n(\cdot)$ the corresponding Weyl function. If the operators
$R_n\in[\cH_n]$  satisfy
       \begin{equation}\label{III.2.2_11}
R^{-1}_n\in[\cH_n]
\qquad \text{and}\qquad \sup_n \|R_n({\widetilde
M}'_n(a))^{-1}R_n^{*}\|<\infty, \quad n\in\N,
       \end{equation}
 then the direct sum $\Pi=\oplus_{n=1}^\infty \Pi_n$ of boundary triplets
 \eqref{III.2.2_08} forms a boundary triplet for $A^* =\oplus_{n=1}^\infty S_n^*$.
       \end{corollary}
      \begin{proof}
Since the Weyl function $M_n(\cdot)$ corresponding to $\Pi_n$ is
given by  \eqref{III.2.2_09}, both conditions \eqref{III.2.2_02}
are  immediate from \eqref{III.2.2_11}. It remains to apply
Theorem \ref{th_III.2.2_01}.
             \end{proof}
  \begin{remark}\label{rem_III.2.2_02}
Corollary \ref{cor_III.2.2_02} is more useful in applications
than Corollary \ref{cor_III.2.2_01}. The reason is that it is more
convenient and easier to select a suitable sequence
$\{R_n\}_{n=1}^\infty$ satisfying \eqref{III.2.2_11} than to find
the operators $({M'_n(a)})^{1/2}$. For instance, to construct
boundary triplets in Theorems \ref{th_bt_2} and \ref{th_bt_1}, we
select $R_n$ being diagonal matrices although $M'_n(a)$,  hence
$(M'_n(a))^{1/2},$ are not diagonal.
         \end{remark}

\subsection{Direct sums of self-similar boundary triplets}\label{subsec_III.3}

In this subsection, we apply Theorem \ref{th_criterion(bt)} to the
special case of symmetric operators $S_n$ that are pairwise
unitarily equivalent  up to multiplicative constants. More
precisely, let $S_1$ be a symmetric operator in $\mathfrak{H}_1,$
$\mathrm{n}_{\pm}(S_1) = n\le \infty$. We assume
that for any $n\in \N$ there exists a unitary operator $U_n$ from
$\mathfrak{H}_n$ onto $\mathfrak{H}_1$ and a constant $\gd_n>0$ such
that (to be precise we set $U_1:=I_{\cH_1}$ and $d_1:=1$)
\begin{equation}\label{III.3_01}
S_n:= d_n^{-2} U^{-1}_n
S_1U_n.
      \end{equation}
First we  suppose that
     \begin{equation}\label{III.3_02}
0<d_*:=\inf_{n\in\N} d_n\leq \sup_{n\in\N} d_n =:d^*<\infty
       \end{equation}
and  reprove one result of Kochubei (cf. \cite[Theorem 3]{Koc_79},
\cite[Lemma 1]{Koc_89}) for this case.
       \begin{lemma}[\cite{Koc_89}]\label{lem_koch}
Let $S_n$ be as above, let  $\Pi_1=\{\cH_1, \Gamma_0^{(1)},
\Gamma_1^{(1)}\}$ be a boundary triplet for $S_1^*$, and
$A=\oplus_{n=1}^{\infty} S_n$. Assume in addition that condition
\eqref{III.3_02} holds. Then:
      \item $(i)$ \ For any
$\alpha\in\R$, a triplet  $\Pi_n := \{\mathcal{H}_n,
\Gamma_0^{(n)}, \Gamma_1^{(n)}\},$ where
      \begin{equation}\label{III.3_03}
\mathcal{H}_n:=\cH_1,\quad\Gamma_0^{(n)}:=d_n^{\alpha-2}\Gamma_0^{(1)}
U_n,\qquad \Gamma_1^{(n)}:= d_n^{-\alpha}\Gamma_1^{(1)}
U_n, \quad n\in\N,
        \end{equation}
forms a boundary triplet for the operator $S_n^*$.
     \item $(ii)$
\ Moreover, $\Pi=\oplus_{n=1}^\infty \Pi_n$ is an (ordinary)
boundary triplet for the operator $A^*.$
       \end{lemma}
            \begin{proof}
$(i)$ Straightforward.

$(ii)$ Let $M_n(\cdot)$ be the Weyl function 
corresponding to the triplet    
$\Pi_n=\{\cH, \Gamma_0^{(n)}, \Gamma_1^{(n)}\}$, $n\in\N$.  It follows from  \eqref{III.3_03} that the Weyl functions
$M_n$ and $M_1$ are connected by
      \begin{equation}\label{3.48A}
M_n(z)=d_n^{2-2\alpha} M_1(d_n^{2}z), \qquad z\in\C_\pm,\qquad n\geq 2.
      \end{equation}
Hence
  \begin{equation}\label{3.48B}
 \|M_n(\I)\|=d_n^{2-2\alpha} \|M_1(\I d_n^{2})\|,\qquad
\|(\im M_n(\I))^{-1}\| = d_n^{2\alpha-2} \|(\im
M_1(\I d_n^{2}))^{-1}\|.
 \end{equation}
 Combining \eqref{3.48B} with \eqref{III.3_02}, we obtain that $\{M_n\}_{n=1}^\infty$ satisfies
\eqref{WF_criterion} since $M_1$ is continuous on $[\I(d_*)^{2}, \I(d^*)^{2}]\subset \C_+$. Theorem \ref{th_criterion(bt)} completes the proof.
         \end{proof}

The following results demonstrate  importance  of  both
inequalities in \eqref{III.3_02} for the direct sum
$\Pi=\oplus_{n=1}^\infty \Pi_n$ to be an (ordinary) boundary
triplet  for $A^*$.
 \begin{lemma}\label{lem_III.3_02}
Let  $S_1$ be a closed densely defined symmetric operator in $\mathfrak{H}_1$ with
$\mathrm{n}_{\pm}(S_1) = n< \infty$, let $\Pi_1=\{\mathcal{H}_1,
\Gamma_0^{(1)}, \Gamma_1^{(1)}\}$ be a boundary triplet for
$S_{1}^*$ and $M_1(\cdot)$  the corresponding Weyl function.
Let also  $S_n$, $n\in\N$, be defined by (\ref{III.3_01}) and suppose that $\{d_n\}_{n=1}^\infty$ satisfies
$\gd_*=0$ and $\gd^*<\infty$. Then:
%
\item $(i)$ A direct sum $\Pi=\oplus_{n=1}^\infty\Pi_n$ of triplets $\Pi_n=\{\mathcal{H}_n,
\Gamma_0^{(n)}, \Gamma_1^{(n)}\}$,
where
\begin{equation}\label{III.3_05}
 \mathcal{H}_n=\cH_1,\qquad \Gamma_0^{(n)}=\Gamma_0^{(1)}U_n,\qquad
 \Gamma_1^{(n)}= d_n^{-2} \Gamma_1^{(1)} U_n,
      \end{equation}
 forms an ordinary boundary triplet for the  operator $A^* =
\oplus_{n=1}^{\infty}S_n^*$
if and only if
      \begin{equation}\label{III.3_06}
C_+:= -\lim_{y\downarrow 0}\frac{M_1(\I y)}{\I y}\in[\cH_1]\qquad \mathrm{and}\qquad 0\in\rho( C_+).
      \end{equation}
\item $(ii)$ A direct sum $\Pi=\oplus_{n=1}^\infty\Pi_n$ of triplets $\Pi_n=\{\mathcal{H}_n,
\Gamma_0^{(n)}, \Gamma_1^{(n)}\}$,
where
\begin{equation}\label{III.3_07}
 \mathcal{H}=\cH_1,\qquad \Gamma_0^{(n)}= d_n^{-2}\Gamma_0^{(1)}U_n,\qquad
 \Gamma_1^{(n)}= \Gamma_1^{(1)} U_n,
       \end{equation}
forms an ordinary boundary triplet for  $A^* =
\oplus_{n=1}^{\infty}S_n^*$  if and only if
\begin{equation}\label{III.3_08}
C_-:=-\lim_{y\downarrow 0}\I y M_1(\I y)\in [\cH_1] \qquad \mathrm{and}\qquad 0\in\rho(C_-).
    \end{equation}
            \end{lemma}
      \begin{proof}
$(i)$ By (\ref{3.48A}), we get $ M_n(\I)=d_n^{-2} M_1(\I
d_n^{2})$. Since $d_*=0$, by Proposition \ref{prop_III.1_01},
$\Pi$ is a generalized boundary triplet for $A^*$ if and only if
$C_+\in [\cH_1].$  Moreover, by Theorem \ref{th_criterion(bt)},
$\Pi$ is an ordinary boundary triplet precisely if in addition
$0\in\rho(C_+)$.

$(ii)$ The proof is  similar  to that of $(i)$ if one notices that
$M_n(\I)=d_n^{2} M_1(\I d_n^{2})$.
      \end{proof}
      \begin{remark}\label{rem_III.3_05}
      Let $\Sigma_{M_1}(\cdot)$ be the  spectral measure of  $M_1(\cdot)$
      (see Section \ref{sss_II.1.3_weylsec}). Then the operators $C_+$ and $C_-$
      can  easily be expressed in terms of $\Sigma_{M_1}(\cdot)$.
      Namely, condition (\ref{III.3_06}) means that the limit $M_1(0):=M_1(+\I0)$ exists, moreover,
      $M_1(0)=0$,  and the following integral converges
      \[
      C_+=\int_\R\frac{d\Sigma_{M_1}(t)}{t^2}\in [\cH_1].
      \]
      Besides, we note that $C_-=\Sigma_{M_1}(\{0\}).$
      \end{remark}
         \begin{corollary}\label{cor_III.3_05}
Let $S_n$ be as in Lemma \ref{lem_III.3_02}, let
$\Pi_1=\{\cH_1,\Gamma^{(1)}_0,\Gamma_1^{(1)}\}$ be a boundary
triplet for $S^*_1$ and  ${S_{1}}_0:= S_1^*\lceil
\ker(\Gamma^{(1)}_0)$. Assume that $\gd_*=0$ and $\gd^*<\infty$. Assume also that $S_1$ is a simple symmetric operator.
The direct sum $\Pi=\oplus_{n=1}^\infty\Pi_n$ of boundary
triplets defined by (\ref{III.3_07}) (by (\ref{III.3_05})) is an
ordinary boundary triplet for $A^*$ if and only if
    \begin{equation}\label{3.60}
\dim(\ker S_{10})=\mathrm{n}_{\pm}(S_1), \qquad \bigl(\mathrm{respectively},\
\ \dim(\ker S_{11})=\mathrm{n}_{\pm}(S_1)\bigr).
      \end{equation}
    \end{corollary}
      \begin{proof}
 Let us prove the first equality in (\ref{3.60})  assuming that
the direct sum $\Pi=\oplus_{n=1}^\infty\Pi_n$ of boundary triplets
(\ref{III.3_07}) forms a boundary triplet. By Remark
\ref{rem_III.3_05}, $C_-=\Sigma_{M_1}(\{0\})$ where
$\Sigma_{M_1}(\cdot)$ is a nonorthogonal spectral measure of $M_{1}(\cdot)$.  The latter implies
\[
\dim(\ker S_{10})=\Rank(\Sigma_{M_1}(\{0\}))=\Rank C_-.
\]
 Since $\Pi$ is an
ordinary boundary triplet for $A^*$, Lemma \ref{lem_III.3_02}(ii)
yields $0\in{\rho}(C_-),$ that is, $C_-$ is of maximal rank.
Combining these relations, we get $\dim(\ker S_{10}) = \Rank (C_-)=\dim
\cH_1=\mathrm{n}_\pm(S_1)$.

To prove sufficiency, let us note that $C_-=\Sigma_{M_1}(\{0\})$. Since $\mathrm{n}_\pm(S_1)=\dim \cH_1<\infty$, we obviously get that relations (\ref{3.60}) and (\ref{III.3_08}) are equivalent.
%
         \end{proof}
%
%
%
%

We complete this subsection by considering the situation when
$\gd^*=\infty$.
\begin{lemma}
Let $\gd_*>0$ and $\gd^* =\infty$.  Then:
       \item $(i)$ The
direct sum $\Pi=\oplus_{n=1}^\infty\Pi_n$ of triplets defined by
\eqref{III.3_05} is a generalized boundary triplet for $A^*$, but
not an ordinary boundary triplet for $A^*$,
        \item $(ii)$
$\Pi=\oplus_{n=1}^\infty\Pi_n$  is not a generalized boundary
triplet for $A^*$ if $\Pi_n$ is defined by \eqref{III.3_07}.
    \end{lemma}
\begin{proof}
$(i)$ 
Since $S_1$ is
densely defined, the Weyl function $M_1(\cdot)$ corresponding to
the triplet $\Pi_1$ satisfies (cf. \eqref{WF_intrepr})
    \begin{equation}\label{3.56A}
 s-\lim_{y\uparrow\infty}M_1(\I y)/y=0
        \end{equation}
Let $\Pi_n, \ n\in\N$, be the boundary
triplet for $S_n^*$ defined by \eqref{III.3_05} and $M_n(\cdot)$ the
corresponding Weyl function.
Setting in  \eqref{III.3_03} and \eqref{3.48A}
$\alpha =2$ and  combining these relations with \eqref{III.3_05},
we get $M_n(z)=d_n^{-2} M_1(d_n^{2}z).$
Combining these relations  with \eqref{3.56A}, we obtain
    \begin{equation}\label{3.43}
\sup_n\|M_n(\I)\|=\sup_n d^{-2}_n\|M_1(\I d^{2}_n)\|<\infty.
       \end{equation}
By Proposition \ref{prop_III.1_01}, $\Pi=\oplus_{n=1}^{\infty} \Pi_n$
forms  a generalized boundary triplet for $A^*.$

Further, the above relations yield $\im M_n(z)= d_n^{-2} \im
M_1(d_n^{2}z).$  Hence and from \eqref{3.56A} we get
   \begin{equation}
\sup_n\|\bigl(\im M_n(\I)\bigr)^{-1}\|=\sup_n d^{2}_n\|\bigl(\im
M_1(\I d^{2}_n)\bigr)^{-1}\| = \infty.
     \end{equation}
By Theorem \ref{th_criterion(bt)},  $\Pi=\oplus_{n=1}^{\infty} \Pi_n$
is not  an ordinary  boundary triplet for $A^*.$

$(ii)$  Since $S_1$ is densely defined, the
Weyl function $M_1(\cdot)$  satisfies (cf.\eqref{3.56A})
%
    \begin{equation}\label{3.56AB}
 s-\lim_{y\uparrow\infty}y^{-1} M_1(\I y)^{-1}=0.
        \end{equation}
Let $\Pi_n, \ n\in\N$, be a boundary
triplet for $S_n^*$ defined by \eqref{III.3_07} and $M_n(\cdot)$ the
corresponding Weyl function. It follows from   \eqref{III.3_03} and \eqref{3.48A} (with $\alpha
=0$) that  $M_n(z)=d_n^{2} M_1(d_n^{2}z), \
n\ge 2$.  Hence $\sup_{n}\|M_n(\I)\|=\sup \gd^{2}_n\|M_1(\I
\gd^{2}_n)\|=\infty$.   By Proposition \ref{prop_III.1_01},
$\Pi=\oplus_{n=1}^{\infty} \Pi_n$ is not   a generalized
boundary triplet for $A^*$.
          \end{proof}

%
%
%
\section{Boundary triplets for the operator $\rH_{\min}^*$.}\label{Subsec_IV.1_bt_interactions}
%
%
%
%

In what follows we assume that $\cI=[0,b)\subseteq\R_+, \
0<b\leq +\infty$, is either a bounded interval or
positive semi axis, $X=\{x_n\}_{n=0}^\infty\subset \cI$ is a
strictly increasing sequence,
     \begin{equation}\label{IV_01}
0=x_0<x_1<x_2<\dots<x_n<\dots<b\leq+\infty, \qquad\text{and}\qquad \lim_{n\to\infty}x_n=b.
    \end{equation}
We denote $\gd_n:=x_{n}-x_{n-1}$.
Consider the following symmetric operator in $L^2(\mathcal{I})$
      \begin{equation}\label{IV.1.1_01}
{\rm H}_{\min}=-\frac{\rD^2}{\rD x^2},\qquad \dom(\rH_{\min})=W^{2,2}_0(\cI\setminus X).
       \end{equation}
Clearly, $\rH_{\min}$ is closed and
       \begin{equation}\label{IV.1.1_02}
\rH_{\min}=\oplus_{n=1}^{\infty} \rH_n,\quad \text{where}\quad \rH_n=-\frac{\rD^2}{\rD x^2},\quad \dom(\rH_{n})=W^{2,2}_0[x_{n-1},x_{n}].
        \end{equation}

\textbf{1.} Note that $\rH_{\min}\ge 0.$  It is known (see for
instance \cite{Gor84}) that \emph{Friedrichs' extension} $\rH^{\rm
F}_n$ of $\rH_n$ is defined by the Dirichlet boundary conditions,
i.e., $\dom(\rH_n^{\rm F})=\{f\in W^{2,2}[x_{n-1},x_{n}]:\
f(x_{n-1}+)=f(x_{n}-)=0\}$. Therefore, the Friedrichs' extension
$\rH^{\rm F}$ of $\rH_{\min}$ is $\rH^{\rm F} =
\oplus_{n=1}^{\infty} \rH_n^{\rm F}$, that is
\begin{equation}\label{IV.1.1_04}
\rH_{\rm F}=-\frac{\rD^2}{\rD x^2},\qquad \dom(\rH_\rF)=\{f\in W_2^2(\mathcal{I}\setminus X): \ f(0)=f(x_n+)=f(x_n-)=0,\ n\in\N\}.
\end{equation}
It is easily seen that a triplet $\wt\Pi_n=\{\C^2,
\wt\Gamma_0^{(n)},\wt\Gamma_1^{(n)}\}$ given by
      \begin{equation}\label{IV.1.1_05}
\wt\Gamma_0^{(n)}f:=\left(\begin{array}{c}
                 f(x_{n-1}+)\\
                 - f(x_{n}-)
                       \end{array}\right),\qquad  \wt\Gamma_1^{(n)}f:=\left(\begin{array}{c}
                                                                           f'(x_{n-1}+)\\
                                                                            f'(x_{n}-)
                                                                          \end{array}\right), \qquad f\in
                                                                          W_2^2[x_{n-1},x_n],
\end{equation}
forms a boundary triplet for $\rH_n^*$ satisfying
$\ker(\wt\Gamma_0^{(n)})= \dom(\rH_n^F).$
Moreover, $\rH_n=\gd_n^{-2} U_n^{-1} S_1 U_n$, where 
$S_1:=-\frac{\rD^2}{\rD x^2},\ \dom(S_1)=W_0^{2,2}[0,1]$, and
$(U_n f)(x):= \sqrt{\gd_n}f(\gd_{n}x+x_{n-1}).$ Clearly, $U_n$
isometrically maps $L^2[x_{n-1}, x_{n}]$ onto $L^2[0,1].$ As it
follows from Lemma \ref{lem_koch}, a triplet $\wt{\Pi}=\oplus_{n\in\N}\wt{\Pi}_n$ forms a boundary triplet for the operator
$\rH_{\min}^*:= (\rH_{\min})^*=\rH_{\max}$ whenever
    \begin{equation}\label{IV.1.1_03}
0\ < \ \gd_*=\inf_{n\in\N} \gd_n\ \ \leq\ \ \gd^*=\sup_{n\in
\N}\gd_n\ <\ +\infty.
       \end{equation}
If $d_*=0,$ then the direct sum $\wt\Pi=\oplus_{n=1}^{\infty}
\wt\Pi_n$  of triplets \eqref{IV.1.1_05} is not  a boundary
triplet for $\rH_{\max}$. We regularize the triplet $\wt\Pi$ by
applying Corollary \ref{cor_III.2.2_02} in order to obtain a
direct sum triplet $\Pi=\oplus_{n=1}^{\infty}\Pi_n$  for the
operator $\rH_{\min}^*$, assuming only that
      \begin{equation}\label{IV.1.1_07}
d^*=\sup_{n\in \N} d_n<+\infty,
    \end{equation}
\begin{theorem}\label{th_bt_2}
Assume condition  \eqref{IV.1.1_07} and define the mappings
$\Gamma_j^{(n)}: W_2^2[x_{n-1},x_n]\to\C^2$, $n\in \N$,  $j=0,1$,
 by setting
\begin{equation}\label{IV.1.1_06}
\Gamma_0^{(n)}f:=\left(\begin{array}{c}
                 \gd_n^{1/2}  f(x_{n-1}+)\\
                 -\gd_n^{1/2}  f(x_{n}-)
                       \end{array}\right),\qquad \Gamma_1^{(n)}f:=\left(\begin{array}{c}
                                                                           \frac{\gd_n f'(x_{n-1}+)+(f(x_{n-1}+)-f(x_{n}-))}{\gd_n^{3/2}}\\
                                                                           \frac{\gd_n f'(x_{n}-)+(f(x_{n-1}+)-f(x_{n}-))}{\gd_n^{3/2}}
                                                                          \end{array}\right).
\end{equation}
Then:

\item $(i)$\ For any $n\in \N$ the triplet $\Pi_n=\{\C^2,
\Gamma_0^{(n)},\Gamma_1^{(n)}\}$ is a boundary triplet for
$\rH_{n}^*$.

\item $(ii)$\ The direct sum $\Pi=\oplus_{n=1}^\infty\Pi_n$
is a boundary triplet for the operator $\rH_{\min}^*$.
      \end{theorem}
      \begin{proof}
(i)\ Straightforward.

(ii)\  The  Weyl function $ \widetilde{M}_{n}(\cdot)$
corresponding to the triplet $\wt\Pi_n$ of the form
\eqref{IV.1.1_05} is
           \begin{equation}\label{IV.1.1_09}
 \widetilde{M}_{n}(z)=\left(\begin{array}{cc}
                           -\frac{\sqrt{z}\cos(\sqrt{z}\gd_n)}{\sin(\sqrt{z}\gd_n)} & -\frac{\sqrt{z}}{\sin(\sqrt{z}\gd_n)}\\
                           -\frac{\sqrt{z}}{\sin(\sqrt{z}\gd_n)} & -\frac{\sqrt{z}\cos(\sqrt{z}\gd_n)}{\sin(\sqrt{z}\gd_n)}
                                                                          \end{array}\right),\qquad
                                                                          z\in\C_+.
         \end{equation}
Comparing definitions \eqref{IV.1.1_05} and
\eqref{IV.1.1_06} of triplets $\Pi_n$ and ${\widetilde\Pi}_n$,
respectively,  we get
   \begin{equation}\label{IV.1.1_09A}
\Gamma^{(n)}_0= R_n{\widetilde\Gamma}_0^{(n)},\quad
 \Gamma^{(n)}_1=
R^{-1}_n({\widetilde\Gamma}_1^{(n)}-Q_n{\widetilde\Gamma}_0^{(n)}),\quad
\text{and}\quad M_n(z) = R^{-1}_n(\widetilde{M}_n(z)-Q_n)R^{-1}_n,
       \end{equation}
where
\begin{equation}\label{IV.1.1_09B}
R_n = R_n^* :=
\begin{pmatrix}
\gd^{1/2}_n & 0\\
0         & \gd^{1/2}_n
\end{pmatrix} \qquad
\text{and}\qquad Q_n=\frac{1}{\gd_n}\left(\begin{array}{cc}
                           -1 & -1\\
                           -1 & -1
                           \end{array}\right) = \wt M_n(0).
\end{equation}
 It follows from
\eqref{IV.1.1_09A}, \eqref{IV.1.1_09B}, and \eqref{IV.1.1_09} that
\begin{equation}\label{IV.1.1_08}
M_n(0)=0,\quad
M_n'(0)=R^{-1}_n{\widetilde M}_n'(0)R^{-1}_n=R^{-1}_n
\begin{pmatrix}
\gd_n/3 & -\gd_n/6\\
-\gd_n/6 & \gd_n/3
\end{pmatrix}
 R^{-1}_n=
\begin{pmatrix}
1/3 & -1/6\\
-1/6 & 1/3
\end{pmatrix}.
         \end{equation}
Relations  \eqref{IV.1.1_08} yield conditions \eqref{III.2.2_11}      
One completes the proof by applying Corollary
\ref{cor_III.2.2_02}.
\end{proof}
        \begin{remark}\label{rem_Mih}
Let $\gd_*=0$.
Hence  both families
$\{\widetilde{M}_n(\I)\}_{n\in\N}$ and
$\{\widetilde{M}_n(\I)\}^{-1}_{n\in\N}$ (see \eqref{IV.1.1_09}) are unbounded. By
Proposition \ref{prop_III.1_01}, neither $\wt
\Pi=\oplus_{n=1}^\infty \wt \Pi_n$ no  $\wt
{\Pi^{(1)}} =\oplus_{n=1}^\infty \wt {\Pi^{(1)}}_n$,  where  $\wt\Pi_n = \{\C^2,
\widetilde{\Gamma}_0^{(n)},\widetilde{\Gamma}_1^{(n)}\}$ is
defined by \eqref{IV.1.1_05}  and  $\wt {\Pi^{(1)}}_n:=
\{\C^2,-\widetilde{\Gamma}_1^{(n)}, \widetilde{\Gamma}_0^{(n)}\}$,   
forms a generalized boundary triplet for $\rH_{\min}^*$.
Moreover, by Proposition \ref{prop_III.1_02} $(i)$, the mappings
$\widetilde{\Gamma}_0=\oplus_{n=1}^\infty\widetilde{\Gamma}_0^{(n)}$ and
$\widetilde{\Gamma}_1=\oplus_{n=1}^\infty\widetilde{\Gamma}_1^{(n)}$ are
unbounded. Note that, the latter might be checked by restricting the
mappings $\widetilde{\Gamma}_0$ and $\widetilde{\Gamma}_1$ on $\cN_\I(\rH_{\min})$.

Note also that $\widetilde{\Gamma}_0$ coincides with the mapping $\Gamma^2$
in \cite[Theorem 1]{Mih_93}. Hence the triplet $\Pi$ constructed in \cite[Theorem 1]{Mih_93} is not an ordinary boundary triplet.
\end{remark}
    \begin{remark}
Let us sketch another proof of Theorem  \ref{th_bt_2}.
Simple calculations with account of  \eqref{IV.1.1_07} yield that
the family $\{M_n(\I)\}_{n=1}^\infty$ is bounded. Moreover, it
follows from \eqref{IV.1.1_09} that
\[
 \lim_{n_k\to\infty} M_{n_k}(\I) = \I\lim_{n_k\to\infty} \im
 M_{n_k}(\I)=
\I  \left(\begin{array}{cc}
                           1/3 & -1/6\\
                           -1/6& 1/3
                           \end{array}\right)\quad \text{whenever}\quad \lim_{n_k\to\infty}\gd_{n_k} =0.
\]
Hence, by Theorem \ref{th_criterion(bt)},
$\Pi=\oplus_{n=1}^\infty\Pi_n$ defined by \eqref{IV.1.1_06} forms
a boundary triplet for $\rH_{\min}^*$.
   \end{remark}

\begin{proposition}\label{prop_IV.1.1_01}
Let $\Pi$ be the boundary triplet defined in Theorem \ref{th_bt_2}
and $M(\cdot)$
the corresponding Weyl function. 
If condition \eqref{IV.1.1_07} is satisfied, then
\begin{equation}\label{IV.1.1_10}
M(-a^2)\rightrightarrows -\infty \quad\text{as} \quad a\to+\infty.
\end{equation}
\end{proposition}
\begin{proof}
By Theorem \ref{th_III.1_01} $(ii)$, the Weyl function $M(\cdot)$
has the form $M(z)=\oplus_{n=1}^\infty M_n(z)$, where $M_n(\cdot)$
is defined by \eqref{IV.1.1_09A}, \eqref{IV.1.1_09} and
\eqref{IV.1.1_09B}. Consider the following matrix-function
\begin{equation}\label{IV.1.1_11}
M(-a^2,x) :=\left(\begin{array}{cc}
                           F_a(x) & G_a(x)\\
                           G_a(x) & F_a(x)
                           \end{array}\right),\qquad  x>0,
\end{equation}
where
\[
F_a(x) :=\frac{1}{x^2}-a\frac{\cosh ax}{x \sinh ax} ,\qquad G_a(x)
:= \frac{1}{x^2}-\frac{a}{x \sinh ax}.
\]
It is easy to check that
\[
F_a(x)<0\quad \text{and} \quad G_a(x)>0\quad \text{for} \quad x>0.
\]
Since $\sigma(M(-a^2,x))=\{F_a(x)+G_a(x),F_a(x)-G_a(x)\}$, we get
\[
M(-a^2,x)\leq (F_a(x)+G_a(x))I_2,\qquad x>0.
\]
Further, consider the function
\[
f(x)=\frac{2}{x^2}-\frac{1+\cosh x}{x \sinh x}\ (\ =F_1(x)+G_1(x)).
\]
Note that $f(x)<0$ if $x>0$. Moreover, $f$ is continuous on $\R_+$ and
\[
\lim_{x\to+0}f(x)=-\frac{1}{6},\qquad \lim_{x\to+\infty}f(x) = 0.
\]
Note also that $\lim_{x\to+\infty}x^2f'(x)=1$. Hence $f'(x)>0$ for
$x\geq x_0$ with sufficiently large $x_0\in\R_+$. Since
$F_a(x)+G_a(x)=a^2f(ax)$, for  $a\geq a_0>0$ large enough we
obtain
\[
\sup_{x\in (0,\gd^*)}(F_a(x)+G_a(x))= 
\frac{2}{(\gd^*)^2}-\frac{a}{\gd^*}\cdot\frac{1+\cosh a\gd^*}{ \sinh a\gd^*}\leq -2\frac{a}{\gd^*}+\frac{2}{(\gd^*)^2}.
\]
Note that $M_n(-a^2)=M(-a^2, \gd_n)$. Combining this fact with the last inequality, we  obtain
\begin{equation}\label{IV_est_lsb}
M(-a^2)=\oplus_{n=1}^\infty M_n(-a^2)\leq -\frac{a}{\gd^*}I_{l_2},
\qquad a\geq \max\{a_0,2/\gd^*\}.
\end{equation}
This completes the proof.
\end{proof}

Combining Theorem \ref{th_bt_2} with Proposition
\ref{prop_II.1.2_01}, we arrive at the following parametrization
of the set $\Ext {\rH_{\min}}$ of closed proper extensions of the
operator $\rH_{\min}:$
      \begin{equation}\label{IV.1.1_12'}
\widetilde{\rH}=\rH_{\Theta}:=\rH_{\min}^* \lceil
\dom(\rH_\Theta),\qquad \dom(\rH_\Theta)=\{f\in
\dom(\rH_{\min}^*): \ \{\Gamma_0f,\Gamma_1f\}\in\Theta\},
     \end{equation}
where $\Theta\in \widetilde{\mathcal{C}}(l_2)$ and $\Gamma_0,\
\Gamma_1$ are defined by \eqref{IV.1.1_06}.
      \begin{theorem}\label{th_IV.1.1_02}
Let $\Pi=\oplus_{n=1}^\infty \Pi_n$ be a boundary triplet for
$\rH_{\min}^*$ defined in Theorem \ref{th_bt_2}, $\Theta,
\widetilde{\Theta}\in \widetilde{\mathcal{C}}(\cH)$, and
$\rH_\Theta, \rH_{\widetilde{\Theta}}\in \Ext {\rH_{\min}}$
proper extensions of $\rH_{\min}$ defined by \eqref{IV.1.1_12'}.
Then:
\item $(i)$ \ The operator $\rH_\Theta$ is symmetric
(self-adjoint) if and only if so is $\Theta$, and
$\mathrm{n}_\pm(\rH_{\min})=\mathrm{n}_\pm(\Theta)$.
\item $(ii)$ \ The self-adjoint (symmetric) operator $\rH_\Theta$
is lower semibounded if and only if so is
 $\Theta$.
\item $(iii)$ \ For any $p\in(0,\infty]$,
$z\in\rho(\rH_{\Theta})\cap\rho(\rH_{\widetilde{\Theta}})$,
and $\zeta\in\rho(\Theta)\cap\rho(\widetilde{\Theta})$ the following equivalence holds 
\[
(\rH_{\Theta}-z)^{-1}-(\rH_{\widetilde{\Theta}}-z)^{-1}\in \mathfrak{S}_p\quad\Longleftrightarrow \quad (\Theta-\zeta)^{-1}-(\widetilde{\Theta}-\zeta)^{-1}\in \mathfrak{S}_p.
\]
\item $(iv)$ \ The  operator $\rH_\Theta = \rH_\Theta^*$ has
discrete spectrum if and only if $\gd_n\searrow 0$ and $\Theta$
has discrete spectrum.
      \end{theorem}
       \begin{proof}
$(i)$ is immediate from Proposition \ref{prop_II.1.2_01}.

$(ii)$ Combining Propositions \ref{prop_II.1.5_02} with
Proposition \ref{prop_IV.1.1_01} yields the first statement. Then
the second one is implied by estimate \eqref{IV_est_lsb}.

$(iii)$ is implied by  Proposition
\ref{prop_II.1.4_02}.

$(iv)$  First we show that conditions are sufficient. Indeed, the operator
\begin{equation}\label{h_0}
\rH_0:=\rH_{\min}^*\lceil\ker(\Gamma_0)=\oplus_{n\in\N}\rH_{n0},\qquad \rH_{n0}:=\rH_n^*\lceil\ker (\gG_0^{(n)}),
\end{equation}
 has discrete spectrum if
$\lim_{n\to\infty}\gd_n=0$.
Moreover, the Krein resolvent formula and discreteness of
$\sigma(\Theta)$ implies
$\mathcal{R}_{H_{\Theta}}(z)-\mathcal{R}_{H_0}(z)\in
\mathfrak{S}_\infty,\ z\in\C_+,$ and hence
$\mathcal{R}_{H_{\Theta}}(z)\in \mathfrak{S}_\infty$.

Let us show  that condition $\gd_n\searrow 0$ is necessary for
discreteness of  $\sigma(H_{\Theta}).$   Without loss of
generality assume that $0\in\rho(\rH_\Theta)$. Assume also that
$\limsup_{n\to\infty}\gd_n>0$ and $\rH_{\Theta}$ has discrete
spectrum. Then there exists a sequence
$\{\gd_{n_k}\}_{k=1}^\infty$ such that $\gd_{n_k}\geq d_*/2 > 0.$
For $\varepsilon \in (0, d_*/2)$, define the function
\[
\varphi_\varepsilon(\cdot) \in W_2^2(\R),\qquad
\varphi_\varepsilon(x)= \left\{\begin{array}{c}
1,\quad  \varepsilon\leq x\leq d_*-\varepsilon,\\
0, \quad \quad \quad x\notin [0,d_*].
\end{array}\right.
\]
Note that $\varphi_k(x):=P_{\cI}\varphi_\varepsilon(x+x_{n_k})\in
\dom(\rH_{\Theta})$, where $P_{\cI}$ is the orthoprojection in
$L^2(\R)$ onto $L^2(\cI)$. Moreover, $\|\varphi_k\|_{L^2}\equiv
const$ and $\|\rH_{\Theta}\varphi_k\|_{L^2}\equiv const$. Since
the functions $\varphi_k(\cdot)$ have disjoint supports, the
operator $(\rH_{\Theta})^{-1}$ is not compact. Contradiction.
        \end{proof}
\begin{corollary}
$\rH_{\Theta}$ is nonnegative if and only if the linear relation $\Theta$ is nonnegative. Moreover,
 if $a$ is large enough, then $\rH_\Theta\ge -a^2$ whenever  $\Theta \ge -\frac{a}{\gd^*}I_{l_2}$.
\end{corollary}
\begin{proof}
Since $M(0)=0$, by \cite[Theorem 4]{DM91}, we get the first part. Moreover, we have the estimate $M(-a^2)\leq -a/d^* I$ (see the proof of Proposition \ref{prop_IV.1.1_01}), and Krein's formula \eqref{II.1.4_01} completes the proof.
\end{proof}

\textbf{2.} Alongside boundary triplet \eqref{IV.1.1_06} consider
another boundary triplet. Namely, define $\wt\Pi_n=\{\cH,\
\wt\Gamma_0^{(n)},\ \wt\Gamma_1^{(n)}\}$ for the operator
$\rH_n^*, \ n\in \N,$ by setting
    \begin{equation}\label{IV.1.1_12B}
\cH=\C^2,\qquad\wt\Gamma_0^{(n)}f:=\left(\begin{array}{c}
                   f(x_{n-1}+)\\
                   f'(x_{n}-)
                       \end{array}\right),\qquad \wt\Gamma_1^{(n)}f:=\left(\begin{array}{c}
                                                                           f'(x_{n-1}+)\\
                                                                           f(x_{n}-)
                                                                          \end{array}\right), \quad f\in W_2^2[x_{n-1},x_n].
\end{equation}
In the following theorem we regularize the family
$\{\wt\Pi_n\}_{n=1}^\infty$ 
in such a way  that the direct sum
of new boundary triplets $\Pi_n$  is
already a boundary triplet for $\rH_{\min}^*=\oplus_{n=1}^\infty \rH_n^*$ if $d^*<\infty$.

         \begin{theorem}\label{th_bt_1}
Assume condition  \eqref{IV.1.1_07} and define the mappings
$\Gamma_j^{(n)}: W_2^2[x_{n-1},x_n]\to\C^2$, $n\in \N$,  $j=0,1$,
 by setting
    \begin{equation}\label{IV.1.1_12}
\Gamma_0^{(n)}f:=\left(\begin{array}{c}
                 \gd_n^{1/2}  f(x_{n-1}+)\\
                 \gd_n^{3/2}  f'(x_{n}-)
                       \end{array}\right),\quad \Gamma_1^{(n)}f:=\left(\begin{array}{c}
                                                                           \frac{f'(x_{n-1}+)-f'(x_{n}-)}{\gd_n^{1/2}}\\
                                                                           \frac{f(x_{n}-)-f(x_{n-1}+)-\gd_nf'(x_{n}-)}{\gd_n^{3/2}}
                                                                          \end{array}\right).
\end{equation}
Then:
\item
$(i)$\  For any $n\in \N$ the triplet  $\Pi_n=\{\C^2,\Gamma_0^{(n)},\Gamma_1^{(n)}\}$ is a boundary
triplet for $\rH_n^*$.
\item
$(ii)$\ The direct sum  $\Pi=\oplus_{n=1}^\infty \Pi_n$
is a boundary triplet for the operator $\rH_{\max}=\rH_{\min}^*$.
      \end{theorem}
\begin{proof}
$(i)$ Straightforward.

$(ii)$\ The  Weyl function of $\rH_n^*$ corresponding to the
triplet $\wt\Pi_n$ defined by \eqref{IV.1.1_12A} is
\begin{equation}\label{IV.1.1_14}
\widetilde{M}_{n}(z)=\left(\begin{array}{cc}
                           \frac{\sqrt{z}\sin(\sqrt{z}\gd_n)}{\cos(\sqrt{z}\gd_n)} & \frac{1}{\cos(\sqrt{z}\gd_n)}\\
                           \frac{1}{\cos(\sqrt{z}\gd_n)} & \frac{\sin(\sqrt{z}\gd_n)}{\sqrt{z}\cos(\sqrt{z}\gd_n)}
                                                                          \end{array}\right).
\end{equation}
Comparing definitions \eqref{IV.1.1_12A} and \eqref{IV.1.1_12}, we
get that the triplets $\Pi_n$ and ${\widetilde\Pi}_n$ are
connected by \eqref{IV.1.1_09A}, where the matrices $R_n$ and
$Q_n$ are defined by
       \begin{equation}\label{IV.1.1_19}
 R_n :=\left(\begin{array}{cc}
                           \gd_n^{1/2} & 0\\
                           0 & \gd_n^{3/2}
                           \end{array}\right)\quad\text{and}\quad
Q_n :=\wt M_n(0)=\left(\begin{array}{cc}
                           0 & 1\\
                           1 & \gd_n
                           \end{array}\right).
\end{equation}
Hence $M_n(z)=R_n^{-1}(\widetilde{M}_n(z)-Q_n)R_n^{-1}$ is the Weyl function
corresponding to the triplet $\Pi_n$.
It follows from  \eqref{IV.1.1_14} and \eqref{IV.1.1_19}  that
   \begin{equation}\label{IV.1.1_08B}
 M_{n}(0)=0,
\qquad M_{n}'(0)=R^{-1}_n{\widetilde M}'_n(0)R^{-1}_n =R^{-1}_n\begin{pmatrix}
\gd_n & \gd_n^2/2\\
\gd_n^2/2 &  \gd_n^3/3
\end{pmatrix}R^{-1}_n =
\begin{pmatrix}
1 & 1/2\\
1/2 &1/3
\end{pmatrix}.
         \end{equation}
One completes the proof by applying Theorem \ref{th_III.2.2_01}.
            \end{proof}
      \begin{remark}\label{rem_IV.1.1_01}
Clearly, all statements  of Theorem \ref{th_IV.1.1_02} with
exception of $(ii)$  remain valid for the boundary triplet
$\Pi=\oplus_1^{\infty}\Pi_n$ with $\Pi_n$ defined by
\eqref{IV.1.1_12} in place of  \eqref{IV.1.1_06}.
         \end{remark}
\begin{corollary}\label{cor4.8}
Let $\widetilde{\Pi}_n$ be a boundary triplet for $\rH_n^*$
defined by \eqref{IV.1.1_12A} and  $\wt {\Pi^{(1)}}_n:=
\{\C^2,-\widetilde{\Gamma}_1^{(n)}, \widetilde{\Gamma}_0^{(n)}\}$.
Let also $\widetilde{\Pi}=\oplus_{n=1}^\infty\widetilde{\Pi}_n$
and $\wt {\Pi^{(1)}} :=\oplus_{n=1}^\infty \wt {\Pi^{(1)}}_n$ be
direct sums of boundary triplets and $\gd_*=0$. Then:

\item $(i)$ $\widetilde{\Pi}$ and $\widetilde{\Pi}^{(1)}$ are generalized boundary triplets for $\rH_{\min}^*$.

\item $(ii)$ $\widetilde{\Pi}$ and $\widetilde{\Pi}^{(1)}$ are not
ordinary boundary triplets for $\rH_{\min}^*$.

\item $(iii)$ The operators $(\rH_{\min})_{*0}$ and
$(\rH_{\min})_{*1}$ (see \eqref{III.1_03}) are self-adjoint and $(\rH_{\min})_{*j} = \oplus_{n=1}^\infty \rH_{nj}$.  

\item $(iv)$ The mappings $\wt{\Gamma}_0$ and $\wt{\Gamma}_1$ are
closed  and unbounded on $\mathfrak{H}_+=\dom(\rH_{\min}^*)$.

\item $(v)$ $(\rH_{\min})_{*0}$ and $(\rH_{\min})_{*1}$
are not transversal.
\end{corollary}
\begin{proof}
$(i)$ It follows from \eqref{IV.1.1_14} that the  families
$\{\widetilde{M}_n(\I)\}_{n=1}^{\infty}$ and
$\{\widetilde{M}^{-1}_n(\I)\}_{n=1}^{\infty}$  are bounded if
$\gd^*<\infty$. It remains to apply Proposition
\ref{prop_III.1_01}.

$(ii)$ If $\lim_{k\to\infty}\gd_{n_k} = 0,$ then $
\lim_{k\to\infty}\im \widetilde{M}_{n_k}(\I)=\im
\left(\begin{array}{cc}
                           0 & 1\\
                           1 & 0
                           \end{array}\right)=\left(\begin{array}{cc}
                           0 & 0\\
                           0 & 0
                           \end{array}\right)$.
Thus, the second of conditions \eqref{WF_criterion} is violated,
hence neither $\wt{\Pi}$ no $\wt{\Pi^{(1)}}$ forms a
boundary triplet for $\rH_{\min}^*.$ 

$(iii)$ follows from $(i)$ and Theorem \ref{th_III.1_01} $(vi)$.

$(iv)$ Clearly, $\wt{\Gamma}_0$ and $\wt{\Gamma}_1$ are unitarily
equivalent. Hence $\wt{\Gamma_0}$ and $\wt{\Gamma}_1$ might be
bounded only simultaneously. Combining $(ii)$ with Proposition
\ref{prop3.8}, we conclude that both $\wt{\Gamma}_0$ and
$\wt{\Gamma}_1$ are unbounded. Further, by
Theorem \ref{th_III.1_01} $(iv)$, $\wt{\Gamma}_j$
is closable. Since, by $(iii)$,  $\ker(\wt{\Gamma}_j) =
\oplus_{n=1}^\infty \dom(\rH_{nj})$ is closed in $\mathfrak{H}_+$
and $\ran (\wt{\Gamma}_j)=\cH$ is closed, the mapping
$\wt{\Gamma}_j$ is closed.

$(v)$ follows from $(iii)$ and Proposition
\ref{prop_III.1_02}(ii).
      \end{proof}
            \begin{remark}\label{rem_Mih2}
 Corollary \ref{cor4.8} shows that condition
$C_1<\infty$ in Proposition \ref{prop_III.1_02} is only sufficient
for $\Pi=\oplus_{n\in\N}\Pi_n$ to form a generalized boundary triplet.
       \end{remark}
%
%
\section{Schr\"odinger operators with $\delta$-interactions}\label{Subsec_IV.2_delta}
%
%
Let $\cI=[0,b)$ and let $X=\{x_n\}_{n=1}^\infty$ be defined by \eqref{IV_01}. In what follows we will always assume that condition \eqref{IV.1.1_07} is satisfied, i.e. $\gd^*=\sup_{n}\gd_n<\infty$.\\
The main object of this section is the formal differential expression
\begin{equation}\label{IV.2.0_01}
\ell_{X,\alpha}:=-\frac{\rD^2}{\rD x^2}+\sum_{n=1}^\infty \alpha_n\delta(x-x_n),\qquad \gA_n\in\R.
\end{equation}
In  $L^2(\cI)$, one associates with \eqref{IV.2.0_01} a symmetric
differential operator
\begin{equation}\label{IV.2.0_02}
\rH^0_{X,\alpha}:=-\frac{\rD^2}{\rD x^2},\qquad \dom(\rH^0_{X,\gA})=\{f\in W^{2,2}_{\comp}(\cI\setminus X): \ \begin{array}{c}
f'(0)=0,\ f(x_n+)=f(x_n-)\\ f'(x_n+)-f'(x_n-)=\alpha_n f(x_n)
\end{array} \}.
    \end{equation}
Denote by $\rH_{X,\gA}$ the closure of $\rH^0_{X,\gA}$,
$\rH_{X,\gA}=\overline{\rH^0_{X,\gA}}$.

\subsection{Parametrization of the operator $\rH_{X,\gA}$}\label{sss_IV.2.1_boun_op}

Let $\Pi^1=\{\cH,\Gamma_0^1,\Gamma_1^1\}$ and
$\Pi^2=\{\cH,\Gamma_0^2,\Gamma_1^2\}$ be the boundary triplets
defined in Theorems \ref{th_bt_2} and \ref{th_bt_1}, respectively.
By Proposition   \ref{prop_II.1.2_01}, the extension
$\rH_{X,\gA} (\in\Ext {\rH_{\min}})$ admits two representations
      \begin{equation}\label{IV.1.1_12A}
\rH_{X,\gA}  =\rH_{\Theta_j}:=\rH_{\min}^* \lceil
\dom(\rH_{\Theta_j}),\quad \dom(\rH_{\Theta_j})=\{f\in
\dom(\rH_{\min}^*): \
\{\Gamma^j_0f,\Gamma^j_1f\}\in\Theta_j\},\quad j=1, 2.
     \end{equation}
(cf. \eqref{IV.1.1_12'}) with  closed symmetric linear relation
$\Theta_j\in \widetilde{\mathcal{C}}(\cH),\ j=1,2$. We show
that $\Theta_2$ as well as the operator part $\Theta'_1$ of
$\Theta_1$ is a Jacobi matrix.

\textbf{1. The first parametrization.} We begin with the triplet
$\Pi^2=\{\cH,\Gamma_0^2,\Gamma_1^2\}$ constructed in Theorem
\ref{th_bt_1}. For any $\alpha$ the operators $\rH_{X,\gA}$ and
$\rH^{(2)}_0:=\rH^*_{\min}\lceil \ker(\Gamma^2_0)$ are disjoint.
Hence $\Theta_2$ in \eqref{IV.1.1_12A} is a (closed) operator in
$\cH=l_2(\N)$. 
More precisely,
consider the Jacobi matrix
\begin{equation}\label{IV.2.1_01}
B_{X,\gA}=\left(%
\begin{array}{cccccc}
  0& -\gd_1^{-2} & 0 & 0& 0  &  \dots\\
  -\gd_1^{-2} & -\gd_1^{-2}& \gd_1^{-3/2}\gd_2^{-1/2}& 0 & 0&  \dots\\
  0 & \gd_1^{-3/2}\gd_2^{-1/2} & \alpha_1\gd_2^{-1} & -\gd_2^{-2} & 0&   \dots\\
  0 & 0 & -\gd_2^{-2} & -\gd_2^{-2} & \gd_2^{-3/2}\gd_3^{-1/2}&  \dots\\
  0 & 0 & 0 & \gd_2^{-3/2}\gd_3^{-1/2} & \alpha_2\gd_3^{-1}&   \dots\\
\dots& \dots&\dots&\dots&\dots&\dots\\
 \end{array}%
\right).
\end{equation}
Let $\tau_{X,\gA}$ be a second order difference expression associated with \eqref{IV.2.1_01}.
One defines the corresponding minimal symmetric operator in $l_2$ by (see \cite{Akh, Ber68})
\begin{gather}\label{IV.2.1_02}
B^0_{X,\gA}f:=\tau_{X,\gA}f,\qquad f\in\dom(B^0_{X,\gA}):=l_{2,0},
\quad\text{and}\quad B_{X,\gA}=\overline{B^0_{X,\gA}}.
\end{gather}
 Recall that $B_{X,\gA}$\footnote{Usually  we will identify the Jacobi matrix with (closed) minimal symmetric operator
 associated with it. Namely, we denote by $B_{X,\gA}$ the Jacobi matrix \eqref{IV.2.1_01}
 as well as the minimal closed symmetric operator \eqref{IV.2.1_02}.}
 has equal deficiency indices and $\mathrm{n}_+(B_{X,\gA})=\mathrm{n}_-(B_{X,\gA})\leq 1$.

Note that $B_{X,\gA}$ admits a representation
      \begin{equation}\label{IV.2.2_01A}
B_{X,\gA}=R^{-1}_X(\widetilde{B}_{\gA}-Q_X)R^{-1}_X,\quad \text{where}\quad \widetilde{B}_{\gA} := \left(%
\begin{array}{cccccc}
  0& 0 & 0 & 0& 0  & \dots\\
  0 & 0& 1& 0 & 0&  \dots\\
  0 & 1 & \alpha_1 &   0& 0& \dots\\
  0 & 0 & 0 & 0 & 1&  \dots\\
  0 & 0 & 0 & 1 & \alpha_2&  \dots\\
\dots& \dots&\dots&\dots&\dots&\dots\\
 \end{array}%
\right)
\end{equation}
and $R_X=\oplus_{n=1}^\infty R_n,\ Q_X=\oplus_{n=1}^\infty Q_n$ are defined by \eqref{IV.1.1_19}.
\begin{proposition}\label{prop_IV.2.1_01}
Let $\Pi^2=\{\cH,\Gamma_0^2,\Gamma_1^2\}$ be the boundary triplet
for $\rH_{\min}^*$ constructed in Theorem \ref{th_bt_1} and let
$B_{X,\gA}$ be the minimal Jacobi operator defined by
\eqref{IV.2.1_01}--\eqref{IV.2.1_02}. Then $\Theta_2=B_{X,\gA}$,
i.e.,
   \begin{equation*}
\rH_{X,\gA}=\rH_{B_{X,\gA}}=\rH_{\min}^*\lceil\dom(\rH_{B_{X,\gA}}),\qquad \dom(\rH_{B_{X,\gA}})=\{f\in W^{2,2}(\cI\setminus X):\Gamma_1^2f=B_{X,\gA}\Gamma_0^2f\}.
     \end{equation*}
     \end{proposition}
  \begin{proof}
 Let $f\in W_{\comp}^{2,2}(\cI\setminus X)$. Then $f\in \dom(\rH_{X,\gA})$ if and only
 if $\wt{\Gamma}_1^2f=\wt{B}_\alpha\wt{\Gamma}_0^2f.$
 Here
$\wt{\Gamma}_j^2 :=\oplus_{n\in\N}\wt{\Gamma}_j^{(n)}$ where
$\wt{\Gamma}_j^{(n)},\ j=0,1,$ are defined  by
\eqref{IV.1.1_12A}, and $\wt{B}_\alpha$ is defined  by
\eqref{IV.2.2_01A}.
 Combining  \eqref{IV.1.1_09A}, \eqref{IV.1.1_19}  with  \eqref{IV.2.2_01A},
we rewrite the equality
$\wt{\Gamma}_1^2f=\wt{B}_\alpha\wt{\Gamma}_0^2f$  as
$\Gamma_1^2f=B_{X,\gA}\Gamma_0^2f$.

 Taking the closures one completes  the
proof.
      \end{proof}
\begin{remark}
Note that the matrix \eqref{IV.2.1_01} has negative off-diagonal
entries, although, in the classical theory of Jacobi operators,
off-diagonal entries are assumed to be positive. But it is known
(see, for instance, \cite{Tes_98}) that the (minimal) operator
$B_{X,\alpha}$
is unitarily equivalent to the minimal Jacobi operator associated
with the matrix
   \begin{equation}\label{IV.2.1_01'}
B'_{X,\alpha} :=\left(%
\begin{array}{cccccc}
  0& \gd_1^{-2} & 0 & 0& 0  &  \dots\\
  \gd_1^{-2} & -\gd_1^{-2}& \gd_1^{-3/2}\gd_2^{-1/2}& 0 & 0&   \dots\\
  0 & \gd_1^{-3/2}\gd_2^{-1/2} & \alpha_1\gd_2^{-1} & \gd_2^{-2} & 0&   \dots\\
  0 & 0 & \gd_2^{-2} & -\gd_2^{-2} & \gd_2^{-3/2}\gd_3^{-1/2}&  \dots\\
  0 & 0 & 0 & \gd_2^{-3/2}\gd_3^{-1/2} & \alpha_2\gd_3^{-1}&   \dots\\
\dots& \dots&\dots&\dots&\dots&\dots\\
 \end{array}%
\right).
      \end{equation}
In the sequel we will identify the operators $B_{X,\alpha}$ and
$B'_{X,\alpha}$ when investigating those spectral properties of the operator $\rH_{X,\gA}$, which are invariant under
unitary transformations.
        \end{remark}
\textbf{2. The second parametrization.} 
Let us consider the boundary triplet
$\Pi^1=\{\cH,\Gamma_0^1,\Gamma_1^1\}$ constructed in Theorem
\ref{th_bt_2}. Now the operators $\rH_{X,\gA}$ and
$\rH_0^{(1)}:=\rH^*_{\min}\lceil \ker(\Gamma^1_0)$     
are not disjoint, hence by  Proposition \ref{prop_II.1.2_01}(ii),
the corresponding linear relation $\Theta_1$ in \eqref{IV.1.1_12A}
is not an operator, i.e.  has a nontrivial multivalued part,
$\mul\Theta_1:=\{f\in\cH:\ \{0,f\}\in\Theta_1\}\neq \{0\}$.

Let $f\in W^{2,2}_{\comp}(\cI\setminus X)$. Then $\Gamma_0^1f,\Gamma_1^1f\in l_{2,0}$ and
$f\in\dom(\rH_{X,\gA})$ if and only if $C_{X,\gA}\Gamma_1f=D_{X,\gA}\Gamma_0f$, where
     \begin{gather}\label{CD_delta_01}
C_{X,\gA} := CR_X,\qquad D_{X,\gA} := (D_\gA-CQ_X)R_X^{-1},\\
C := \left(%
\begin{array}{cccccc}
  0& 0 & 0 & 0& 0  & \dots\\
  0 & 0& 0& 0 & 0&  \dots\\
  0 & -1 & 1 &   0& 0& \dots\\
  0 & 0 & 0 & 0 & 0&  \dots\\
  0 & 0 & 0 & -1 & 1&  \dots\\
\dots& \dots&\dots&\dots&\dots&\dots\\
 \end{array}%
\right),
           \quad D_\gA := \left(%
\begin{array}{cccccc}
  1& 0 & 0 & 0& 0  & \dots\\
  0 & 1& 1& 0 & 0&  \dots\\
  0 & 0 & \alpha_1 &   0& 0& \dots\\
  0 & 0 & 0 & 1 & 1&  \dots\\
  0 & 0 & 0 & 0 & \alpha_2&  \dots\\
\dots& \dots&\dots&\dots&\dots&\dots\\
 \end{array}%
\right),\label{CD_delta_02}
\end{gather}
and $R_X=\oplus_{n=1}^\infty R_n,\ Q_X=\oplus_{n=1}^\infty Q_n$ are defined by \eqref{IV.1.1_09B}.

Define a linear relation $\Theta_1^0$ by
\begin{equation}\label{V.1_theta0}
\Theta_1^0=\{\{f,g\} \in l_{2,0}\oplus l_{2,0}:\
D_{X,\gA}f=C_{X,\gA}g\}.
\end{equation}
Hence we obviously get
\begin{equation}\label{IV.2.1_10}
\rH_{X,\gA}^0=\rH_{\min}^*\lceil \dom(\rH_{X,\gA}^0),\qquad
\dom(\rH_{X,\gA}^0)=\{f\in W^{2,2}_{\comp}(\cI\setminus X):\
\{\Gamma_0^1f,\Gamma_1^1f\}\in\Theta_1^0\}.
       \end{equation}
Straightforward calculations show that $\Theta_1^0$ is symmetric.
Moreover, \eqref{IV.2.1_10} implies that the closure of
$\Theta_1^0$ is $\Theta_1$. Hence $\Theta_1$ is a closed symmetric
linear relation. 
Therefore (see Subsection \ref{sss_II.1.1_lr}), 
$\Theta_1$ admits the  representation
    \begin{equation}\label{theta_1_rep}
\Theta_1= \Theta_1^{\op}\oplus \Theta_1^\infty, \quad
\cH=\cH_{\op}\oplus \cH_\infty,\quad \cH_{\op}=\overline{\dom
(\Theta_1)}=\overline{\dom (\Theta_1^{\op})}, \quad \cH_\infty :=
\mul\Theta_1,
      \end{equation}
where $\Theta_1^{\op}(\in\mathcal{C}(\cH_{\op}))$ is the operator
part of $\Theta_1.$
Moreover, it follows from \eqref{CD_delta_01}  that                
\begin{equation}\label{theta_1_rep'}
\mul\Theta_1=\ker (C_{X,\gA})=\overline{R_X^{-1}(\ker C)}, \qquad
\Theta_1^\infty=\{\{0,f\}:\ f\in \mul\Theta_1\}.
   \end{equation}
Since  
$\cH_{\op}=\overline{\ran (R_XC^*)},$ the system
$\{\mathbf{f}_n\}_{n=1}^\infty, $\ $\mathbf{f}_n:=
\frac{\sqrt{\gd_{n}}\mathrm{e}_{2n}-\sqrt{\gd_{n+1}}\mathrm{e}_{2n+1}}{\sqrt{\gd_{n}+\gd_{n+1}}}$,
forms the orthonormal basis in $\cH_{\op}$. Next we show that the
operator part $\Theta_1^{\op}$  of $\Theta_1$ is unitarily
equivalent to the minimal Jacobi operator
         \begin{equation}\label{IV.2.1_05}
B_{X,\gA}=\left(\begin{array}{cccc}
r_1^{-2}\bigl(\alpha_1+\frac{1}{\gd_1}+\frac{1}{\gd_2}\bigr) & -(r_1r_2\gd_2)^{-1} & 0&   \dots\\
-(r_1r_2\gd_2)^{-1} &r_2^{-2}\bigl(\alpha_2+\frac{1}{\gd_2}+\frac{1}{\gd_3}\bigr) & -(r_2r_3\gd_3)^{-1} &  \dots\\
0 & -(r_2r_3\gd_3)^{-1} & r_3^{-2}\bigl(\alpha_3+\frac{1}{\gd_3}+\frac{1}{\gd_4}\bigr)&  \dots\\
\dots & \dots& \dots & \dots
\end{array}\right) ,
    \end{equation}
where
$r_n:=\sqrt{\gd_{n}+\gd_{n+1}}, \ n\in\N$.
Observe first that
\begin{equation}\label{IV.2.1_04}
B_{X,\gA}=\widetilde{R}_X^{-1}(B_X+\mathcal{A}_\gA)\widetilde{R}_X^{-1}, \quad
\text{where}\quad 
         \end{equation}
        \begin{equation}\label{IV.2.1_03}
\widetilde{R}_X=\diag(r_n),\quad \mathcal{A}_\gA :=
\diag(\alpha_n),\quad B_X = \left(\begin{array}{cccc}
\frac{1}{\gd_1}+\frac{1}{\gd_2} & -\frac{1}{\gd_2} & 0&   \dots\\
-\frac{1}{\gd_2} &\frac{1}{\gd_2}+\frac{1}{\gd_3} & -\frac{1}{\gd_3}&   \dots\\
0 & -\frac{1}{\gd_3} & \frac{1}{\gd_3} + \frac{1}{\gd_4}&  \dots\\
\dots & \dots & \dots&  \dots
        \end{array}\right).
\end{equation}
Further, let us show that $\{\mathbf{f}_n\}_{n=1}^\infty \subset \dom(\Theta_1^{\op})$. Assume that there exists $\mathbf{g}_n$ such that $\{\mathbf{f}_n,\mathbf{g}_n\}\in \Theta_1^{\op}$, i.e., $\mathbf{g}_n=\Theta_1^{\op}\mathbf{f}_n$. The latter yields $\mathbf{g}_n\in \cH_{\op}$ and hence $\mathbf{g}_n=\sum_{k=1}^\infty g_{n,k} \mathbf{f}_k$. Moreover, after straightforward calculations we obtain
\begin{gather}
D_{X,\gA}\mathbf{f}_1=r_1^{-1}\bigl(-(\gA_1+\gd_1^{-1}+\gd_{2}^{-1})\mathrm{e}_{3}+\gd_{2}^{-1}\mathrm{e}_{5}\bigr),\notag\\
D_{X,\gA}\mathbf{f}_n=r_n^{-1}\bigl(\gd_n^{-1}\mathrm{e}_{2n-1}-(\gA_n+\gd_n^{-1}+\gd_{n+1}^{-1})\mathrm{e}_{2n+1}+\gd_{n+1}^{-1}\mathrm{e}_{2n+3}\bigr),\quad n\geq 2 \notag\\
C_{X,\gA}\mathbf{g}_n=-\sum_{k=1}^\infty g_{n,k}r_k\mathrm{e}_{2k+1},\quad n\geq 1.\notag 
\end{gather}
Hence $\{\mathbf{f}_n,\mathbf{g}_n\}\in\Theta$, i.e., equality $D_{X,\gA}\mathbf{f}_n=C_{X,\gA}\mathbf{g}_n$ holds, if and only if
\[
g_{n,n-1}=-\frac{1}{\gd_{n}r_{n-1}r_{n}},\quad
g_{n,n}=\frac{1}{r_{n}^2}\bigl(\alpha_n+\frac{1}{\gd_n}+\frac{1}{\gd_{n+1}}\bigr),\quad
g_{n,n+1}=-\frac{1}{\gd_{n+1}r_{n}r_{n+1}},\qquad n\geq 2,
\]
and $g_{n,k}=0$ for all $k\notin\{n-1,n,n+1\}$. Hence $\mathbf{f}_n \in \dom(\Theta_1^{\op})$ and in the basis $\{\mathbf{f}_n\}_{n=1}^\infty$ the matrix representation of the operator $\Theta_1^{\op}$ coincides with the matrix $B_{X,\gA}$ defined by \eqref{IV.2.1_05}.
Since the operator $B_{X,\gA}$ of
the form \eqref{IV.2.1_02} and \eqref{IV.2.1_05} is closed, we conclude that $\Theta_1^{\op}$ and $B_{X,\gA}$ are unitarily equivalent.

Let us summarize the above considerations in the following
proposition.
       \begin{proposition}\label{prop_IV.2.1_02}
Let $\Pi^1=\{\cH,\Gamma_0^1,\Gamma_1^1\}$ be the boundary triplet
constructed in Theorem \ref{th_bt_2} and let the linear relation
$\Theta_1$ be defined by \eqref{IV.1.1_12A}. Then $\Theta_1$
admits  representation \eqref{theta_1_rep}, where the "pure"
relation $\Theta_1^\infty$ is determined by \eqref{theta_1_rep'}
and \eqref{CD_delta_02}, and the operator part $\Theta_1^{\op}$ is
unitarily equivalent to the minimal Jacobi operator $B_{X,\gA}$ of
the form \eqref{IV.2.1_02} and \eqref{IV.2.1_05}.
\end{proposition}

\subsection{Self--adjontness}\label{sss_IV.2.2_sa_delta}

\textbf{1.} \ We begin with a result that reduces the property of
$\rH_{X,\gA}$ to be self--adjoint to that of the corresponding Jacobi
matrices $B_{X,\gA}$.
     \begin{theorem}\label{th_delta_sa}
The operator $\rH_{X,\gA}$ has equal deficiency indices and
$\mathrm{n}_+(\rH_{X,\gA})=\mathrm{n}_-(\rH_{X,\gA})\leq 1$.
Moreover, $\mathrm{n}_\pm(\rH_{X,\gA})=\mathrm{n}_\pm(B_{X,\gA})$,
where $B_{X,\gA}$ is the minimal operator associated with the Jacobi  matrix either \eqref{IV.2.1_01} or
\eqref{IV.2.1_05}. In particular, $\rH_{X,\gA}$ is self-adjoint if
and only if $B_{X,\gA}$ is.
      \end{theorem}
 \begin{proof}
Combining  Theorem \ref{th_IV.1.1_02} $(i)$ with Propositions
\ref{prop_IV.2.1_01} and \ref{prop_IV.2.1_02}, we arrive at
the equality
$\mathrm{n}_\pm(\rH_{X,\gA})=\mathrm{n}_\pm(B_{X,\gA}).$ It remans
to note that
for Jacobi matrices $\mathrm{n}_{\pm}(B_{X,\gA})\leq1$ (see
\cite{Akh, Ber68}).
       \end{proof}
The following result is immediate from Theorem  \ref{th_delta_sa}
though  we don't know its direct proof.
    \begin{corollary}
Let $B^{(1)}_{X,\gA}$ and $B^{(2)}_{X,\gA}$ be the minimal Jacobi
operators associated with \eqref{IV.2.1_05} and
\eqref{IV.2.1_01}, respectively. Then
$\mathrm{n}_\pm(B^{(1)}_{X,\gA})=
\mathrm{n}_\pm(B^{(2)}_{X,\gA}).$  In particular,
$B^{(1)}_{X,\gA}$ is self-adjoint if and only if so is
$B^{(2)}_{X,\gA}.$
        \end{corollary}
     \begin{remark}\label{rem_ext_jacobi}
 It was found out by Shubin Christ and Stolz  \cite{Chr_Sto_94} that
the operator $\rH_{X,\gA}$ may be symmetric with
$\mathrm{n}_{\pm}(H_{X,\gA})=1$ even if $\cI=\R_+$. In this case
the set of self-adjoint extensions of $\rH_{X,\gA}$ can be
described in terms of the classical Sturm--Liouville theory (for detail
see \cite{bsw}). Theorem \ref{th_delta_sa} enables us to describe
self-adjoint extensions of $\rH_{X,\gA}$ in a different way. More precisely,
consider the boundary triplet $\Pi^2$ defined in Theorem \ref{th_bt_1}. By Theorem \ref{th_delta_sa}, $\rH_{X,\gA}$ is symmetric if and only if the Jacobi operator
$B_{X,\gA}$ of the form \eqref{IV.2.1_01}--\eqref{IV.2.1_02} is also symmetric. By Proposition \ref{prop_II.1.2_01}, the mapping 
\begin{equation*}
\widetilde{B}_{X,\gA}\to \rH_{\widetilde{B}_{X,\gA}}:=\rH_{\min}^*\lceil\dom\rH_{\widetilde{B}_{X,\gA}},
\quad \dom\rH_{\widetilde{B}_{X,\gA}}:=\ker (\Gamma_1^2-\widetilde{B}_{X,\gA}\Gamma_0^2)
\end{equation*}
establish a bijective correspondence between the sets of self-adjoint extensions of $B_{X,\gA}$ and $\rH_{X,\gA}$.
\end{remark}
Using various criteria of self-adjointness of Jacobi matrices (see
e.g. \cite{Akh, Ber68, KosMir99, KosMir01}), we obtain necessary
and sufficient conditions for the operator $\rH_{X,\gA}$ to be
self-adjoint (symmetric) in $L^2(\cI)$. We emphasize that
different parameterizations \eqref{IV.2.1_01} and
\eqref{IV.2.1_05} of $\rH_{X,\gA}$ lead to different criteria.
     \begin{proposition}\label{cor_delta_carleman}
The Hamiltonian $\rH_{X,\gA}$ is self-adjoint for any
$\gA=\{\alpha_n\}_{n=1}^\infty\subset\R$ whenever
       \begin{equation}\label{IV.2.2_01}
\sum_{n=1}^\infty \gd_n^2=\infty.
\end{equation}
      \end{proposition}
       \begin{proof}
Let $B_{X,\gA}$ be the minimal Jacobi operator of the form
\eqref{IV.2.1_01'}, \eqref{IV.2.1_02}. By Carleman's theorem
\cite{Akh}, \cite[Chapter VII.1.2]{Ber68}, $B_{X,\gA}$ is
self-adjoint provided that
\begin{equation}\label{IV.2.2_02}
\sum_{n=1}^\infty (\gd_n^2+\gd_n^{3/2}\gd_{n+1}^{1/2})=\infty.
\end{equation}
Clearly,
$\gd_n^2<\gd_n^2+\gd_n^{3/2}\gd_{n+1}^{1/2}\leq\frac{7}{4}\gd_n^2
+\frac{1}{4}\gd_{n+1}^2$ and hence relations  \eqref{IV.2.2_01}
and \eqref{IV.2.2_02} are equivalent.

One completes the proof by
applying  Theorem \ref{th_delta_sa}.
   \end{proof}
If $\limsup_n \gd_n>0$, then condition \eqref{IV.2.2_01} is obviously
satisfied and Proposition \ref{cor_delta_carleman} yields the
following improvement of the result of Gesztesy and Kirsch (cf.
\cite[Theorem 3.1]{Ges_Kir_85}).
      \begin{corollary}[\cite{Ges_Kir_85}]\label{rem_IV.2.2_01}
If $\limsup_{n} \gd_n>0$ (in particular, $\gd_* =\liminf_{n}\gd_n
>0$), then $\rH_{X,\gA}$ is self-adjoint.
     \end{corollary}
In fact, Gesztesy and Kirsch \cite{Ges_Kir_85} established
self-adjointness for the operator $\rH_{X,\gA,q}$ (see
\eqref{I_01}--\eqref{I_03}) for a wide class of unbounded
potentials assuming only  $d_*>0$ . Note also that under assumption $\gd_*>0$ Corollary
\ref{rem_IV.2.2_01} was reproved  by Kochubei \cite{Koc_89} in the
framework of boundary triplets approach.

\textbf{2.} \ If  $\mathcal{I}=\R_+$ and
condition \eqref{IV.2.2_01} is violated, then  the operator
$\rH_{X,\gA}$ might be symmetric with nontrivial deficiency
indices $\mathrm{n}_\pm(H_{X,\gA})=1$. In particular,
this is the case when $\mathcal{I}=\R_+,\ \gd_n=1/n,$ and
$\alpha_n=-(2n+1)$ (see \cite[Remark on pp.
495--496]{Chr_Sto_94}). Our next result is partially inspired by
the example of C. Shubin Chtist and G. Stolz, and it also shows that Proposition
\ref{cor_delta_carleman} is sharp.
        \begin{proposition}\label{cor_delta_ber}
Let $\{\gd_n\}_{n=1}^\infty\in l_2,$\   
$\gd_n\ge 0,$  and
      \begin{equation}\label{IV.2.2_03}
\gd_{n-1}\gd_{n+1}\geq \gd_n^2,\quad n\in\N.
      \end{equation}
If, in addition,  the strengths $\alpha_n$ of
$\delta$-interactions satisfy
\begin{equation}\label{IV.2.2_05}
\sum_{n=1}^\infty \gd_{n+1}\left|\alpha_n+\frac{1}{\gd_{n}}+\frac{1}{\gd_{n+1}}\right|<\infty,
\end{equation}
then the operator $\rH_{X,\gA}$ is symmetric with $\mathrm{n}_\pm(\rH_{X,\gA})=1$.
\end{proposition}
              \begin{proof}
Consider the Jacobi  matrix \eqref{IV.2.1_05}. To apply \cite[Theorem 1]{KosMir01} we denote $a_n:=r_n^{-2}(\gA_n+1/\gd_n+1/\gd_{n+1})$
and  $b_n:=(r_nr_{n+1}d_{n+1})^{-1}$, $n\in\N$, and define the
sequence $\{c_n\}_{n=1}^{\infty}$ as follows
\[
c_1:=b_1,\qquad c_2:=1,\qquad c_{n+1}:=-\frac{b_{n-1}}{b_n}c_{n-1}, \qquad n\in \N.
\]
It is easily seen that
\[
c_{n+1}=(-1)^{n+1}r_{n+1}\frac{\gd_{n+1}\ \gd_{n-1}\cdot\dots}{\gd_{n}\ \gd_{n-2}\cdot\dots}\cdot \widetilde{c},\quad n\in\N;
\qquad \widetilde{c}:=\left\{\begin{array}{cc}
                                           {c_1}{r^{-1}_1}, & n=2k+1,\\
                                           {c_2}{r^{-1}_2}, & n=2k.
                                            \end{array}\right.
\]
Due to \eqref{IV.2.2_03}, we obtain
\begin{equation}\label{IV.2.2_06}
\frac{\gd_{n+1}\ \gd_{n-1}\cdot\dots}{\gd_{n}\
\gd_{n-2}\cdot\dots}=
\sqrt{\gd_{n+2}}\frac{\gd_{n+1}}{\sqrt{\gd_{n+2}\gd_{n}}}\frac{\gd_{n-1}\cdot\dots}{\sqrt{\gd_{n}\gd_{n-2}}\cdot\dots}\leq
C\sqrt{\gd_{n+2}}, \qquad n\in \N.
\end{equation}
Therefore,
    \[
 |c_{n+1}|\leq
\widetilde{c}Cr_{n+1}\sqrt{\gd_{n+2}},
     \]
and hence $\{c_n\}_{n=1}^\infty\in l_2$.
On the other hand, it follows from \eqref{IV.2.2_05} and
\eqref{IV.2.2_06} that $\sum_{n=1}^\infty |a_n|c_n^2<\infty$. By
\cite[Theorem 1]{KosMir01}, this inequality together with the
inclusion  $\{c_n\}_{n=1}^\infty\in l_2$ yields
$\mathrm{n}_{\pm}(B_{X,\gA})=1.$
It remains to apply Theorem \ref{th_delta_sa}.
         \end{proof}
         \begin{remark}\label{rem_mih_sa}
         Note that in the case $\cI=\R_+$ the self-adjointness of $\rH_{X,\gA}$ for arbitrary $\gA\subset \R$ was erroneously stated in \cite{Mih_93, Mih_94a}.
         \end{remark}
Let us present sufficient conditions for self-adjointness in the
case when \eqref{IV.2.2_01} does not hold.
\begin{proposition}\label{prop_IV.2.2_03}
Assume that \eqref{IV.2.2_01} does not hold. Let also
$\gA=\{\gA_n\}_{n=1}^\infty$ and $X=\{x_n\}_{n=1}^\infty$ satisfy
one of the following conditions: \item $(i)$
      \begin{equation}\label{5.20}
\sum_{n=1}^\infty
|\alpha_n|\gd_n\gd_{n+1}r_{n-1}r_{n+1}=\infty,\qquad
r_n=\sqrt{d_n+d_{n+1}}.
     \end{equation}
\item $(ii)$ There exists a positive constant $C_1>0$ such that
\begin{equation}\label{IV.2.2_09}
\alpha_n+\frac{1}{\gd_n}\left(1+\frac{r_{n}}{r_{n-1}}
\right)+\frac{1}{\gd_{n+1}}\left(1+\frac{r_{n}}{r_{n+1}}
\right)\leq C_1(\gd_n+\gd_{n+1}),\qquad n\in\N.
\end{equation}
\item $(iii)$ There exists a positive constant $C_2>0$ such that
\begin{equation}\label{IV.2.2_10}
\alpha_n+\frac{1}{\gd_n}\left(1-\frac{r_{n}}{r_{n-1}}
\right)+\frac{1}{\gd_{n+1}}\left(1-\frac{r_{n}}{r_{n+1}}
\right)\geq -C_2(\gd_n+\gd_{n+1}) 
,\qquad n\in\N.
\end{equation}
Then the operator $\rH_{X,\gA}$ is self-adjoint in $L^2(\cI)$.
        \end{proposition}
                \begin{proof}
$(i)$ Since $\{\gd_n\}_{n=1}^\infty\in l_2$, we get
$\sum_{n=1}^\infty (\gd_n + \gd_{n+1})r_{n-1}r_{n+1}<
C\sum_{n=1}^\infty \gd^2_n < \infty$.
Applying the Dennis-Wall test (\cite[p.25, Problem 2]{Akh}) to
 matrix \eqref{IV.2.1_05}, we obtain that \eqref{5.20} yields
self-adjointness of the minimal operator $B_{X,\alpha}$ associated
with \eqref{IV.2.1_05}. By Theorem \ref{th_delta_sa}, $\rH_{X,\gA}
= \rH_{X,\gA}^*$.

$(ii)-(iii)$ Applying \cite[Theorem VII.1.4]{Ber68} (see
also \cite[Problem 3, p.37]{Akh}) to the Jacobi matrix
\eqref{IV.2.1_05},  we obtain that conditions \eqref{IV.2.2_09} and \eqref{IV.2.2_10} guarantee self-adjointness of $B_{X,\alpha}$. Theorem  \ref{th_delta_sa} completes the proof.
          \end{proof}
Conditions $(i)$--$(iii)$ show that if $\rH_{X,\gA}$ is
self-adjoint, then the coefficients $\alpha_n$ cannot tend to
$\infty$ very fast. Let us demonstrate this by considering an
example.

      \begin{example}\label{example_IV.2.2_01}
Let $\mathcal{I}=\R_+,\ x_0=0, \ x_{n}-x_{n-1}=\gd_n:=1/n,\ n\in\N $. Consider the operator
\begin{equation}\label{IV.2.2_11}
\rH_{A}:=-\frac{\rD^2}{\rD x^2}+\sum_{n=1}^\infty\alpha_n\delta(x-x_n).
      \end{equation}
Clearly,  $\{\gd_n\}_{n=1}^\infty\in l_2$, i.e.,
condition \eqref{IV.2.2_01} is violated. Applying Propositions
\ref{cor_delta_ber} and \ref{prop_IV.2.2_03}, after
straightforward calculations we obtain:
\item  $(i)$ \ \ \ If $\sum_{n=1}^\infty
\frac{|\alpha_n|}{n^3}=\infty$, then the  operator $\rH_{A}$  is
self-adjoint (cf. Proposition \ref{prop_IV.2.2_03} $(i)$).\ \item
$(ii)$ \ \ If $\alpha_n\leq -4\bigl(n+\frac12\bigr)+O(n^{-1})$,
then $\rH_{A}$ is self-adjoint (cf. Proposition
\ref{prop_IV.2.2_03} $(ii)$). \item $(iii)$ \ If $\alpha_n\geq
-\frac{C}{n},\ n\in \N,\ C\equiv const>0$, then $\rH_{A}$ is
self-adjoint (cf. Proposition \ref{prop_IV.2.2_03} $(iii)$).
\item  $(iv)$ \ If $\alpha_n=-2n-1+O(n^{-\varepsilon})$ with some $ \varepsilon>0$, then $\mathrm{n}_\pm(\rH_A)=1$ (cf. Proposition \ref{cor_delta_ber}).
\end{example}
Conditions $(ii)$ and $(iii)$ show that there is a gap
between conditions of self-adjointness.  Moreover, $(iii)$ shows
that for the case of positive interactions $\alpha_n$ the operator
$\rH_{A}$ is self-adjoint.
We can extend $(iv)$ as follows.
               \begin{proposition}\label{prop_IV.2.2_04}
Let the Hamiltonian $\rH_{A}$ be the same as in Example
\ref{example_IV.2.2_01}. If
        \begin{equation}\label{alpha_pi}
\alpha_n=a\left(n+\frac12\right)+O(n^{-1}),\qquad a\in(-4,0),
     \end{equation}
then the  operator $\rH_{A}$   
is symmetric with $\mathrm{n}_{\pm}(\rH_{A})=1$.
      \end{proposition}
       \begin{proof}
Define the sequence
\begin{equation}\label{5.25}
\widetilde{r}_{n+1}:=\frac{d_{n+1}}{\widetilde{r}_n},\qquad \widetilde{r}_1:=1,\qquad d_n=\frac1n,\qquad n\in\N.
\end{equation}
Then
\begin{equation}\label{5.26}
\widetilde{r}_{n+1}=\frac{n(n-2)\cdot\dots}{(n+1)(n-1)\cdot\dots}=\frac{n!!}{(n+1)!!}
\end{equation}
Let us estimate $\widetilde{r}_n$. 
Observe that
\[
(2k-1)!!=2^k\frac{\Gamma(k+\frac12)}{\Gamma(1/2)},\qquad (2k)!!=2^k\Gamma(k+1),
\]
where $\Gamma(\cdot)$ is the classical $\Gamma$-function. Using the asymptotic of $\Gamma(\cdot)$, we get
\begin{equation}\label{5.27}
(4k+1)\widetilde{r}^2_{2k}
=\frac{4}{\pi}\bigl(1+O(k^{-2})\bigr),\qquad (4k+3)\widetilde{r}_{2k+1}^2
=\pi\bigl(1+O(k^{-2})\bigr),\qquad k\to\infty.
\end{equation}
Indeed, consider the first equality in \eqref{5.27}. Since
$\Gamma(1/2)=\sqrt \pi$ and
\begin{gather}
\Gamma(k)=\sqrt{2\pi}e^{-k}k^{k-1/2}\left(1+\frac{1}{12k}+O(k^{-2})\right),\notag\\
\left(1+\frac{1}{k}\right)^k=\mathrm{e}\left(1-\frac{1}{2k}+O(k^{-2})\right)\notag,
\end{gather}
we obtain
\begin{eqnarray}
(4k+1)\widetilde{r}^2_{2k}=(4k+1)\frac{\Gamma(k+1/2)^2}{\pi\Gamma(k+1)^2}=(4k+1)\frac{\mathrm{e}}{\pi}\frac{(k+1/2)^{2k}(1+\frac{1}{6(2k+1)}+O(k^{-2}))^2}{(k+1)^{2k+1}(1+\frac{1}{6(2k+2)}+O(k^{-2}))^2}\notag\\
=\frac{\mathrm{e}}{\pi}\frac{4k+1}{k+1/2}\left(1+\frac{1}{2k+1}\right)^{-(2k+1)}\left(1+O(k^{-2})\right)=\frac{4}{\pi}\left(1+O(k^{-2})\right), \qquad k\to\infty.
\end{eqnarray}
Further, define $\alpha^0:= \{\gA_n^0\}_{n=1}^\infty$
by setting
\[
\alpha_n^0:=\left\{\begin{array}{cc}
                           -(4k+1)+\frac{4}{\pi}\bigl(1+\frac{a}{2}\bigr)\widetilde{r}_{n}^{-2},& n=2k,\\
                           -(4k+3)+\pi\bigl(1+\frac{a}{2}\bigr)\widetilde{r}_{n}^{-2},& n=2k+1.
                           \end{array}
                                  \right.
\]
Clearly, by \eqref{5.27}, $\alpha^0$ satisfies \eqref{alpha_pi}.
Moreover, for this choise of $\gA^0$ we get
\[
B_{X,\gA^0} + \mathcal{A}_{\gA^0} = \widetilde{R}_1^{-1}J_{a}\widetilde{R}_1^{-1},\quad 
\]
where  $B_{X,\gA^0}$ is defined by \eqref{IV.2.1_05},\
 $\mathcal{A}_{\gA^0} = \diag(\alpha^0_n),$  and
\[
\widetilde{R}_1:=\diag (\widetilde{r}_n),\quad \text{and}\quad  J_{a} :=\left(%
\begin{array}{ccccc}
  \frac{4}{\pi}\bigl(1+\frac{a}{2}\bigr)& 1 & 0 & 0&  \dots\\
  1 & \pi\bigl(1+\frac{a}{2}\bigr)& 1& 0 &   \dots\\
  0 & 1 & \frac{4}{\pi}\bigl(1+\frac{a}{2}\bigr) &   1&  \dots\\
  0 & 0 & 1 & \pi\bigl(1+\frac{a}{2}\bigr) &   \dots\\
\dots& \dots&\dots&\dots&\dots\\
 \end{array}%
\right).
\]
The Floquet   determinant (see, for instance, \cite[\S
7.1]{Tes_98}) of the  peridic Jacobi matrix
 $J_{a}$ is
$\Delta_a(\lambda)=-2+(\lambda-\frac{4}{\pi}\bigl(1+\frac{a}{2}\bigr))(\lambda-\pi\bigl(1+\frac{a}{2}\bigr))$.
Note that all solutions of $\tau_{a}f=0$ are bounded if
$|\Delta_a(0)|< 2$ (here $\tau_{a}$ is a difference expression
associated with the matrix $J_{a}$). The latter is equivalent to
the inequality $0<|1+\frac{a}{2}|<1$. Moreover, all solutions of
$\tau_{-2}f=0$ are bounded too.
Therefore, all solutions of $\tau_{a}f=0$ are bounded if
\[
|2+a|<2.
\]
Furthermore, $g$ solves $\tau_{X,\gA}y=0$ precisely when $\widetilde{R}_X \widetilde{R}_{1}g$ solves $\tau_{a}f=0$. By \eqref{5.26}--\eqref{5.27} and \eqref{IV.2.1_03}, we get $\{r_n\widetilde{r}_n\}_{n\in\N}\in l_2$. Hence all solutions of the equation $\tau_{X,\gA}y=0$ are $l_2$ solutions, that is 
the operator $B_{X,\gA^0}$ is symmetric  with $\mathrm{n}_{\pm}(B_{X,\gA^0})=1$.
Since bounded perturbations do not change the deficiency indices of $B_{X,\gA}$, we complete the proof by applying Theorem \ref{th_IV.1.1_02} $(i)$.
\end{proof}

\subsection{Resolvent comparability}\label{sec_delta_rescom}

Let us fix  $X=\{x_n\}_1^\infty\subset\cI$ and consider
Hamiltonians
$\rH_{X,\gA_1}$ and $\rH_{X,\gA_2}$   
corresponding  the  strengths $\gA_1=\{\gA_n^{(1)}\}_{n=1}^\infty$
and $\gA_2=\{\gA_n^{(2)}\}_{n=1}^\infty$, respectively.
       \begin{proposition}\label{th_delta_res}
Suppose $\rH_{X,\gA_1}$ and $\rH_{X,\gA_2}$  are self-adjoint and
${B}_{X,\gA_1}$ and ${B}_{X,\gA_2}$  the corresponding
(self-adjoint) Jacobi operators  defined either by
\eqref{IV.2.1_01} or \eqref{IV.2.1_05}. Then for any
$z\in\rho(\rH_{X,\gA_1})\cap\rho(\rH_{X,\gA_2})$ and
$p\in(0,\infty)\cup\{\infty\}$ the inclusion
     \begin{equation}\label{IV.2.3.01A}
(\rH_{X,\gA_1}-z)^{-1}-(\rH_{X,\gA_2}-z)^{-1}\in \mathfrak{S}_p
     \end{equation}
 is equivalent to the inclusion
      \begin{equation}\label{IV.2.3.01}
({B}_{X,\gA_1}-\I)^{-1}-({B}_{X,\gA_2}-\I)^{-1}\in \mathfrak{S}_p.
     \end{equation}
          \end{proposition}
        \begin{proof}
Combining Theorem \ref{th_IV.1.1_02} with Proposition
\ref{prop_IV.2.1_02}, we get  
the result with ${B}_{X,\gA_j}$ defined by
\eqref{IV.2.1_05}.
The result  with the matrices defined by \eqref{IV.2.1_01} is
implied by combining Proposition \ref{prop_IV.2.1_01}
with  Remark  \ref{rem_IV.1.1_01}.             
              \end{proof}

Next we present simple sufficient condition.
      \begin{corollary}\label{col_rc_2}
If
$\left\{\frac{\alpha_n^{(1)}-\alpha_n^{(2)}}{\gd_{n+1}}\right\}_{n=1}^\infty
\in l_p, \ p\in (0,\infty)$ ($\in c_0$, $p=\infty$), 
 then
inclusion \eqref{IV.2.3.01A} holds.
     \end{corollary}
    \begin{proof}
Clearly,  $l_{2,0}\subset \dom(B_{X,\gA_1})\cap
\dom(B_{X,\gA_2}).$ On the other hand, for any $f\in l_{2,0}$
\eqref{IV.2.2_01A} yields
\[
B_{X,\gA_2}f-B_{X,\gA_1}f=R^{-1}_X\bigl(\widetilde{B}_{\gA_1} -
\widetilde{B}_{\gA_2}\bigr)R^{-1}_Xf=
\oplus_{n=1}^\infty\left(\begin{array}{cc}
                                   \frac{\alpha_n^{(1)}-\alpha_n^{(2)}}{\gd_{n+1}}&0\\
                                                            0&                    0
                             \end{array}\right)f.
\]
Hence and due to the assumption,
$\overline{B_{X,\gA_2}-B_{X,\gA_1}}\in \mathfrak{S}_p\subset
[\cH]$ and  $\dom(B_{X,\gA_1}) = \dom(B_{X,\gA_2}).$ 
It remains to apply Proposition   \ref{prop_II.1.4_02}.
     \end{proof}
 In the case $d_*>0$, the resolvent comparability criterion was obtained in \cite{Koc_89}
 (see also \cite{Mih_94a}).
We omit the corresponding proof, though it can be extracted from
Proposition  \ref{th_delta_res}.
      \begin{corollary}[\cite{Koc_89, Mih_94a}]\label{col_rc_1}
If\  $0<d_*\leq d^*<\infty$, then  \eqref{IV.2.3.01A} is
equivalent to the inclusion
    \begin{equation}\label{IV.2.3.02}
({\alpha_n^{(1)}-\I})^{-1} - ({\alpha_n^{(2)}-\I})^{-2}\in
l_p,\qquad
p\in(0,\infty),\qquad (\in
c_0,\quad\text{if}\quad p=\infty).
\end{equation}
Moreover, if $\{\alpha_n^{(j)}\}_{n=1}^\infty\in l_\infty$, then \eqref{IV.2.3.02} holds precisely when  $\{\alpha_n^{(1)}-\alpha_n^{(2)}\}_{n=1}^\infty\in l_p \ \ (\in c_0)$.
\end{corollary}

%
\subsection{Operators with discrete spectrum}\label{sec_delta_disc}
%

Combining the results of Section \ref{sss_IV.2.1_boun_op} with
Theorem \ref{th_IV.1.1_02}, we obtain the discreteness criterion
for the Hamiltonian $\rH_{X,\gA}$.
          \begin{theorem}\label{th_disc_d}
 Let $B_{X,\gA}$ be the minimal Jacobi operator defined either by
\eqref{IV.2.1_01} or \eqref{IV.2.1_05}.
\item $(i)$  If $\mathrm{n}_\pm(B_{X,\gA})=1$,  then any self-adjoint
extension of $\rH_{X,\gA}$ has discrete spectrum.
\item $(ii)$ If $B_{X,\gA} = B_{X,\gA}^*$, then the Hamiltonian
$\rH_{X,\gA} (= \rH_{X,\gA}^*)$ has discrete spectrum if and only
if
\begin{description}
 \item $\bullet$\quad $\lim_{n\to\infty}\gd_n=0$, and
 \item $\bullet$\quad  $B_{X,\gA}$ has discrete spectrum.
\end{description}
      \end{theorem}
\begin{proof}
$1)$ To be precise, let
$B_{X,\gA}$ be defined by \eqref{IV.2.1_01}. Since
$\mathrm{n}_{\pm}(B_{X,\gA})=1$, any self-adjoint extension of
$B_{X,\gA}$ has discrete spectrum (see \cite{Akh, Ber68}).
Moreover, by Corollary \ref{rem_IV.2.2_01}, $\lim_{n\to
\infty}\gd_n= 0$. Hence the operator $\rH_0$ defined by
\eqref{h_0} has discrete spectrum too. The Krein resolvent formula
\eqref{II.1.4_01} implies that any  self-adjoint extension of
$\rH_{X,\gA}$ is discrete.

$2)$ follows from Theorem \ref{th_IV.1.1_02} $(iv)$ and Remark
\ref{rem_IV.1.1_01}.
        \end{proof}

Next we present some sufficient conditions for self-adjoint
Hamiltonian $\rH_{X,\gA}$ to be discrete.
       \begin{proposition}\label{prop_IV.2.4_01}
Assume that 
the operator $B_{X,\gA}$ defined by \eqref{IV.2.1_01}--\eqref{IV.2.1_02} is self-adjoint and $\lim_{n\to\infty}\gd_n=0$. If
    \begin{equation}\label{prop_chihara_1}
\lim_{n\to\infty}\frac{|\gA_n|}{\gd_n}=\infty\qquad \text{and}\qquad\lim_{n\to\infty}\frac{1}{\gd_n\gA_n}>-\frac{1}{4},
\end{equation}
then the operator $\rH_{X,\gA}$ has discrete spectrum.
\end{proposition}
\begin{proof}
Applying \cite[Theorem 8]{Chi62} to the operator $B'_{X,\gA}$ of
the form \eqref{IV.2.1_01'}, we obtain that the spectrum of
$B'_{X,\gA}$ is discrete provided that $\lim_{n\to\infty}\gd_n=0$
and conditions \eqref{prop_chihara_1}
are satisfied. Theorem \ref{th_disc_d} completes the proof.
       \end{proof}
Proposition \ref{prop_IV.2.4_01} enables us to construct
Hamiltonians  $\rH_{X,\gA}$ with discrete spectrum, which is not
lower semibounded.
     \begin{example}\label{example_IV.2.4_01}
$(a)$ Let $\cI=\R_+, \ x_n=\sqrt{n},\ n\in\N$. Then $\gd_n=\frac{1}{\sqrt{n}+\sqrt{n+1}}\approx \frac{1}{2\sqrt{n}}$ and, by Proposition \ref{cor_delta_carleman}, the operator $\rH_{X,\gA}$ is self-adjoint for arbitrary $\gA=\{\gA_n\}_{n=1}^\infty\subset \R$. Consider the operator
\[
\rH_{\varepsilon}:=-\frac{\rD^2}{\rD x^2}+\sum_{n=1}^\infty n^{-\varepsilon}\ \delta(x-\sqrt{n}),\qquad \varepsilon\in (0,1/2).
\]
Clearly, conditions \eqref{prop_chihara_1} hold and hence the operator $\rH_{\varepsilon}$ is discrete if $\varepsilon \in(0,1/2)$.

$(b)$ Again, let $\cI=\R_+, \ x_n=\sqrt{n},\ n\in\N$.
Define $\alpha_n=-C\sqrt{n}, \ C\equiv const\in \R$. By Proposition \ref{prop_IV.2.4_01}, the operator
\[
\rH_{C}:=-\frac{\rD^2}{\rD x^2}-\sum_{n=1}^\infty C\sqrt{n}\ \delta(x-\sqrt{n}),
\]
has discrete spectrum if $C>8$. Moreover, the operator $\rH_C$ is
not lower semibounded since so is the operator considered in
Proposition \ref{prop_IV.2.5_01} (see below).
     \end{example}
     \begin{remark}\label{rem_mih_discr}
It was stated in \cite{Mih_94b} that the spectrum $\sigma(\rH_{X,\gA})$ of $\rH_{X,\gA}$ is not discrete whenever $\gA\in l^\infty$. However, Example \ref{example_IV.2.4_01} $(a)$ shows that $\sigma(\rH_{X,\gA})$ may be discrete even if $\lim_{n\to\infty}\gA_n=0$.
     \end{remark}
         \begin{proposition}\label{prop_IV.2.4_02}
Let 
the 
operator $B_{X,\gA}$ defined by
\eqref{IV.2.1_05} be self-adjoint and
$\lim_{n\to\infty}\gd_n=0$. If
    \begin{gather}
\lim_{n\to\infty}\frac{\left|\gA_n+1/\gd_n+1/\gd_{n+1}\right|}{\gd_n
+\gd_{n+1}}=\infty\notag,   \\
\lim_{n\to\infty}\left(\gA_n\gd_{n+1}+1+\frac{\gd_{n+1}}{\gd_{n}}\right)^{-1}
\left(\gA_{n+1}\gd_{n+1}+1+\frac{\gd_{n+1}}{\gd_{n+2}}\right)^{-1}<\frac{1}{4},
\label{prop_chihara_2}
     \end{gather}
then the operator $\rH_{X,\gA}$ has  discrete spectrum.
\end{proposition}

\begin{proof}
Applying  \cite[Theorem 8]{Chi62} to the Jacobi matrix $B_{X,\gA}$
of the form \eqref{IV.2.1_05} we get that $B_{X,\gA}$ is discrete.
Since $\lim_{n\to\infty}\gd_n=0$, by Theorem \ref{th_disc_d}  so
is $\rH_{X,\gA}.$
     \end{proof}
\begin{remark}
In the case $\lim_{n\to\infty}\frac{d_n}{d_{n+1}}=1$, Proposition
\ref{prop_IV.2.4_01} follows from Proposition
\ref{prop_IV.2.4_02}. Let us also note that the second of
conditions \eqref{prop_chihara_1}  (of conditions
 \eqref{prop_chihara_2}) is sharp. In
\cite{Szw03}, under additional mild assumptions on coefficients it
is shown that the operator $B_{X,\gA}$ has absolutely continuous
spectrum if the limit in \eqref{prop_chihara_1} is less than
$-\frac14$ $($ resp. greater than $\frac14$
$)$ and $\{\gd_n\}_{n\in\N}\notin l_2$.
        \end{remark}
         \begin{proposition}\label{prop_IV.2.4_03}
Assume that $\lim_{n\to\infty}\gd_n=0$ and
       \begin{equation}\label{prop_cojuhari}
\lim_{n\to\infty}\frac{1}{\bigl(\gd_n+\gd_{n+1}\bigr)}\left(\gA_n+\frac{1}{\gd_n}+\frac{1}{\gd_{n+1}}-\frac{r_{n-1}}{\gd_n r_{n}}-\frac{r_{n+1}}{\gd_{n+1} r_{n}}\right)=+\infty,
      \end{equation}
where $r_n=\sqrt{\gd_n+\gd_{n+1}}$. Then the operator $\rH_{X,\gA}$ is self-adjoint and has discrete spectrum.
     \end{proposition}
\begin{proof}
By Proposition \ref{prop_IV.2.2_03} $(iii)$, the operator $B_{X,\gA}$ defined by 
\eqref{IV.2.1_05} is self-adjoint. By \cite[Theorem
3.1]{CojJan07}, \eqref{prop_cojuhari} yields discreteness of
$B_{X,\gA}$. It remains to apply Theorem \ref{th_disc_d}.
      \end{proof}

\subsection{Semiboundedness}\label{sec_delta_sembd}

We start with general criterion of semiboundedness.
     \begin{theorem}\label{th5.22}
Let the minimal Jacobi operator $B_{X,\gA}$ be defined by
\eqref{IV.2.1_02} and \eqref{IV.2.1_05}. Then the operator
$\rH_{X,\gA}$ is lower semibounded 
if and only if $B_{X,\gA}$ is lower semibounded. 
     \end{theorem}
   \begin{proof}
According to \eqref{IV.1.1_12A} $\rH_{X,\gA}= \rH_{\Theta_1}$. By
Theorem \ref{th_IV.1.1_02} $(ii)$, the operator $\rH_{X,\gA}=
\rH_{\Theta_1}$
is lower semibounded if and only if
$\Theta_1$  is.  
It remains to note that by  Proposition \ref{prop_IV.2.1_02},
the operator part $\Theta_1^{\op}$ of $\Theta_1$  
is unitarily equivalent to the operator $B_{X,\gA}$ defined by
\eqref{IV.2.1_02} and\eqref{IV.2.1_05}.
         \end{proof}
Let us present several conditions for semiboundedness in terms of
$X=\{x_n\}_1^\infty$ and $\alpha=\{\alpha_n\}_1^\infty$. The following result has been
obtained in \cite{Bra85} using the form method.
        \begin{corollary}[\cite{Bra85}]\label{cor_IV.2.5_01}
Let $d_*>0$. Then the operator $\rH_{X,\gA}$ is lower semibounded if and only if
   \begin{equation}\label{sb_brasche}
\inf_{n\to\infty}\alpha_n>-\infty.
   \end{equation}
    \end{corollary}
\begin{proof}
Since $d_*>0$, the operators $B_X,\ R_X,\ R_X^{-1}$ in
\eqref{IV.2.1_04} are bounded. Therefore, $B_{X,\gA}$ is
semibounded if and only if so is  $\mathcal{A}_\gA,$ that is the
sequence $\gA=\{\gA_n\}_{n=1}^\infty.$
      \end{proof}
In the case $d_*=0$ the situation becomes  more complicated.
Indeed, condition \eqref{sb_brasche} is no longer necessary for
lower semiboundedness  
(see \cite[Example 2]{Bra85}). Moreover, we will show that
\eqref{sb_brasche} is no longer sufficient
(cf. \cite[Corollary 2]{Mih_93} 
where the opposite statement is announced). Moreover,
$\rH_{X,\gA}$ might be non-semibounded below even if
$\inf_{n\to\infty}\alpha_n=0$.

We begin with the following
sufficient condition.
   \begin{corollary}\label{cor_IV.2.5_02}
The Hamiltonian  $\rH_{X,\gA}$ is semibounded below whenever
        \begin{equation}\label{IV.2.5.07}
\inf_{n\to\infty}\frac{\alpha_n}{\gd_n+\gd_{n+1}}>-\infty,
   \end{equation}
         \end{corollary}
\begin{proof}
The matrix $B_X$ in \eqref{IV.2.1_03} admits the representation
$B_X=(I-U^*)D_X^{-1}(I-U)$, where $D_X := \diag(\gd_n)$ and $U$ is
unilateral shift in $l_2$. Hence $B_X$ is nonnegative, $B_X\geq
0$, and we get
\[
B_{X,\gA}=\widetilde{R}_X^{-1}(B_X +
\mathcal{A}_\gA)\widetilde{R}_X^{-1}\geq
\widetilde{R}_X^{-1}\mathcal{A}_\gA \widetilde{R}_X^{-1},
\]
Since  $\widetilde{R}_X=\diag(r_n)$ and  $\mathcal{A}_\gA =
\diag(\alpha_n)$ we obtain lower semiboundedness of $B_{X,\gA}$ by
combining the last inequality with condition \eqref{IV.2.5.07}.
Theorem \ref{th5.22} completes the proof.
       \end{proof}
           \begin{remark}\label{rem_IV.2.5_01}
In the case $\gd_*>0$, condition \eqref{IV.2.5.07} is equivalent to
\eqref{sb_brasche} and hence is also necessary for semiboundedness
of $\rH_{X,\gA}$.
If $\gd_*=0$, then \eqref{IV.2.5.07} is only sufficient (see
\cite[Example 2]{Bra85}).
          \end{remark}

Note that condition \eqref{IV.2.5.07} may be violated even if
$\alpha_n\to 0$. Next example shows that in this case the operator
$\rH_{X,\gA}$ might  be non-semibounded below.
\begin{proposition}\label{prop_IV.2.5_01}
Let $\cI=\R_+$ and $x_n=\sqrt{n}$. If $\alpha_n= -n^{-\varepsilon}$ with $\varepsilon\in[0,1/2)$, then the operator $\rH_{X,\gA}$ is self-adjoint and not semibounded below in $L^2(\cI)$.
\end{proposition}
\begin{proof}
Note that $\gd_n=\sqrt{n}-\sqrt{n-1}=\frac{1}{\sqrt{n-1}+\sqrt{n}}\asymp \frac{1}{2\sqrt{n}}$ as $n\to\infty$.
Hence, by Proposition \ref{cor_delta_carleman}, the operator $\rH_{X,\gA}$ is self-adjoint.

By Proposition \ref{prop_IV.2.1_02}, $\rH_{X,\gA} =\rH_{\Theta_1}$, where the
operator part $\Theta'_1$ of  $\Theta_1$  is unitarily equivalent
to the Jacobi matrix $B_{X,\gA}$ of the form \eqref{IV.2.1_05}. 
Clearly,
 $B_{X,\gA}$ admits the following representation
\[
B_{X,\gA}=\widetilde{R}_X^{-1}(B_X+\mathcal{A}_\gA)\widetilde{R}_X^{-1}=\widetilde{R}_X^{-1}\bigl[D_X^{-1/2}\bigl(J_{\per}+UK^2U^*+UK+KU^*+\widetilde{\mathcal{A}}_\gA\bigr)D_X^{-1/2}\bigr] \widetilde{R}_X^{-1},
\]
where $D_X=\diag(d_n)$,  $U$ is unilateral shift in $l_2$, and 
\[
J_{\per}=\left(%
\begin{array}{cccc}
  2& 1 & 0 &   \dots\\
  1 & 2& 1&    \dots\\
  0 & 1 & 2 &    \dots\\
\dots&\dots&\dots&\dots\\
 \end{array}%
\right),\quad K=\diag(k_n), \quad k_n:=\frac{\sqrt{d_{n-1}}}{\sqrt{d_n}}-1,\quad \widetilde{\mathcal{A}}_\gA=\mathcal{A}_\gA D_X=\diag(\widetilde{\gA}_n).
\]
Note that $\widetilde{\gA}_n=\gA_n\gd_n\asymp- n^{-(1/2+\varepsilon)}$ and $k_n=O(n^{-1})$ as $n\to\infty$.
Since $\varepsilon\in [0,1/2)$, the sum $\frac{1}{2}\widetilde{\gA}_n+k_n$ is negative for $n$ large enough.
Therefore, 
$\widetilde{B}_{X,\gA}:=\widetilde{R}_X^{-1}\bigl[D_X^{-1/2}\bigl(J_{\per}
+ \frac{1}{2}\widetilde{\mathcal{A}}_\gA\bigr)D_X^{-1/2}\bigr]
\widetilde{R}_X^{-1}$ is lower semibounded if so is $B_{X,\gA}$.

Let $f_N=(f_1,\dots,f_{2N},0,0,\dots)$, where $f_{2n}=1$, $f_{2n-1}=-1$, $n\in\{1,\dots,N\}$. Then
we get
\begin{gather*}
(\widetilde{\mathcal{A}}_\gA f_N,f_N)=\sum_{n=1}^{2N}\gA_n\gd_n=-\sum_{n=1}^{2N} \frac{n^{-\varepsilon}}{\sqrt{n}+\sqrt{n-1}}\geq -\sum_{n=1}^{2N} n^{-\varepsilon-1/2},
\\
(J_{\per}f_N,f_N)=2,\qquad \|\widetilde{R}_X D_X^{1/2} f_N\|^2=\sum_{n=1}^{2N}\gd_n(\sqrt{\gd_n}+\sqrt{\gd_{n+1}})^2\geq 
\sum_{n=1}^{2N}\frac{1}{n+1}=\sum_{n=2}^{2N+1}\frac{1}{n}. 
\end{gather*}
Therefore,
\[
\inf_{f\neq 0}\frac{(\widetilde{B}_{X,\gA}f,f)}{\|f\|^2}\leq \frac{\bigl((J_{\per} + \frac{1}{2}\widetilde{\mathcal{A}}_\gA)f_N,f_N\bigr)}{\|\widetilde{R}_X D_X^{1/2} f_N\|^2}\leq -\frac{\sum_{n=1}^{2N} n^{-\varepsilon-1/2}}{\sum_{n=2}^{2N+1}n^{-1}}\asymp 
-\frac{(2N)^{1/2-\varepsilon}}{\log (2N+1)},\qquad N\to \infty.
\]
Since $\varepsilon\in[0,1/2)$, the operator $\widetilde{B}_{X,\gA}$ is not lower semibounded and hence so is $B_{X,\gA}$. By Theorem \ref{th5.22}, $\rH_{X,\gA}$ is not lower semibounded too.
\end{proof}
\begin{remark}\label{rem_IV.2.5_02}
The matrix $B_{X,\gA}$ in Proposition \ref{prop_IV.2.5_01} can be considered as an unbounded Jacobi matrix with periodically modulated entries \cite{JanNab01,JanNab03}. But in the above situation we cannot apply the criteria of Janas and Naboko \cite[\S 2]{JanNab03} since $\sigma_{ac}(J_{\per})=[0,2]$.
In the proof of Proposition \ref{prop_IV.2.5_01} we follow the
line of \cite[Example 3.2]{JanNab03}.
\end{remark}
\begin{remark}\label{rem_IV.2.5_03}
$(i)$ In \cite[Theorem 3.2]{Mih_93},
it was announced (without proof) that 
$\rH_{X,\gA}$ is lower semibounded if $\cI=\R_+$ and \eqref{sb_brasche} holds. However, by Proposition \ref{prop_IV.2.5_01}, $\rH_{X,\gA}$ may be not lower semibounded even in the case $\lim_{n\to\infty}\gA_n=0$.

$(ii)$ Using the form method, semiboundedness of the the operator
$\rH_{X,\gA}$ has been studied by Brasche (see \cite{Bra85} and
references therein). In the case when all strength $\gA_n$ are
negative, he  obtained a criterion for the operator $\rH_{X,\gA}$
to be lower semibounded \cite[Theorem 3]{Bra85}. Note also that
Proposition \ref{prop_IV.2.5_01} can be extracted from
\cite[Theorem 3]{Bra85}.

Semiboundedness and discreteness  of the operator
$\rH_{X,\gA}$ will be treated by  using the form method in our
forthcoming paper.
\end{remark}

\section{Operators with $\delta'$-interactions}\label{sec_delta'}

Let $\cI$ and $X$ 
be as in Section \ref{Subsec_IV.1_bt_interactions} and let  $\gB=\{\gB_n\}_{n=1}^\infty\subset \R$. Consider the following operator in $L^2(\cI)$
\begin{gather}
\rH^0_{X,\gB}:=-\frac{\rD^2}{\rD x^2},\qquad \notag \\
\dom(\rH^0_{X,\gA})=\{f\in W^{2,2}_{\comp}(\cI\setminus X): \begin{array}{c}
         f'(0+)=0,\ f'(x_n+)=f'(x_n-)\\
         f(x_n+)-f(x_n-)=\gB_n f'(x_n)
         \end{array},\ x_n\in X
 \}.\label{delta'}
\end{gather}
Note that $\rH_{X,\gB}^0$ is symmetric in $L^2(\cI)$. Denote its closure by $\rH_{X,\gB}$, $\rH_{X,\gB}=\overline{\rH^0_{X,\gB}}$.  The Hamiltonian $\rH_{X,\gB}$ is known in the literature as the Hamiltonian of $\delta'$-interactions with strengths $\gB_n$ at points $x_n$ (see \cite{Alb_Ges_88, Alb_Kur_00, Gol_Man_09, Ges_Hol_87, S86}) and it is associated with the formal differential expression
\begin{equation}\label{IV.3.0_01}
\ell_{X,\gB}:=-\frac{\rD^2}{\rD x^2}+\sum_{n=1}^\infty \gB_n(\cdot,\delta_n')\delta_n', \qquad \gB_n\in\R,
\end{equation}
where $\delta_n':=\delta'(x-x_n)$.

In what follows we always assume that $\beta_n\neq 0,\ n\in\N$, and $\gd^*<\infty$.

\subsection{Parametrization of the operator $\rH_{X,\gB}$}\label{sec_delta'_B}

Following the line of reasoning of Subsection
\ref{sss_IV.2.1_boun_op}, we treat  $\rH_{X,\gB}$ as
an extension of $\rH_{\min}$ defined by
\eqref{I_06}. As in Subsection
\ref{sss_IV.2.1_boun_op} we consider two parameterizations of
$\rH_{X,\gB}$ corresponding to the boundary triplets constructed
in Theorems  \ref{th_bt_2} and \ref{th_bt_1}.

\textbf{1. The first parametrization.} We begin with the triplet
$\Pi^1=\{\cH,\gG_0^1,\gG_1^1\}$ constructed in Theorem
\ref{th_bt_2}  and denote  by $\Theta_1$  the
linear relation parameterizing the operator $\rH_{X,\gB}$ in the
triplet $\Pi_1$ according to \eqref{IV.1.1_12'}. Since
$\beta_n\neq 0$, $n\in\N$, the operator $\rH_{X,\gB}$ is disjoint
with the operator $\rH_0:=\rH_{\min}^*\lceil\ker(\Gamma_0^1)$
(cf. \eqref{h_0} and \eqref{IV.1.1_06}).
Therefore, by Proposition \ref{prop_II.1.2_01},
the linear relation $\Theta_1$ 
is a closed (not necessarily densely
defined) operator.

Consider the following 
Jacobi matrix
\begin{equation}\label{IV.3.1_02}
B_{X,\gB}:=\left(\begin{array}{cccccc}
\gd_1^{-2} & \gd_1^{-2} & 0& 0& 0 & \dots\\
\gd_1^{-2} & \frac{\gd_1^{-1}}{\beta_1}+\gd_1^{-2} & \frac{\gd_1^{-1/2}\gd_2^{-1/2}}{\beta_1}& 0& 0 & \dots\\
0 & \frac{\gd_1^{-1/2}\gd_2^{-1/2}}{\beta_1} & \frac{\gd_2^{-1}}{\beta_1}+\gd_2^{-2}& \gd_2^{-2}& 0 & \dots\\
0 & 0 & \gd_2^{-2}& \frac{\gd_2^{-1}}{\beta_2}+\gd_2^{-2}& \frac{\gd_2^{-1/2}\gd_3^{-1/2}}{\beta_2} & \dots\\
0 & 0 & 0& \frac{\gd_2^{-1/2}\gd_3^{-1/2}}{\beta_2}& \frac{\gd_3^{-1}}{\beta_2}+\gd_3^{-2} & \dots\\
\dots & \dots & \dots& \dots& \dots & \dots
\end{array}\right).
\end{equation}
Note that $B_{X,\gB}$ admits the representation
\begin{equation}\label{IV.3.1_02B}
B_{X,\gB}=R_X^{-1}(\widetilde{B}_{\gB}-Q_X)R_X^{-1},\qquad
\widetilde{B}_{\gB}=\left(\begin{array}{ccccc}
0 & 0 & 0& 0&  \dots\\
0 & \frac{1}{\beta_1} & \frac{1}{\beta_1}& 0&  \dots\\
0 & \frac{1}{\beta_1} & \frac{1}{\beta_1}& 0&  \dots\\
0 & 0 & 0& \frac{1}{\beta_2}& \dots\\
\dots & \dots & \dots&  \dots & \dots
\end{array}\right),
\end{equation}
where $R_X=\oplus_{n_1}^\infty R_n, \
Q_X=\oplus_{n=1}^\infty Q_n$ are determined by \eqref{IV.1.1_09B}.
Arguing as in the proof of Proposition \ref{prop_IV.2.1_01}, we arrive at the following proposition. 
\begin{proposition}\label{prop_IV.3.1_01}
Let $\Pi^1=\{\cH,\Gamma_0^1,\Gamma_1^1\}$ be the boundary triplet
constructed in Theorem \ref{th_bt_2} and let $B_{X,\gB}$ be the
minimal closed symmetric operator associated with the matrix
\eqref{IV.3.1_02}. Then $\Theta_1$ is densely
defined, $\Theta_1\in \cC(\cH),$ and  $\Theta_1=B_{X,\gB}$, that
is
\begin{equation}\label{IV.3.1_03}
\rH_{X,\gB}=\rH_{B_{X,\gB}}:=\rH_{\min}^*\lceil\dom(\rH_{B_{X,\gB}}),\qquad \dom\rH_{B_{X,\gB}}:=\{f\in\dom(\rH_{\min}^*):\ \Gamma_1^1=B_{X,\gB}\Gamma_0^1\}.
\end{equation}
\end{proposition}

\textbf{2. The second parametrization.}
Consider now the boundary triplet $\Pi^2=\{\cH,\Gamma_0^2,\Gamma_1^2\}$ constructed in Theorem \ref{th_bt_1}.
 Further, consider another  Jacobi matrix
\begin{equation}\label{IV.3.1_04}
B_{X,\gB}=\left(\begin{array}{cccccc}
0 & -\gd_1^{-2} & 0& 0& 0 & \dots\\
-\gd_1^{-2} & -(\beta_1+\gd_1)\gd_1^{-3} & \gd_1^{-3/2}\gd_2^{-1/2}& 0& 0 & \dots\\
0 & \gd_1^{-3/2}\gd_2^{-1/2} & 0& -\gd_2^{-2}& 0 & \dots\\
0 & 0 & -\gd_2^{-2}& -(\beta_2+\gd_2)\gd_2^{-3}& \gd_2^{-3/2}\gd_3^{-1/2} & \dots\\
0 & 0 & 0& \gd_2^{-3/2}\gd_3^{-1/2} & 0 & \dots\\
\dots & \dots & \dots& \dots& \dots & \dots
\end{array}\right).
\end{equation}
Though we denote by  $B_{X,\gB}$ two different
Jacobi matrices \eqref{IV.3.1_02B} and \eqref{IV.3.1_04}, it will
not lead to misunderstanding in the sequel. Using the boundary
triplet $\Pi^2=\{\cH,\Gamma_0^2,\Gamma_1^2\}$,
after straightforward calculations we arrive at the following
parametrization of $\rH_{X,\gB}$.
\begin{proposition}\label{prop_IV.3.1_02}
Let $\Pi^2=\{\cH,\Gamma_0^2,\Gamma_1^2\}$ be the boundary triplet
constructed in Theorem \ref{th_bt_1} and let
$B_{X,\gB}$ be the minimal closed symmetric operator associated
with the Jacobi matrix \eqref{IV.3.1_04}. Then
\begin{equation}\label{IV.3.1_05}
\rH_{X,\gB}=\rH_{B_{X,\gB}}:=\rH_{\min}^*\lceil\dom\rH_{B_{X,\gB}},\qquad \dom\rH_{B_{X,\gB}}=\{f\in\dom(\rH_{\min}^*):\
\Gamma_1^2=B_{X,\gB}\Gamma_0^2\}.
\end{equation}
\end{proposition}

\subsection{Self--adjointness}\label{sec_delta'_sa}

The following result gives a self-adjointness criterion for the
operator with $\delta'$-interactions on $X$.
     \begin{theorem}\label{th_delta'_sa}
 The operator $\rH_{X,\gB}$ has equal deficiency indices and $\mathrm{n}_+(\rH_{X,\gB}) = \mathrm{n}_-(\rH_{X,\gB})\leq1$.
Moreover, $\rH_{X,\gB}$ is self-adjoint if and only if at least
one of the following conditions is satisfied:
\item $ (i)$\quad $\sum_{n=1}^{\infty} \gd_n =\infty$, i.e.,
$\cI=\R_+$. \item $(ii)$\quad $\sum_{n=1}^\infty \left[\gd_{n+1}
\left|\sum_{i=1}^n (\beta_i+\gd_i)\right|^2\right]=\infty$.
\end{theorem}
\begin{proof}
Combining Theorem \ref{th_IV.1.1_02} $(i)$ with Proposition
\ref{prop_IV.3.1_01}, we get
$\mathrm{n}_\pm(\rH_{X,\gB})=\mathrm{n}_\pm(B_{X,\gB})$.
Since $B_{X,\gB}$ is a minimal Jacobi operator,
$\mathrm{n}_+(\rH_{X,\gB}) = \mathrm{n}_-(\rH_{X,\gB})\leq1$.

Further, consider the Jacobi matrix $B_{X,\gB}$
defined by \eqref{IV.3.1_02}. One can check that
$B_{X,\gB}$ admits the representation
\eqref{A_01}. Namely,
\begin{equation}\label{IV.3.1_02C}
B_{X,\gB}=R_X^{-1}(I+U)D_{X,\gB}^{-1}(I+U^*)R_X^{-1},\quad
D_{X,\gB}:=\left(\begin{array}{ccccc}
\gd_1 & 0 & 0& 0&  \dots\\
0 & \beta_1 & 0& 0& \dots\\
0 & 0 & \gd_2&  0 & \dots\\
0 & 0 & 0& \beta_2&  \dots\\
\dots &  \dots& \dots& \dots & \dots
\end{array}\right),
\end{equation}
where $U$ is unilateral shift in
$l_2$ and $R_X=\oplus_{n=1}^\infty R_n$ is
defined by \eqref{IV.1.1_09B}. In other words, $B_{X,\gB}$
coincides with $J_{m,l}$ defined by \eqref{A_01} if we set
\begin{equation}\label{IV.3.1_02D}
l_{2n-1}:=\gd_n,\quad l_{2n}:=\gB_n,\qquad m_{2n-1}=m_{2n}:=\gd_n,\qquad n\in\N.
\end{equation}
Therefore, the corresponding difference equation $\tau_{X,\gB}y=0$ has the following linearly independent solutions (cf. \cite[formulas (0.9), p.236]{Akh})
\[
\begin{array}{c}
P(0):=\{p_n\}_{n=1}^\infty,\qquad p_{2n-1}=-p_{2n}=\sqrt{\gd_n}\\
Q(0):=\{q_n\}_{n=1}^\infty,\qquad q_{2n-1}=-\sqrt{\gd_n}\sum_{k=1}^{n-1} (\gB_k+\gd_k),\quad q_{2n}=-q_{2n-1}+\gd_n^{3/2}
\end{array}
,\qquad n\in\N.
\]
The operator $B_{X,\gB}$ is symmetric with $\mathrm{n}_\pm(B_{X,\gB})=1$ precisely when $P(0),\ Q(0)\in l_2$ (cf. \cite{Akh, Ber68}). The latter holds if and only if both conditions $(i)$ and $(ii)$ are not satisfied.
\end{proof}
Condition $(i)$ of Theorem \ref{th_delta'_sa} immediately yields the following result of Buschmann, Stolz and Weidmann \cite[Theorem 4.7]{bsw}. 
\begin{corollary}[\cite{bsw}]\label{cor_bsw}
If \ $\cI=\R_+$, then the operator $\rH_{X,\gB}$ with $\delta'$-interactions is self-adjoint.
\end{corollary}

\begin{remark}\label{rem_IV.3.2_01}
In the case $d_*=0$, the structure of the boundary matrices
$B_{X,\gA}$ and $B_{X,\gB}$ that correspond to operators with
$\delta$- and $\delta'$-interactions, respectively, is completely
different. Therefore,  the spectral properties of the operators $\rH_{X,\gA}$ and $\rH_{X,\gB}$
are substantially different (cf. Proposition \ref{cor_delta_ber}
and Corollary \ref{cor_bsw}). Moreover, for the Hamiltonian
$\rH_{X,\gB}$ Theorem \ref{th_delta'_sa}
gives simple self-adjointness criterion  formulated in terms of
both $X$ and $\gB$, although for the Hamiltonian  $\rH_{X,\gA}$ we
have only necessary and sufficient conditions.
     \end{remark}

\subsection{Resolvent comparability}\label{sec_delta'_rc}

Let us fix $X\subset\cI$ and assume that
$d^*<\infty$. Consider the Hamiltonians $\rH_{X,\gB^{(1)}}$ and
$\rH_{X,\gB^{(2)}}$ \eqref{delta'} with strengths $\gB=\gB^{(1)}$
and $\gB=\gB^{(2)}$, respectively.

\begin{proposition}\label{th_delta'_res}
Suppose $\rH_{X,\gB^{(1)}}$ and
$\rH_{X,\gB^{(2)}}$  are self-adjoint. Let also $B_{X,\gB^{(1)}}$ and $B_{X,\gB^{(2)}}$ be the
corresponding (self-adjoint) Jacobi operators defined either by
\eqref{IV.3.1_02} or by \eqref{IV.3.1_04}.
Then:
\item $(i)$ \quad For any  $p\in(0,\infty]$ and for any $z\in \rho(\rH_{X,\gB^{(1)}})\cap\rho(\rH_{X,\gB^{(2)}})$ the inclusion
\begin{equation}\label{IV.3.3_01}
(\rH_{X,\gB^{(1)}}-z)^{-1}-(\rH_{X,\gB^{(2)}}-z)^{-1}\in \mathfrak{S}_p
\end{equation}
 is equivalent to the inclusion
\begin{equation}\label{IV.3.3_02}
(B_{X,\gB^{(1)}}-\I)^{-1}-(B_{X,\gB^{(2)}}-\I)^{-1}\in \mathfrak{S}_p.
\end{equation}
\item $(ii)$ \quad If
\[
    \left\{\left(\frac{1}{\gB_n^{(1)}}-\frac{1}{\gB_n^{(2)}}\right)\left(\frac{1}{\gd_{n}}+\frac{1}{\gd_{n+1}}\right)\right\}_{n=1}^\infty \in l_p, \quad p\in (0,\infty)\qquad \ (\in c_0, \ p=\infty),
 \]
     then \eqref{IV.3.3_01} holds.
\item $(iii)$ \quad If
\[
\left\{\frac{\gB_n^{(1)}-\gB_n^{(2)}}{\gd_{n}^{3}}\right\}_{n=1}^\infty \in l_p, \quad p\in (0,\infty)\qquad \ (\in c_0, \ p=\infty),
\]
 then \eqref{IV.3.3_01} holds.
\end{proposition}
\begin{proof}
$(i)$ follows from Theorem \ref{th_IV.1.1_02} and Propositions \ref{prop_IV.3.1_01} and \ref{prop_IV.3.1_02}.

Proof of $(ii)$ and $(iii)$ is similar to the proof of Corollary \ref{col_rc_2}.
We  only emphasize  that for proving  $(ii)$ we use parametrization \eqref{IV.3.1_02},  while
for proving  $(iii)$ we exploit  parametrization \eqref{IV.3.1_04}
of the Hamiltonians $\rH_{X,\gB^{(1)}}$ and $\rH_{X,\gB^{(2)}}$.
\end{proof}
In the case $d_*>0$, the resolvent comparability criterion was obtained in \cite{Mih_96b}.
\begin{corollary}[\cite{Mih_96b}]\label{cor_rc_mih}
If $0<\gd_*\leq\gd^*<\infty$, then \eqref{IV.3.3_01} is equivalent to the  inclusion
\begin{equation}\label{IV.3.3_03}
(\gB_n^{(1)}-\I)^{-1}-(\gB_n^{(2)}-\I)^{-1}\in l_p, \qquad
p\in(0,\infty),\qquad (\in c_0,\quad p=\infty).
\end{equation}
\end{corollary}
The proof of Corollary \ref{cor_rc_mih} can be extracted from
Proposition  \ref{th_delta'_res} $(i)$ and we omit it.

%

\subsection{Operators with discrete spectrum}\label{sec_delta'_disc}

Following the line of Subsection \ref{sec_delta_disc}, we begin with the criterion for the operator $\rH_{X,\gB}$ to have purely discrete spectrum.
          \begin{theorem}\label{th_disc_d'}
 Let $B_{X,\gB}$ be the minimal Jacobi operator defined either by
\eqref{IV.3.1_02} or by \eqref{IV.3.1_04}.
\item $(i)$  If $\mathrm{n}_\pm(B_{X,\gA})=1$, i.e., both conditions of Theorem
\ref{th_delta'_sa} are not satisfied, then any self-adjoint
extension of $\rH_{X,\gB}$ has discrete spectrum.
\item $(ii)$ If $B_{X,\gB} = B_{X,\gB}^*$, then the Hamiltonian
$\rH_{X,\gB} (= \rH_{X,\gB}^*)$ has discrete spectrum if and only
if
\begin{description}
 \item $\bullet$\quad $\lim_{n\to\infty}\gd_n=0$, and
 \item $\bullet$\quad  $B_{X,\gB}$ has discrete spectrum.
\end{description}
      \end{theorem}
      \begin{proof}
      Easily follows from Theorem \ref{th_IV.1.1_02} and the results of Subsection \ref{sec_delta'_B}.
\end{proof}
Let us first present several simple necessary conditions for the operator $\rH_{X,\gB}$ to have purely discrete spectrum.
\begin{proposition}\label{prop_d'_disc1}
Let $\cI=\R_+$, $\gd_n\to 0$. 
If there exists a positive constant $C> 0$ such that at least one of the following conditions is satisfied:
\item $(i)$ \quad $\gB_n\geq -C\gd_n^3$, \quad $n\in\N$,
\item $(ii)$ \quad $\gB_n^{-}\le -C(\gd_n^{-1}+\gd_{n+1}^{-1})$, \quad $n\in\N$, \qquad ($\gB_n^{-}:=\gB_n$ if $\gB_n<0$ and $\gB_n^{-}:=-\infty$ if $\gB_n>0$),\\
then the spectrum of the operator $\rH_{X,\gB}$ is not discrete.
\end{proposition}
\begin{proof}
First, assume that $\gB_n>0$, $n\in\N$. Consider the matrix \eqref{IV.3.1_02}. Since $B_{X,\gB}$ admits the representation \eqref{IV.3.1_02C}, we can apply the discreteness criterion of Kac and Krein (Theorem \ref{th_Append_01}). However, by \eqref{IV.3.1_02D}, neither $\{m_n\}_{n=1}^\infty$ nor $\{l_n\}_{n=1}^\infty$ is in $l_1$ if $\{\gd_n\}_{n=1}^\infty\notin l_1$. Hence, by Remark \ref{rem_A01}, the spectrum of $B_{X,\gB}$ is not discrete. Applying Theorem \ref{th_disc_d'}, we conclude that the spectrum of $\rH_{X,\gB}$ is not discrete. 

Consider now the matrix $B_{X,\gB}$ defined by
\eqref{IV.3.1_04} and assume that condition $(i)$ is satisfied,
i.e., $\gB_n\ge -C\gd_n^3$, $n\in \N$, with some positive constant
$C>0$. Setting $\wt{\gB_n}:=\gB_n$ if $\gB_n>0$ and
$\wt{\gB_n}:=C\gd_n^3$ if $\gB_n<0$, we obtain
$\{(\gB_n-\wt{\gB}_n)\gd^{-3}_n\}_{n=1}^\infty\in l_\infty$ and,
by Proposition \ref{th_delta'_res}$(iii)$,
$B_{X,\gB}$ is a bounded perturbation of $B_{X,\wt{\gB}}$.
Therefore, the spectra of $B_{X,\gB}$ and
$B_{X,\wt{\gB}}$ are discrete only simultaneously.
However, as it is already proved, the spectrum of
$B_{X,\wt{\gB}}$ is not discrete since $\wt{\gB_n}>0$, $n\in\N$.
Theorem \ref{th_disc_d'} $(ii)$ completes the proof.

Assume now that condition $(ii)$ holds. Then the matrix $B_{X,\gB}$ of the form \eqref{IV.3.1_02} is a bounded perturbation of the matrix $B_{X,|\gB|}$, where $|\gB|:=\{|\gB_n|\}_{n=1}^\infty$, since
\[
\left\{\left(\frac{1}{\gB_n}-\frac{1}{|\gB_n|}\right)\left(\frac{1}{\gd_{n}}+\frac{1}{\gd_{n+1}}\right) \right\}_{n=1}^\infty \in l_\infty.
\]
Therefore, $(i)$ implies that the spectrum of $B_{X,\gB}$ is not discrete and hence the spectrum of  $\rH_{X,\gB}$ is not discrete.
\end{proof}
\begin{corollary}\label{cor_d'_disc1}
If $\cI=\R_+$ and $\gB_n>0$ for all $n\in\N$, then the spectrum of $\rH_{X,\gB}$ is not discrete.
\end{corollary}

The following result gives sufficient condition for the operator $\rH_{X,\gB}$ to have discrete spectrum.
\begin{proposition}\label{prop_d'_disc2}
Assume $\gB_n+\gd_n\ge0$ for all $n\in\N$.
\item $(i)$ Let $\cI=[0,b)$ be a bounded interval and let $X$ and $\gB$ be such that the Hamiltonian $\rH_{X,\gB}$ is self-adjoint. Then $\rH_{X,\gB}$ has discrete spectrum if  and only if
\begin{equation}\label{IV.3.4_08}
\lim_{n\to\infty}(b-x_n)\sum_{j=1}^n(\gB_j+\gd_j)=0.
\end{equation}
\item
$(ii)$ Let $\cI=\R_+$. Then the Hamiltonian $\rH_{X,\gB}\ (=\rH_{X,\gB}^*)$ has discrete spectrum if and only if 
\begin{equation}\label{IV.3.4_09}
\lim_{n\to\infty}x_n\sum_{j=n}^\infty\gd_j^3=0\quad
\text{and}\quad
\lim_{n\to\infty}x_n\sum_{j=n}^\infty(\gB_j+\gd_j)=0.
\end{equation}
\end{proposition}
\begin{proof} 
Consider the minimal symmetric operator associated with the Jacobi matrix \eqref{IV.3.1_04}. First note that it is unitarily equivalent to the Jacobi operator with positive offdiagonal entries,
\begin{equation}\label{IV.3.4_01A}
B_{X,\gB}'=\left(\begin{array}{cccccc}
0 & \gd_1^{-2} & 0& 0& 0 & \dots\\
\gd_1^{-2} & -(\beta_1+\gd_1)\gd_1^{-3} & \gd_1^{-3/2}\gd_2^{-1/2}& 0& 0 & \dots\\
0 & \gd_1^{-3/2}\gd_2^{-1/2} & 0& \gd_2^{-2}& 0 & \dots\\
0 & 0 & \gd_2^{-2}& -(\beta_2+\gd_2)\gd_2^{-3}& \gd_2^{-3/2}\gd_3^{-1/2} & \dots\\
0 & 0 & 0& \gd_2^{-3/2}\gd_3^{-1/2} & 0 & \dots\\
\dots & \dots & \dots& \dots& \dots & \dots
\end{array}\right).
\end{equation}
Further, consider the orthogonal decomposition
\[
l_2=\cH_1\oplus\cH_2, \qquad \cH_1=\Span\{e_{2n-1}\}_{n\in\N},\quad \cH_2=\Span\{e_{2n}\}_{n\in\N}.
\]
Define the unitary operators
\begin{equation}\label{VI_Vop}
V_j:\cH_j\to l_2,\quad (j=1,2),\quad  V_1(\mathrm{e}_{2n-1})=\mathrm{e}_n\quad \text{and} \quad V_2(\mathrm{e}_{2n})=\mathrm{e}_n,\quad n\in\N.
\end{equation}
Then the operator $\widetilde{B}_{X,\gB}:=VB'_{X,\gB}V^{-1}$ with $V:=V_1\oplus V_2$ admits the representation
\[
\widetilde{B}_{X,\gB}=\left(\begin{array}{cc}
D_X^{-1/2}&0\\
0&D_X^{-3/2}
\end{array}\right)\left(\begin{array}{cc}
0_{\cH_1}& I+U\\
I+U^*& -(\mathcal{B}_\gB+D_X)
\end{array}\right)\left(\begin{array}{cc}
D_X^{-1/2}&0\\
0&D_X^{-3/2}
\end{array}\right),
\]
where
\[
\mathcal{B}_\gB=\diag(\beta_n) ,\qquad D_X=\diag(\gd_n), 
\]
and $U$ is unilateral shift.
Since $B_{X,\gB}'$ is symmetric and $\dim \ker B'_{X,\gB}\leq 1$, the inverse operator $(\widetilde{B}_{X,\gB})^{-1}$ is closed on $\cH\ominus\ker(B_{X,\gB}')$ and is given by the following matrix 
\[
(\widetilde{B}_{X,\gB})^{-1}=\left(\begin{array}{cc}
D_X^{1/2}&0\\
0&D_X^{3/2}
\end{array}\right)\left(\begin{array}{cc}
-(I+U^*)^{-1}(\mathcal{B}_\gB+D_X)(I+U)^{-1}&(I+U^*)^{-1}\\
(I+U)^{-1}&0
\end{array}\right)\left(\begin{array}{cc}
D_X^{1/2}&0\\
0&D_X^{3/2}
\end{array}\right).
\]
Therefore, the operator $(B_{X,\gB})^{-1}$ is compact precisely when
the spectra of operators
\begin{gather}\label{IV.3.4_01B'}
J_\gB:=D_X^{-1/2}(I+U)(\mathcal{B}_\gB+D_X)^{-1}(I+U^*)D_X^{-1/2},\\
J_X:=D_X^{-1/2}(I+U)D_X^{-3}(I+U^*)D_X^{-1/2}, \label{IV.3.4_01}
\end{gather}
are purely discrete. Without loss of generality we can assume that
$\gB_n+\gd_n> 0$ for all $n\in\N$. Indeed, in the opposite case we
can choose $\widetilde{\gB}_n$ satisfying the assumption of
Proposition \ref{prop_d'_disc2} and such that
$\widetilde{\gB}_n+\gd_n> 0$, $n\in\N$, and
$\{(\widetilde{\gB}_n-\gB_n)\gd_n^{-3}\}_{n=1}^\infty\in c_0$.
By Proposition \ref{th_delta'_res}$(iii)$,
$B_{X,\widetilde{\gB}}$ is a bounded perturbation of $B_{X,\gB}$
and hence the operators $\rH_{X,\gB}$ and
$\rH_{X,\widetilde{\gB}}$ have discrete spectrum simultaneously.

As in Subsection \ref{ss_II_krein},
with  $J_X$ and $J_\gB$ we associate the functions
\begin{gather}\label{IV.3.4_07}
\mathcal{M}_X(x)=\sum_{y_{n-1}<x}\gd_n,\quad y_n-y_{n-1}=\gd_n^3,\qquad
\mathcal{M}_\gB(x)=\sum_{z_{n-1}<x}\gd_n,
\quad
z_n-z_{n-1}=\gB_n+\gd_n,
\end{gather}
respectively. Here $x>0$ and $y_0=z_0=0$.

We begin with the case of a finite interval $\cI$, i.e.,
assume that $\sum_{n\in\N}\gd_n<\infty$. Then
$\sum_{n\in\N}\gd_n^3<\infty$ and hence the string with the mass
$\mathcal{M}_X$ is regular. Therefore,  $\sigma(J_X)$ is discrete
(see \cite[Section 11.8]{KK71}). Moreover, by
Theorem \ref{th_Append_01}, the operator $J_\gB$ has discrete
spectrum precisely when \eqref{IV.3.4_08} holds.

Assume now that $\cI=\R_+$, i.e., $\sum_{n\in\N}\gd_n=\infty$. By Theorem \ref{th_Append_01}, $\sigma(J_X)$ is discrete
if and only if $\{\gd_n^3\}_{n=1}^\infty\in l_1$ and the first condition in \eqref{IV.3.4_09} holds.
Further, $\sigma(J_\gB)$ is discrete precisely when  $\{\gB_n+\gd_n\}_{n=1}^\infty\in l_1$ and the
function $\mathcal{M}_\gB$ also satisfies the second condition in \eqref{A_05}, that is the second condition in \eqref{IV.3.4_09} holds.

Theorem \ref{th_disc_d'} completes the proof.
\end{proof}
\begin{corollary}\label{cor_d'_disc2}
Let $\cI=\R_+$. Then for any $\gB$ the spectrum of the operator $\rH_{X,\gB}$ is not discrete if at least one of the following conditions is satisfied
\item $(i)$ $\{\gd_n\}_{n=1}^\infty\notin l_3$,
\item $(ii)$ $\{\gd_n\}_{n=1}^\infty\in l_3$ and
\begin{equation}\label{6.19}
\lim_{n\to\infty}x_n\sum_{j=n}^\infty \gd_j^3>0.
\end{equation}
\end{corollary}
\begin{proof}
Let $\sigma(\rH_{X,\gB})$ be discrete. Consider the operator $J_X$ defined by \eqref{IV.3.4_01}. It easily follows from the proof of Proposition \ref{prop_d'_disc2} that $\sigma(J_X)$ is discrete.
However, by Theorem \ref{th_Append_01}, $J_X$ has discrete spectrum if and only if $\{\gd_n\}_{n=1}^\infty\in l_3$ and the limit in \eqref{6.19} equals $0$.
\end{proof}
Let us illustrate the above results by the following example.
\begin{example}\label{example_IV.3.4_01}
Let $\cI=\R_+$. Consider the Hamiltonian
\[
\rH_\gB=-\frac{\rD^2}{\rD x^2}+\sum_{n=1}^\infty\gB_n(\cdot,\delta'(x-n^\varepsilon))\delta'(x-n^\varepsilon),\qquad 0<\varepsilon<1.
\]
First note that, by Theorem \ref{th_delta'_sa} (see also \cite[Theorem 4.7]{bsw}), the operator $\rH_\gB$ is self-adjoint for any $\gB=\{\gB_n\}_{n=1}^\infty\subset\R$. Since $x_n=n^\varepsilon$, we get $\gd_n\asymp n^{\varepsilon-1}$ and $\sum_{j=1}^n\gd_j^3\asymp n^{3\varepsilon-2}$. Therefore, the following is true:
\item $(i)$ If $\varepsilon \geq 1/2$, then for any $\gB$ the spectrum of $\rH_{\gB}$ is not discrete.
\item $(ii)$ If $\varepsilon <1/2$ and either $\gB_n^-\geq -C n^{3\varepsilon-3},\ n\in\N$ or $\gB_n^-\leq -C n^{1-\varepsilon},\ n\in\N$, with some positive constant $C>0$, then the spectrum of $\rH_{\gB}$ is not discrete.
\item $(iii)$ Assume $\varepsilon<1/2$ and $\gB_n+\gd_n=\gB_n+n^\varepsilon-(n-1)^\varepsilon\geq 0$, $n\in\N$. Then
 the operator $\rH_\gB$ has discrete spectrum if and only if
\[
 \lim_{n\to\infty}n^{\varepsilon}\sum_{j=n}^\infty(\gB_j+j^\varepsilon-(j-1)^\varepsilon)=0.
\]
\end{example}

\subsection{Semiboundedness}\label{sec_delta'_sembd}

Combining Theorem \ref{th_IV.1.1_02} $(iii)$ with Proposition \ref{prop_IV.3.1_01}, we arrive at the following result.
\begin{theorem}\label{th_deltapr_sembd}
The operator $\rH_{X,\gB}$ with $\delta'$-interactions on $X$ is lower semibounded if and only if
the Jacobi operator $B_{X,\gB}$ of the form \eqref{IV.3.1_02} is lower semibounded.
\end{theorem}
\begin{proposition}\label{cor_IV.3.4_01}
For the operator $\rH_{X,\gB}$ to be lower semibounded it is necessary that
\begin{equation}\label{IV.3.5_01A}
\frac{1}{\gB_n}\geq -C_1 \gd_n- \frac{1}{\gd_{n}},\quad\mathrm{and}\quad
\frac{1}{\gB_n}\geq -C_1 \gd_{n+1} - \frac{1}{\gd_{n+1}},\qquad n\in\N,
\end{equation}
and it is sufficient that
\begin{equation}\label{IV.3.5_01B}
\frac{1}{\gB_n}\geq -C_2\min\{\gd_{n},\gd_{n+1}\},\qquad n\in\N,
\end{equation}
with some positive constants $C_1,\ C_2>0$ independent of $n\in\N$. 
\end{proposition}
\begin{proof}
By Theorem \ref{th_deltapr_sembd}, $\rH_{X,\gB}$
is lower semibounded if and only if the matrix (\ref{IV.3.1_02})
is lower semibounded. First, consider the representation
(\ref{IV.3.1_02C}). Let $V_1$ and $V_2$ be the unitary mappings
defined by (\ref{VI_Vop}) and $V := V_1\oplus V_2$. Then it is easy to check that
\[
VR_{X}V^{-1}=\left(\begin{array}{cc}
D_X& 0\\
0& D_X
\end{array}\right), \quad V(I+U)V^{-1}=\left(\begin{array}{cc}
                                          I  & U\\
                                          I  & I
                                      \end{array}\right),\quad VD_{X,\gB}V^{-1}=\left(\begin{array}{cc}
                                                                                   D_{X}& 0\\
                                                                                      0& \mathcal{B}_\gB
                                                                                  \end{array}\right),
\]
where $D_X:=\diag(\gd_n)$, $\mathcal{B}_\gB=\diag(\gB_n)$, $I=I_{l_2}$, and $U$ is unilateral shift in $l_2$. After straightforward calculations we obtain
\[
\widetilde{B}_{X,\gB}:=VB_{X,\gB}'V^{-1}=\left(\begin{array}{cc}
                                          D^{-2}_X+D_X^{-1/2}U\mathcal{B}^{-1}_\gB U^* D_X^{-1/2}& D_X^{-1/2}U \mathcal{B}^{-1}_\gB D_X^{-1/2} +D_X^{-2}\\
                                          D_X^{-1/2}\mathcal{B}^{-1}_\gB U^* D_X^{-1/2} +D_X^{-2}  & D_X^{-2}+\mathcal{B}^{-1}_\gB D_X^{-1}
                                      \end{array}\right),
\]
where $U_+$ is unilateral shift in $\cH_+$. Therefore, inequalities
\[
D_X^{-2}+\mathcal{B}^{-1}_\gB D_X^{-1}\geq -C_1 I, \qquad D^{-2}_X+D_X^{-1/2}U\mathcal{B}^{-1}_\gB U^* D_X^{-1/2}\geq -C_1 I
\]
are necessary for the operator $B_{X,\gB}$ to be lower semibounded. The latter is equivalent to
(\ref{IV.3.5_01A}).

To prove sufficiency we use the representation
\eqref{IV.3.1_02B} of $B_{X,\gB}$. By
\eqref{IV.1.1_09B}, $Q_X\leq 0$ and hence the
operator $B_{X,\gB}$ is lower semibounded whenever the
operator $R_X^{-1}\widetilde{B}_\gB R_X^{-1}$ is lower semibounded.  The latter is
equivalent to the validity of the following
inequalities
\[
\left(\begin{array}{cc}
\frac{1}{\beta_n}& \frac{1}{\beta_n}\\
\frac{1}{\beta_n}& \frac{1}{\beta_n}
\end{array}\right) \geq -\widetilde{C}_2 \left(\begin{array}{cc}
                                  \gd_n & 0\\
                                  0 & \gd_{n+1}
                                             \end{array}\right),\qquad n\in\N,
\]
with the constant $\widetilde{C}_2>0$ independent of $n\in\N$. Thus condition \eqref{IV.3.5_01B} is sufficient for lower semiboundedness.
The proof is completed.
\end{proof}
\begin{corollary}
Let $0<\gd_*\leq\gd^*<\infty$. Then the
Hamiltonian $\rH_{X,\gB}$ is lower semibounded if and only if
$\{\frac{1}{\gB_n}\}_{n=1}^\infty$ is lower semibounded.
\end{corollary}

\section{Operators with $\delta$-interactions and semibounded potentials}\label{Sec_V}

The results of Section \ref{Subsec_IV.2_delta} are stable under perturbations by  $L^\infty$ potentials $q$ since deficiency indices, discreteness, and lower semiboundedness are stable under bounded perturbation. In particular, the results of Section \ref{Subsec_IV.2_delta} hold true for operators
\begin{equation}\label{VII_01}
\rH_{X,\gA,q}=-\frac{\rD^2}{\rD x^2}+q(x)+\sum_{n=1}^\infty\alpha_n \delta(x-x_n),\qquad q\in L^\infty(\cI).
\end{equation}
Moreover, it follows from \cite[Theorem 3.1]{Ges_Kir_85} that self-adjointness is stable under perturbations by lower semibounded potentials  if $\gd_*>0$.

The main aim of this section is to show that in the case $\gd_*=0$ the situation is substantially different. Namely, we will show that self-adjointness of the operators with $\delta$-interactions is not stable under perturbations by positive potentials $q$ if $\gd_*=0$. 

Let $\mathcal{I}=\R_+$, $x_0=0$, $x_{n}-x_{n-1}=\gd_n:=\frac{1}{n}$, $ n\in\N$. Set
\begin{equation}\label{VII_02}
q_a(x):=a^2\sum_{n=1}^\infty n^2\chi_{(x_{n-1},x_n)}(x),\qquad a\in\R_+.
\end{equation}
Consider the operator
\begin{equation}\label{VII_03}
\rH_{X,\gA,q_a}=-\frac{\rD^2}{\rD x^2}+q_a(x)+a^2\sum_{n=1}^\infty\alpha_n \delta(x-x_n).
\end{equation}
The corresponding  minimal symmetric operator $\rH_{\min}$ has the form
   \begin{equation}\label{VII_04}
\rH_{\min}=\oplus_{n=1}^\infty\rH_n,\qquad
\rH_{n}:=-\frac{\rD^2}{\rD x^2}+a^2n^2,\quad
\dom(\rH_n)=W^{2,2}_0[x_{n-1},x_n].
     \end{equation}
In the following proposition we construct a boundary triplet for
$\rH_{\min}^*$.
\begin{proposition}\label{prop_VII_01}
For $f\in W_2^2[x_{n-1},x_n]$, define the mappings $\Gamma_j^{(n)}:W_2^2[x_{n-1},x_n]\to\C^2$,
\begin{equation}\label{VII_05}
\Gamma_0^{(n)}f:=\left(\begin{array}{c}
                \gd_n^{1/2}  f(x_{n-1}+)\\
                 -\gd_n^{1/2}  f(x_{n}-)
                       \end{array}\right),\quad \Gamma_1^{(n)}f:=\left(\begin{array}{c}
                                                                           \frac{\gd_n f'(x_{n-1}+)+(\varepsilon_1f(x_{n-1}+)-\varepsilon_2f(x_{n}-))}{\gd_n^{3/2}}\\
                                                                           \frac{\gd_n f'(x_{n}-)+(\varepsilon_1f(x_{n-1}+)-\varepsilon_2f(x_{n}-))}{\gd_n^{3/2}}
                                                                          \end{array}\right),
\end{equation}
where
\begin{equation}\label{VII_06}
\gd_n=\frac{1}{n},\qquad \varepsilon_1=\varepsilon_1(a):=a\frac{\cosh a}{\sinh a},\qquad \varepsilon_2=\varepsilon_2(a):=\frac{a}{\sinh a}.
\end{equation}
Then:
\item $(i)$\ For any $n\in \N$ the triplet $\Pi_n=\{\C^2,\Gamma_0^{(n)},\Gamma_1^{(n)}\}$ is a boundary triplet for
$\rH_{n}^*$.
\item $(ii)$ \ The direct sum $\Pi=\oplus_{n=1}^\infty\Pi_n$
is a boundary triplet for the operator $\rH_{\min}^*$.
\end{proposition}
\begin{proof}
$(i)$ Straightforward.

$(ii)$ Note that the triplet $\widetilde{\Pi}_n=\{\C^2,
\widetilde{\Gamma}_0^{(n)}, \widetilde{\Gamma}_1^{(n)}\}$ defined
by \eqref{IV.1.1_05} forms a boundary triplet for the operator
$\rH^*_{n}$ defined by  \eqref{VII_04}. The corresponding Weyl
function  $\widetilde{M}_{n}(\cdot)$ is
       \begin{equation}\label{VII_07}
 \widetilde{M}_{n}(z) = -\frac{\sqrt{z-a^2n^2}}{\sin\sqrt{z/n^2-a^2}}\left(\begin{array}{cc}
                           \cos\sqrt{z/n^2-a^2} & 1 \\
                           1 & \cos\sqrt{z/n^2-a^2}
                                                                          \end{array}\right),\qquad z\in\C_+.
        \end{equation}

It is easily seen that  $\wt\Pi := \oplus_{n=1}^\infty\wt\Pi_n$ is
not an ordinary boundary triplet for $\rH_{\min}^*$. On the other
hand, triplets $\wt\Pi_n$ and $\Pi_n$ of the form
\eqref{IV.1.1_05} and \eqref{VII_05}, respectively, are
connected by 
       \begin{equation}\label{VII_08}
\Gamma_0^{(n)}=R_n\widetilde{\Gamma}_0^{(n)},\quad \Gamma_1^{(n)}
= R_n^{-1}(\widetilde{\Gamma}_1^{(n)}-
Q_n\widetilde{\Gamma}_0^{(n)}),
     \end{equation}
where
\[
Q_n:=\widetilde{M}_n(0)=\left(\begin{array}{cc}
                           -n\varepsilon_1(a) 
                           & -n\varepsilon_2(a) \\ 
                           -n\varepsilon_2(a) 
                           & -n\varepsilon_1(a) 
                                                                          \end{array}\right)\qquad\text{and}\qquad R_n=\left(\begin{array}{cc}
                           n^{-1/2}
                           & 0\\
                           0 & n^{-1/2} 
                                                                          \end{array}\right).
\]
The corresponding Weyl functions $M_n(\cdot)$ and $\wt M_n(\cdot)$
are connected by $M_n(z)=R_n^{-1}(\widetilde{M}_n(z)-Q_n)R_n^{-1}.$
Clearly, relations  \eqref{VII_08}  coincide  with
\eqref{III.2.2_08}. Moreover, direct  calculations show that
\[
M_n(0)=0,\quad M_n'(0)=R_n^{-1}\widetilde{M}_n'(0)R_n^{-1}=a^{-2}\left(\begin{array}{cc}
                          (a-\varepsilon_1(a))(\varepsilon_1(a)-1)
                           & \varepsilon_2(a)- \varepsilon_1(a)\varepsilon_2(a)\\ 
                          \varepsilon_2(a)- \varepsilon_1(a)\varepsilon_2(a)& 
                         (a-\varepsilon_1(a))(\varepsilon_1(a)-1)
                                                                          \end{array}\right).
\]
 Therefore, by  Corollary  \ref{cor_III.2.2_02},
the direct sum $\Pi=\oplus_{n=1}^\infty\Pi_n$ forms a boundary
triplet for $\rH_{\min}^*$.   
         \end{proof}
Arguing as in Subsection \ref{sss_IV.2.1_boun_op}, we obtain that the operator $\rH_{X,\gA,q_a}$ admits the representation
\[
\rH_{X,\gA,q_a}=\rH_{\Theta}:=\rH_{\min}^*\lceil\dom(\rH_{\Theta}),\qquad \dom\rH_\Theta:=\{f\in\dom(\rH_{\min}^*):\ \{\Gamma_0, \Gamma_1\} \in \Theta\},
\]
where $\Gamma_0=\oplus_{n=1}^\infty\Gamma_0^{(n)}$ and $\Gamma_1=\oplus_{n=1}^\infty\Gamma_1^{(n)}$ are defined by \eqref{VII_05} and
the operator part $\Theta_{\op}$ of the linear relation $\Theta\in \widetilde{\mathcal{C}}(\cH)$ is unitary equivalent to the following Jacobi matrix
    \begin{gather}
B_{X,\gA,q_a}=\widetilde{R}_X^{-1}(B_{X}(a)+\mathcal{A}_\gA)\widetilde{R}_X^{-1},\qquad
B_{X}(a)=\left(\begin{array}{ccccc}
3\varepsilon_1(a) & 2\varepsilon_2(a) & 0& 0&  \dots\\
2\varepsilon_2(a) & 5\varepsilon_1(a) & 3\varepsilon_2(a)& 0&  \dots\\
0 & 3\varepsilon_2(a) & 7\varepsilon_1(a)& 4\varepsilon_2(a)&  \dots\\
0 & 0 & 4\varepsilon_2(a)& 9\varepsilon_1(a)&  \dots\\
\dots & \dots & \dots& \dots& \dots
\end{array}\right),\notag\\
\quad \text{and} \quad \widetilde{R}_X=\diag(\widetilde{r}_n),\quad \widetilde{r}_n:=\sqrt{\frac{1}{n}+\frac{1}{n+1}},\qquad \mathcal{A}_\gA=\diag(\alpha_n).\label{VII_09}
  \end{gather}
Thus we arrive at the following result.
      \begin{proposition}\label{prop_VII_02}
Let $q_a$ be defined by \eqref{VII_02} and let $B_{X,\gA,q_a}$ be the minimal
symmetric operator associated with the Jacobi matrix \eqref{VII_09}. Then the operator
$\rH_{X,\gA,q_a}$ has equal deficiency indices and
$\mathrm{n}_\pm(\rH_{X,\gA,q_a}) =
\mathrm{n}_\pm(B_{X,\gA,q_a})\leq 1$. In particular,
$\rH_{X,\gA,q_a}$ is self-adjoint if and only if so is $B_{X,\gA,q_a}$.
     \end{proposition}
Proof is straightforward and we omit it.

Let us consider $\varepsilon_1(a), \ a>0$. Since $\lim_{a\to0}\varepsilon_1(a)=1$ and $\varepsilon_1(a)\approx a$ as $a\to +\infty$, there exists $a_0>0$ such that
\begin{equation}\label{VII_10}
\varepsilon_1(a_0)=2.
\end{equation}
\begin{corollary}\label{cor_VII_01}
Let $\cI=\R_+$, $\gd_n=1/n$, and $\alpha_n=-4n-2$, $n\in\N$.
 \item $(i)$ The Hamiltonian 
 \[
\rH_{X,\gA,0}=-\frac{d^2}{dx^2}-\sum_{n=1}^\infty (4n+2) \delta(x-x_n),
\]
 is self-adjoint.
 \item $(ii)$ Let $a_0$ be defined by \eqref{VII_10}. 
 Then the Hamiltonian 
\[
\rH_{X,\gA,q_a}=-\frac{d^2}{dx^2}+a_0^2\sum_{n=1}^\infty n^2\chi_{(x_{n-1},x_n)}-\sum_{n=1}^\infty (4n+2) \delta(x-x_n),
\]
 is symmetric with $\mathrm{n}_\pm(\rH_{X,\gA,q_a})=1$.
\end{corollary}
\begin{proof}
$(i)$ follows from Example \ref{example_IV.2.2_01} $(ii)$.

$(ii)$  Consider the matrix $B_{X,\gA, q_a}$ with $a=a_0$. Clearly, $\alpha_n=-\varepsilon_1(a_0)(2n+1)$ and hence the diagonal entries of $B_{X,\gA, q_a}$ equal zero. The offdiagonal entries $b_n=n\frac{\varepsilon_2(a_0)}{\widetilde{r}_n\widetilde{r}_{n+1}}$ satisfies $b_n\approx \varepsilon_2(a_0)n^{2}/4$ and hence $\{b_n^{-1}\}_{n=1}^\infty\in l_1$. Moreover, $b_{n-1}b_{n+1}\leq b_n^2$ holds for all $n\in\N$. Therefore, Berezanskii's test \cite[Theorem VII.1.5]{Ber68} implies $\mathrm{n}_\pm(B_{X,\gA,q_a})=1$. By Proposition \ref{prop_VII_02},   $\mathrm{n}_\pm(\rH_{X,\gA,q_a})=1$.
\end{proof}

\ack{The authors thank M. Derevyagin,  L. Oridoroga, and G. Teschl
for useful discussions. We are grateful to V. Derkach, F. Gesztesy, and A. Kochubei for careful reading of the manuscript and helpful remarks.

AK gratefully acknowledges the financial support from the Junior Research
Fellowship Programme of the Erwin Schr\"odinger Institute for
Mathematical Physics and from the IRCSET Postdoctoral Fellowships Programme.}


\quad
\\
Aleksey Kostenko, \\
\emph{School of Mathematical Sciences},\\
\emph{DIT Kevin Street},
\emph{Dublin 8, IRELAND}\\
\emph{e-mail:} duzer80$@$gmail.com\\
\\
Mark Malamud, \\
\emph{Institute of Applied Mathematics and Mechanics, NAS of Ukraine,}\\
\emph{R. Luxemburg str., 74, Donetsk 83114, UKRAINE}\\
\emph{e-mail:}\ mmm$@$telenet.dn.ua

\end{document}